\documentclass{article}
\usepackage[utf8]{inputenc}

\usepackage{import}
\usepackage{mathrsfs}
\usepackage{enumerate}
\usepackage{amsthm}
\usepackage{amsmath}
\usepackage{thm-restate}
\usepackage[hidelinks]{hyperref}
\usepackage{cleveref}
\usepackage{amsfonts}
\usepackage{amssymb}
\usepackage{dsfont}
\usepackage{stmaryrd}
\usepackage{esint}
\usepackage{mathtools}
\usepackage{float}
\usepackage{accents}
\usepackage{verbatim}
\usepackage[toc]{appendix}
\usepackage{parskip}
\usepackage{relsize}
\usepackage{xspace}
\usepackage{bm}
\hypersetup{
    colorlinks,
    linkcolor={blue!80!black},
    citecolor={blue!80!black},
    urlcolor={blue!80!black}
}
\usepackage{tikz}
\usetikzlibrary{arrows.meta,positioning,matrix}
\usepackage{pgfplots}
\pgfplotsset{compat=1.18}
\usepackage{xcolor}
\usepackage{mleftright}

\usepackage{titlesec}
\titleformat*{\section}{\LARGE\bfseries}

\usepackage[top=1.0in,bottom=1.3in,left=1.15in,right=1.15in,twoside]{geometry}

\usepackage[backend=biber,url=false,isbn=false,doi=false,eprint=false]{biblatex}
\counterwithin{figure}{section}
\numberwithin{equation}{section}

\newtheorem{theorem}{Theorem}[section]
\newtheorem*{theorem*}{Theorem}
\newtheorem{lemma}[theorem]{Lemma}
\newtheorem{corollary}[theorem]{Corollary}
\newtheorem{definition}[theorem]{Definition}
\newtheorem{remark}[theorem]{Remark}

\newtheorem{claim}[theorem]{Claim}
\newtheorem{proposition}[theorem]{Proposition}

\newenvironment{claimproof}
  {%
   \begin{proof}}
  {\end{proof}}

\newcommand{\ZBC}{\ensuremath{\mathrm{\textsf{ZBC}}}\xspace}
\newcommand{\NZBC}{\ensuremath{\mathrm{\textsf{NZBC}}}\xspace}

\newcommand{\KdV}{\ensuremath{\mathrm{\textsf{KdV}}}\xspace}
\newcommand{\mKdV}{\ensuremath{\mathrm{\textsf{mKdV}}}\xspace}
\newcommand{\ZS}{\ensuremath{\mathrm{\textsf{ZS}}}\xspace}
\newcommand{\GP}{\ensuremath{\mathrm{\textsf{GP}}}\xspace}
\newcommand{\NLS}{\ensuremath{\mathrm{\textsf{NLS}}}\xspace}
\newcommand{\dNLS}{\ensuremath{\mathrm{\textsf{dNLS}}}\xspace}

\DeclareMathOperator{\R}{\mathbb{R}}

\DeclareMathOperator{\C}{\mathbb{C}}

\DeclareMathOperator{\SSS}{\mathbb{S}}
\DeclareMathOperator{\N}{\mathbb{N}}
\DeclareMathOperator{\Z}{\mathbb{Z}}

\DeclareMathOperator{\K}{\mathbb{K}}

\DeclareMathOperator{\supp}{supp}

\DeclareMathOperator{\Real}{Re}
\DeclareMathOperator{\Imag}{Im}
\DeclareMathOperator{\sign}{sign}

\DeclareMathOperator{\loc}{loc}

\DeclareMathOperator{\sech}{sech}

\DeclareMathOperator{\PV}{P.\!V.}

\NewDocumentCommand{\opnorm}{sO{}m}{%
  \IfBooleanTF{#1}{
    \left|\opnormkern\left|\opnormkern\left|
    #3
    \right|\opnormkern\right|\opnormkern\right|
  }{
    \mathopen{#2|\opnormkern #2|\opnormkern #2|}
    #3
    \mathclose{#2|\opnormkern #2|\opnormkern #2|}
  }%
}
\newcommand{\opnormkern}{\mkern-1.5mu\relax}

\DeclarePairedDelimiter{\ceil}{\lceil}{\rceil}
\DeclarePairedDelimiter{\floor}{\lfloor}{\rfloor}
\makeatletter
  \let\oldceil\ceil
  \def\ceil{\@ifstar{\oldceil}{\oldceil*}}
  \let\oldfloor\floor
  \def\floor{\@ifstar{\oldfloor}{\oldfloor*}}
\makeatother
\newcommand*\dd{\mathop{}\!\mathrm{d}}

\DeclareFontFamily{U}{mathx}{\hyphenchar\font45}
\DeclareFontShape{U}{mathx}{m}{n}{
      <5> <6> <7> <8> <9> <10>
      <10.95> <12> <14.4> <17.28> <20.74> <24.88>
      mathx10
      }{}
\DeclareSymbolFont{mathx}{U}{mathx}{m}{n}
\DeclareFontSubstitution{U}{mathx}{m}{n}
\DeclareMathAccent{\widecheck}{0}{mathx}{"71}

\newcommand{\conj}[1]{\overline{#1}}
\newcommand{\ti}[1]{\widetilde{#1}}
\newcommand{\ha}[1]{\widehat{#1}}
\newcommand{\ch}[1]{\widecheck{#1}}

\newcommand{\dg}{
  {\mathrm{d}}
}
\newcommand{\odg}{
  {\mathrm{od}}
}

\usetikzlibrary{tikzmark,arrows.meta,bending,calc}
\newcounter{coverarcarrow}
\newcommand{\overarcarrow}[2][]{\stepcounter{coverarcarrow}%
\overset{\vphantom{+}}{\tikzmarknode{coverarcarrow-\number\value{coverarcarrow}}{#2}}%
\begin{tikzpicture}[overlay,remember picture]
\draw let \p1=($(coverarcarrow-\number\value{coverarcarrow}.north east)-(coverarcarrow-\number\value{coverarcarrow}.south west)$),
 \n1={atan2(1.2ex,\x1)} in 
 [-{Stealth[bend,length={min(0.3*\x1,0.7ex)}]},thin,#1]
 ([yshift=0.2*\y1]coverarcarrow-\number\value{coverarcarrow}.north west) 
 to[out=\n1,in=180] ([yshift=0.83ex]coverarcarrow-\number\value{coverarcarrow}.north)
 to[out=0,in={180-\n1}] ([yshift=0.2*\y1]coverarcarrow-\number\value{coverarcarrow}.north east);
\end{tikzpicture}}

\let\phi=\varphi

\let\epsilon=\varepsilon

\let\del=\partial
\let\mc=\mathcal

\title{Global Well-posedness of the NLS Hierarchy with Nonzero Boundary Condition}

\author{Xian Liao and Robert Wegner}
\date{}

\begin{document}

\maketitle

\begin{abstract}
  We consider the \NLS hierarchy with the nonzero boundary condition $q(t, x) \rightarrow q_\pm \in \SSS^1$ as $x \rightarrow \pm \infty$ and prove that it is global well-posedness for initial data of high regularity. 
  Specifically, we prove well-posedness of the problem for the perturbation $p = q - q_\ast$ from a time-independent front $q_\ast$ connecting $q_-$ to $q_+$.
  
  The equations in the \NLS hierarchy are defined using a recurrence relation derived from the expansion of the logarithmic derivative of the Jost solutions associated to the Lax operator.
  Using this recurrence relation, we are able to determine explicit formulas for all terms in the \NLS hierarchy with at most one factor that is $q_x$, $\conj{q}_x$, or a derivative thereof.
  
  We then view the equation for $p$ as part of a large class of dispersive nonlinear systems, for which we develop a local well-posedness theory in weighted Sobolev spaces.
  This involves certain local smoothing and maximal function estimates, which we establish for a large class of dispersion relations with finitely many critical points.
  Finally, we globalize the solutions using the conserved energies constructed in \cite{KochLiao,KochLiao2022}.
\end{abstract}
\noindent{\sl Keywords:} NLS hierarchy, Gross-Pitaevskii equation, global well-posedness, transmission coefficient, dispersive equations, local smoothing
\\ \noindent{\sl AMS Subject Classification (2020):} 35G55, 35Q51, 35Q55, 37K10
 
\tableofcontents

\mleftright

\section{Introduction}
Consider the defocusing nonlinear Schrödinger equation
\begin{align*} \tag{\NLS} \label{eqn:NLS}
    i q_t + q_{xx} &= 2 q |q|^2
\end{align*}
for a wavefunction $q = q(t, x): \R \times \R \longrightarrow \C$.
While \NLS is most naturally associated with the zero boundary condition
\begin{align*} \tag{\ZBC} \label{eqn:ZBC}
    \lim_{|x| \rightarrow \infty} q(t, x) &= 0 \,,
\end{align*}
we are interested in the nonzero boundary condition
\begin{align} \label{eqn:NZBCt}
    \lim_{x \rightarrow \pm \infty} e^{2 i t} q(t, x) &= q_\pm \text{ where } |q_\pm| = 1 \,.
\end{align}
Because solutions of \NLS--\eqref{eqn:NZBCt} are not stationary at infinity, it is preferable to work instead with the so-called Gross-Pitaevskii equation
\begin{align*} \tag{\GP} \label{eqn:GP}
    i q_t + q_{xx} &= 2 q (|q|^2 - 1) \,,
\end{align*}
which is formally equivalent to \NLS under the transformation $q \mapsto e^{2 i t} q$.
The boundary condition \eqref{eqn:NZBCt} is replaced by
\begin{align*} \tag{\NZBC} \label{eqn:NZBC}
    \lim_{x \rightarrow \pm \infty} q(t, x) &= q_\pm \text{ where } |q_\pm| = 1 \,.
\end{align*}
From now on, we write \NLS for the system \NLS--\ZBC and \GP for the system \GP--\NZBC.
Before stating our main result, we discuss various aspects of \NLS, \GP, and the associated hierarchies.

\textbf{On applications of \NLS and \GP.}
The systems \NLS and \GP, also in their focusing variants, have been intensely studied due to their ubiquity in the analysis of wave phenomena in various physical systems.
We refer to the review papers \cite{Hasegawa} for applications in nonlinear optics and specifically high-speed communications, \cite{Bao} for applications in
Bose-Einstein condensation and plasma physics, and \cite{Berge} for applications in the study of wave collapse.

We remark that \GP is equivalent (under the Madelung transform, see \cite{Mohamad2014,Wegner2023}) to a quantum hydrodynamical system, which is relevant in the study of 
Bose-Einstein condensation \cite{Dalfovo1998,Grant1973}, superfluidity \cite{Feynman,Landau,Loffredo1993}, and quantum semiconductors \cite{Gardner}.

\textbf{Well-posedness results and function spaces in view of \NZBC.}
We start with a brief review of some well-posedness results for \NLS, focusing only on global results on $\R$ due to the breadth of the literature.
For global well-posedness results in Sobolev spaces $H^s(\R)$, we note in descending regularity \cite{BaillonCazenaveFigueira1977} for $s = 2$, \cite{GinibreVelo1979} for $s = 1$, \cite{Tsutsumi1987} for $s = 0$, and \cite{ChristCollianderTao2008} for a priori bounds and existence of weak solutions with $s > - \frac12$. 
In \cite{KochTataru2007,KochTataru2018,KochTataru2012} a priori bounds are given down to $s > - \frac12$.
Besides the Sobolev scale, there are global well-posedness results in Fourier-Lebesgue spaces $\ha{H}^s_r(\R)$ with $s \geq 0$ \cite{Grünrock2005} and modulation spaces $M^s_{p,q}(\R)$, also with $s \geq 0$ \cite{Klaus2023,OhWang2020}. 
We note also \cite{VargasVega2001}, where global well-posedness is shown in certain spaces larger than $L^2(\R)$.
For the reader interested in other settings, we refer to the books \cite{Cazenave2003,Tao2006} as well as the recent program started in \cite{IfrimTataru2023} by M. Ifrim and D. Tataru.

Because \NZBC is not preserved under addition and scalar multiplication of functions, it is not compatible with traditional scales of function spaces, such as the Sobolev scale.
As a result, the well-posedness of \GP is more difficult.
Spaces compatible with \ZBC contain the trivial background solution $q = 0$ to \NLS, while spaces compatible with \NZBC contain, for example, the trivial solution $q = 1$ of \GP, the stationary ``black soliton'' solution $q = \tanh$, or any of the time-dependent ``dark soliton'' solutions $q(t,x)=\Real[\zeta]+i\Imag[\zeta]\tanh(\Imag[\zeta](x+2\Real[\zeta]t))$, $\zeta \in \SSS^1$ to \GP. 
Generally, for any $q_\pm\in \mathbb{S}^1$, there exists a dark soliton profile connecting $q_-$ to $q_+$ (given explicitly in \eqref{eqn:dark-soliton} below). 

P. E. Zhidkov gave a first result in 1987 \cite{Zhidkov1987}, establishing local well-posedness in the Zhidkov space $Z^k(\R^n)$ with $n, k \in \N_{\geq 1}$, which is the closure of $C^\infty_c(\R^n)$ under the norm $\|\cdot\|_{Z^k(\R^n)} = \|u\|_{L^\infty(\R^n)} + \|\nabla u\|_{H^{k-1}}$.
This implies global well-posedness in $Z^1(\R)$ due to conservation of energy (see $\mc{H}^\GP_2$ below).
In \cite{Gallo} this is extended to $n = 2, 3$. A global result in the energy space was obtained by P. Gérard \cite{Gerard2006,Gerard2008} for $n = 1, 2, 3$, and for $n = 4$ under smallness assumptions.
The smallness assumptions for $n = 4$ are lifted in \cite{RomainTadahiroOanaMonica}. 
We note also \cite{Pecher2012} for a global well-posedness result permitting infinite energy in $1 + H^s(\R^3)$, $s \in \big(\frac56, 1)$,
and \cite{Antonelli2023} for more general nonlinearities in the cases $n = 2, 3$.
Recently in \cite{KochLiao,KochLiao2022}, H. Koch and X. Liao proved the global well-posedness of \GP for $s \geq 0$ in a complete metric space $(X^s, d^s)$, which we define below in \eqref{eqn:def-Xs}.
They use the complete integrability to construct conserved energies that control the solution at every regularity $s \geq 0$.
We make use of these energies in our globalization argument in Section \ref{proof:main-thm}.

\textbf{Complete integrability and conserved quantities of \NLS and \GP.}
The systems \NLS and \GP are completely integrable and have Lax pairs. 
A Lax pair for an evolutionary completely integrable PDE is a pair of operators $L$ and $P$, depending on a time-dependent potential function, such that the Lax equation $\del_t L = [P, L] = P L - L P$ holds true if and only if the potential is a solution to the PDE.
In the case of \NLS and \GP, it is shown in \cite{ZakharovShabat} that
\begin{align} \label{eqn:Lax}
    L = L^\NLS = L^\GP &= i \begin{pmatrix}
        \del_x & - q \\ \conj{q} & - \del_x
    \end{pmatrix}
\end{align}
and
\begin{align} \label{eqn:Pax}
    P^\NLS &= i \begin{pmatrix}
        2 \del_x^2 - q \conj{q} & - 2 q \del_x - q_x \\ 2 \conj{q} \del_x + q_x & - 2 \del_x^2 + q \conj{q}
    \end{pmatrix}
    & P^\GP &= i \begin{pmatrix}
        2 \del_x^2 - q \conj{q} + 1 & - 2 q \del_x - q_x \\ 2 \conj{q} \del_x + q_x & - 2 \del_x^2 + q \conj{q} - 1
    \end{pmatrix} \,.
\end{align}
A consequence of the Lax pair formulation is that eigenvalues of the Lax operator are stationary, and furthermore that the time evolution of the scattering data (including the reflection coefficient) is characterized by a \emph{linear} equation.
One may try to  recover the potential from the evolved scattering data. This is indeed possible for \NLS and \GP, where the so-called inverse scattering transform (IST) method can be applied.
For \NLS, the IST method was already well-developed in the seminal paper \cite{AblowitzKaupNewellSegur1974} (see also \cite{BorgheseJenkinsMcLaughlinKenneth2018} and references therein for more recent work). 
Recently, attention toward the nonzero boundary data case \GP, where the IST method is much harder to deploy, has grown, and significant progress has been made. 
We refer to the pioneering work \cite{ZakharovShabat} and the recent review \cite{Prinari}.

Another consequence of complete integrability is the existence of an infinite number of conserved quantities. 
For \NLS, we denote these by $\mc{H}^\NLS_n$ and list the initial ones:
\begin{align}
    \tag{Mass} \mc{H}^\NLS_0(q) &= \int_{\R} q \conj{q} \dd x
    \\ \tag{Momentum} \mc{H}^\NLS_1(q) &= - i \int_{\R} q \conj{q}_x \dd x
    \\ \tag{Energy} \mc{H}^\NLS_2(q) &= \int_{\R} q_x \conj{q}_x + q^2 \conj{q}^2 \dd x
    \\ \nonumber \mc{H}^\NLS_3(q) &= i \int_{\R} q \conj{q}_{xxx} - 4 q^2 \conj{q} \conj{q}_x - q_x q \conj{q}^2 \dd x
    \\ \nonumber \mc{H}^\NLS_4(q) &= \int_{\R} q_{xx} \conj{q}_{xx} - 6 q^2 \conj{q} \conj{q}_{xx} - 5 q^2 \conj{q}_x^2 
    - 6 q q_x \conj{q} \conj{q}_x - q q_{xx} \conj{q}^2 + 2 q^3 \conj{q}^3 \dd x \,.
\end{align}
Using the energy, \NLS can be written in Hamiltonian form
\begin{align*}
    i q_t &= \frac{\delta \mc{H}^\NLS_2(q)}{\delta \conj{q}} \,.
\end{align*}
The Hamiltonian structure is commonly considered for a pair of variables $(q, r)$, where later the identification $r = \conj{q}$ is made.
For convenience, we use $\conj{q}$ from the beginning, in return writing $q \conj{q}$ instead of $|q|^2$ when we intend to highlight the functional dependencies.

The Hamiltonians $\mc{H}^\NLS_n$ are formally still conserved under the flow of \GP, but the even ones are ill-defined while the odd ones are inconvenient in the setting of \NZBC.
We use instead a renormalized sequence of Hamiltonians $\mc{H}^\GP_n$, which are conserved quantities for \GP and compatible with \NZBC. 
Here the initial ones are:
\begin{align*}
    \tag{Mass} \mc{H}^\GP_0(q) &= \int_{\R} q \conj{q} - 1 \dd x
    \\ \tag{Momentum} \mc{H}^\GP_1(q) &= - i \int_{\R} q \conj{q}_x \dd x
    \\ \tag{Energy} \mc{H}^\GP_2(q) &= \int_{\R} q_x \conj{q}_x + (q \conj{q} - 1)^2 \dd x
    \\ \nonumber \mc{H}^\GP_3(q) &= i \int_{\R} q \conj{q}_{xxx} - 4 q^2 \conj{q} \conj{q}_x - q_x q \conj{q}^2 + 4 q \conj{q}_x \dd x
    \\ \nonumber \mc{H}^\GP_4(q) &= \int_{\R} q_{xx} \conj{q}_{xx} - 6 q^2 \conj{q} \conj{q}_{xx} - 5 q^2 \conj{q}_x^2 - 6 q q_x \conj{q} \conj{q}_x
    - q q_{xx} \conj{q}^2 
    \\ \nonumber &\qquad + 2 q^3 \conj{q}^3  + 6 q_x \conj{q}_x - 6 q^2 \conj{q}^2 + 6 q \conj{q} - 2 \dd x \,.
\end{align*}
In the setting of \NZBC there is an additional nontrivial conserved quantity, which plays no role in this work:
\begin{align*}
    \tag{Phase change} \mc{H}^\GP_{-1}(q) &= i \log\left(\frac{q_+}{q_-}\right) \,.
\end{align*}
Using the energy, \GP can be written in Hamiltonian form
\begin{align*}
    i q_t &= \frac{\delta \mc{H}^\GP_2(q)}{\delta \conj{q}} \,.
\end{align*}
The Hamiltonians $\mc{H}^\NLS_n$ and $\mc{H}^\GP_n$ satisfy
\begin{align}
    \label{eqn:NLS-GP-3} \frac{\delta \mathcal{H}^\GP_{2m}(q)}{\delta \conj{q}} &= \sum_{k=0}^m \binom{m-\frac12}{m-k} (- 4)^{m-k} \frac{\delta \mathcal{H}^\NLS_{2k}(q)}{\delta \conj{q}} \\
    \label{eqn:NLS-GP-4} \frac{\delta \mathcal{H}^\GP_{2m+1}(q)}{\delta \conj{q}} &= \sum_{k=0}^m\binom{m}{m-k} (- 4)^{m-k} \frac{\delta \mathcal{H}^\NLS_{2k+1}(q)}{\delta \conj{q}} \,.
\end{align}
More specifically, the densities of $\mc{H}^\GP_n$ are affine linear combinations of the densities of $\mc{H}^\NLS_n$, and vice versa (see \eqref{eqn:NLS-GP-1}--\eqref{eqn:NLS-GP-2} below, and also \cite[(10.25)]{FaddeevTakhtajan}).

\textbf{The \NLS and \GP hierarchies.}
It is natural to consider for $n \geq 0$ the infinite sequence of Hamiltonian PDEs
\begin{align} \label{eqn:NLS-n-alt}
    i q_t &= \frac{\delta \mc{H}^\NLS_n(q)}{\delta \conj{q}} \,.
\end{align}
This is referred to as the \NLS hierarchy. Some caution is needed, as sometimes only the flows where $n$ is even are called the \NLS hierarchy, while the flows where $n$ is odd are called the complex \mKdV hierarchy.
The Hamiltonians $\mc{H}^\NLS_n$ are conserved quantities for every one of these flows. Equivalently, for any $n, m \geq -1$ we have
\begin{align}
    \left\{\mc{H}^\NLS_n, \mc{H}^\NLS_m\right\} &\coloneqq \frac{\delta \mc{H}^\NLS_n}{\delta q} \frac{\delta \mc{H}^\NLS_m}{\delta \conj{q}} - \frac{\delta \mc{H}^\NLS_n}{\delta \conj{q}} \frac{\delta \mc{H}^\NLS_m}{\delta q} 
    = 0 \,,
\end{align}
i.e. they pairwise Poisson commute (see \cite[III.§2]{FaddeevTakhtajan} or \cite[Theorem B.7]{KlausKochLiu2023}). 
As a result, their Hamiltonian flows commute and we can consider instead for a single function $q = q(t_0, t_1, t_2, \dots, x): \R^{(\N)} \times \R \longrightarrow \C$
the infinite hierarchy of PDEs
\begin{align*} \tag{$\NLS_n$} \label{eqn:NLS-n}
    i \del_{t_n} q &= \frac{\delta \mc{H}^\NLS_n(q)}{\delta \conj{q}} \,.
\end{align*}

As is the case for \NLS, solutions to \eqref{eqn:NLS-n} with nonzero boundary data are not stationary at infinity,
and so it is preferable to work instead with the infinite sequence of Hamiltonian PDEs
\begin{align} \label{eqn:GP-n-alt}
    i q_t &= \frac{\delta \mc{H}^\GP_n(q)}{\delta \conj{q}} \,,
\end{align}
which we call the \GP hierarchy.
Of course, the same Poisson commutation relations $\left\{ \mc{H}^\GP_n, \mc{H}^\GP_m \right\} = 0$, $n, m \geq 1$ hold,
so we can again consider instead for a single function $q = q(t_0, t_1, \dots, x)$ the infinite hierarchy of PDEs
\begin{align*} \tag{$\GP_n$} \label{eqn:GP-n}
    i \del_{t_n} q &= \frac{\delta \mc{H}^\GP_n(q)}{\delta \conj{q}} \,.
\end{align*}

\textbf{On applications of the \NLS, \GP, and other integrable hierarchies.}
First and foremost, the study of any particular higher equation of an integrable hierarchy is also the study of every other equation in the hierarchy due to Noether's theorem:
the flows generated by the higher Hamiltonians are symmetries for every other flow, in particular the usually eponymous initial PDE (\NLS, \KdV, etc.) of the hierarchy in consideration.
As such, the well-posedness result for the \NLS hierarchy that we present here represents the construction of an infinite number of symmetries for \GP.
Integrable equations like \NLS and \KdV appear frequently as amplitude equations for the long-wave regime of other equations (see e.g. \cite{Berge,SchneiderWayne2000,TotzWu2012,Whitham1965,Zakharov1968}).
In such cases, they may be modulated using their symmetries, including the higher ones (see \cite{Bambusi} for an example). 
On another note, integrable hierarchies are also relevant in mathematical physics, in connection with algebraic geometry \cite{London-1,London-2}.

\textbf{Renormalization of hierarchies.}
Let us state clearly that taking affine linear combinations of the (densities of the) Hamiltonians in an integrable hierarchy yields flows which are essentially equivalent to those already in the hierarchy.
We make this precise with a proposition.
\begin{proposition}[Equivalence between hierarchies] \label{prop:hierarchy-equivalence}
    Let $A = (A_{n,k}) \in \C^{\N \times \N}$ be an invertible infinite lower triangular matrix and set $\bm{t} = (t_0, t_1, \dots)$. 
    Formally, $q(\bm{t}, x)$ solves \eqref{eqn:NLS-n} for all $n \in \N$ if and only if $u(\bm{t}, x) = q(A \bm{t}, x)$ solves 
    \begin{align*}
        i \del_{t_n} u &= \sum_{k=0}^n A_{n,k} \frac{\delta \mc{H}^\NLS_k(u)}{\delta \conj{u}} \,.
    \end{align*}
    In this sense, the well-posedness theories of the \NLS hierarchy and this renormalized hierarchy are equivalent.
\end{proposition}
\begin{proof}
    The Hamiltonians only depend on $u$, $\conj{u}$ and their spatial derivatives, but not on $t$, so the proof follows from the chain rule.
\end{proof}
The relations \eqref{eqn:NLS-GP-3}--\eqref{eqn:NLS-GP-4} represent precisely such an equivalence between the \NLS and \GP hierarchies. Here
\begin{align*}
    A_{2m,2k} &= \binom{m-\frac12}{m-k} (-4)^{m-k} & A_{2m,2k+1} &= 0
    \\ A_{2m+1,2k} &= 0 & A_{2m+1,2k+1} &= \binom{m}{m-k} (-4)^{m-k} \,,
\end{align*}
and this matrix is invertible (replace $(- 4)^{m-k}$ by $4^{m-k}$).
Indeed, our main result concerns the well-posedness of the \GP hierarchy instead of the \NLS hierarchy,
and in the proof we introduce yet another renormalization for the odd flows (see \eqref{eqn:tGP-n} below).
As an example, in the case $n = 2$ the energy of \NLS is renormalized with the mass to yield the energy of \GP.

\subsection{Main result}
Let $d \in \N$. We denote by $\mc{D}(\R^d; \C) = C_c^\infty(\R^d; \C)$ the space of test functions.
For a function $f \in \mc{D}(\R^d; \C)$, we define the Fourier transform
\begin{align*}
    \ha{f}(\xi) &= \frac{1}{(2 \pi)^{\frac{d}{2}}} \int_{\R^d} e^{- i x \cdot \xi} f(x) \dd x
\end{align*}
and extend it to tempered distributions in the usual manner.
Before stating our main result, we define for $s, s' \in \N$ the (weighted) Sobolev spaces
\begin{align*}
    \big( H^s, \|\cdot\|_{H^s} \big) &= \big( \big\{ u \in L^2: (1 + |\xi|^2)^{\frac{s}{2}} \ha{u} \in L^2 \big\}, \|(1 + |\xi|^2)^{\frac{s}{2}} \, \ha{\cdot} \,\|_{L^2} \big)
    \\ \big( H^{s',1}, \|\cdot\|_{H^{s',1}} \big) &= \big( \big\{ u \in H^{s'}: x u \in H^{s'} \big\}, \|x \cdot\|_{H^{s'}} + \|\cdot\|_{H^{s'}} \big)
\end{align*}
and the energy functionals
\begin{align*}
    E^s(q) &= \||q|^2 - 1\|_{H^{s-1}} + \|q_x\|_{H^{s-1}}
    & E^{s',1}(q) &= \||q|^2 - 1\|_{H^{s'-1,1}} + \|q_x\|_{H^{s'-1,1}} \,.
\end{align*}
The notations $\ceil{\cdot}$ resp. $\floor{\cdot}$ mean $\ceil{m}=\min\{k\in \Z,\quad k\geq m\}$ resp. $\floor{m}=\max\{k\in \Z,\quad k\leq m\}$.
\begin{theorem}[Global well-posedness of the \GP hierarchy] \label{thm:1}
    Let $n \in \N$ with $n \geq 2$ and define $m = \frac{n-1}{2}$. Let $s, s' \in \N$ such that $s \geq 2n + m$ and $\frac{s+m}{2} \leq s' \leq s - n$.
    Equation \eqref{eqn:GP-n} is globally well-posed in the following sense:
    For every $q_\ast \in H^{s+1-\ceil{m}+n}_{\loc}$ with $E^{s+1-\ceil{m}+n,1}(q_\ast) < \infty$ and every $p_0 \in H^s \cap H^{s',1}$, there exists a unique $p \in C(\R; H^s \cap H^{s',1})$ with $p(0) = p_0$ such that $q_\ast + p(t)$ is a solution of \eqref{eqn:GP-n} in the sense of distributions.
    The solution is strong in the sense that $p$ solves the corresponding perturbative formulation (see \eqref{eqn:struct-GP-5} below) strongly in $H^{s-n} \cap H^{s'-n,1}$.
    For every $T > 0$, the map
    \begin{align*}
        H^s \cap H^{s',1} &\longmapsto C_b([-T, T]; H^s \cap H^{s',1})
        \\ p_0 &\longmapsto p
    \end{align*}
    is Lipschitz continuous. Lastly, there exists $C = C(n, s, E^{s+1-\ceil{m}+n,1}(q_\ast), \|p_0\|_{H^s \cap H^{s',1}}) > 0$ such that for all $t \geq 0$ we have
    \begin{align*}
        \|p(t)\|_{L^\infty_x} + \|p_x(t)\|_{H^{s-1}_x} &\leq C
        & \|p(t)\|_{H^{s',1}_x} &\leq C e^{C t} \,.
    \end{align*}
\end{theorem}
\begin{remark}[On the proof] \label{rem:1}
    The difficulty lies in the nonzero boundary condition, which is why we solve for a localized perturbation $p$ of a time-independent front $q_\ast$.
    The main step of the proof is a study of the coefficients appearing in \eqref{eqn:GP-n} in order to show that the equation for $p$ is part of a broad class of nonlinear dispersive PDEs (Proposition \ref{prop:2}) for which we can show local well-posedness (Theorem \ref{thmlwp}) in $H^s \cap H^{s',1}$.
    In \cite{KochLiao, KochLiao2022} the authors construct conserved energies for \GP that control the solution at every level of regularity in certain metric spaces. 
    In particular, the energies $E^s(q)$ remain bounded. These energies are strong enough to prevent finite-time blow up of the local solution $p$ in $H^s \cap H^{s',1}$ (Lemma \ref{lem:12}), allowing us to conclude global well-posedness.
    The reason the $H^{s',1}$-norm of $p$ is not globally bounded is that we cannot directly control it with the energies, instead having to use Grönwall's inequality in an energy estimate.
\end{remark}
\begin{remark}[On the energy assumption and \NZBC]
    The assumption $E^{s+1-\ceil{m}+n,1}(q_\ast) < \infty$ implies that $\lim_{x \rightarrow \pm \infty} q_\ast(x) = q_\pm$ for some $q_\pm \in \SSS^1$.
    Conversely, we can connect any end states $q_\pm$ with a dark soliton profile (see \eqref{eqn:dark-soliton}) that satisfies $E^{s+1-\ceil{m}+n,1}(q_\ast) < \infty$.
    One may be tempted to consider the more general nonzero boundary condition $|q_\ast|^2 - 1 \in L^2$ and $(q_\ast)_x \in L^2$.
    This permits $(q_\ast)_x \not\in L^1$, and hence solutions which do not have limits $q_\pm$ at $\pm \infty$ could be studied.
    It turns out that the well-posedness theory in Section \ref{section:4} requires $|q_\ast|^2 - 1$ and $(q_\ast)_x$ to basically be admissible as substitutes for $p$ and $p_x$. 
    Correspondingly, we can only work in settings where $(q_\ast)_x \in H^{s-\ceil{m}+n,1}$ for some large $s$. 
    This implies $(q_\ast)_x \in L^1$, so our method cannot consider boundary conditions more general than \NZBC.
\end{remark}
\begin{remark}[On the time-independence of $q_\ast$ and stability of solitons]
    It would be natural to choose $q_\ast$ as a moving front, for example a soliton, and phrase the theorem as a stability result.
    Currently we choose $q_\ast$ to be time-independent, which causes terms depending only on $q_\ast$ to appear on the right-hand side in \eqref{eqn:nonlinear} (this also causes $g \neq 0$ in \eqref{eqn:dispersive-nonlinear}).
    Removing these terms would lead to minor improvements of certain estimates and regularity assumptions.
    Since our well-posedness theory can handle these inhomogeneous terms, and for the sake of simplicity, we choose $q_\ast$ to be time-independent.
    Note that when viewed through the lens of stability, our result has the advantage of not requiring smallness of the perturbation $p$.
\end{remark}

\subsection{Organization of the paper and further results}

All sections and appendices are largely self-contained and can in principle be read in any order.

\subsubsection{Overview}
In Subsection \ref{subsection:1.1.2} we motivate Section \ref{section:2}, which is not directly relevant for the proof of the main result and instead is concerned with the rigorous definition and asymptotic expansion of the logarithm of the transmission coefficient associated to the Lax operator $L$ (see \eqref{eqn:Lax}).
In addition, we review some of the literature on the IST method for \NLS/\GP with a nonzero background and elaborate on our contribution.
Section \ref{section:2} concludes with the derivation of the recurrence relation \eqref{eqn:sig-3}, which is fundamental for the computations in Section \ref{section:3}.
    
In Subsection \ref{subsection:1.1.3} we discuss Section \ref{section:3}, which contains an essential result on the structure of the \NLS and \GP hierarchies. 
Using the method of generating functions, we are able to give explicit formulas for all terms in the hierarchies containing at most one factor that is $q_x$, $\conj{q}_x$, or a derivative thereof.
Section \ref{section:3} concludes with the proof of Proposition \ref{prop:2}, which is a perturbative formulation of the \GP hierarchy that fits into the well-posedness theory developed in Section \ref{section:4}.
    
In Subsection \ref{subsection:1.1.4} we discuss Section \ref{section:4}. We begin with a review of different approaches to the well-posedness of the \NLS and \KdV hierarchies. 
Subsequently, we motivate our choice to use techniques developed by C. E. Kenig, G. Ponce and L. Vega in the 1990s.
The key result is the local well-posedness of a large class of dispersive nonlinear systems.
Section \ref{section:4} concludes with the proof of Theorem \ref{thm:1}, our main result.

\subsubsection{Section 2: Direct scattering theory revisited and further developed} \label{subsection:1.1.2} 
In Section \ref{section:2} we recall some aspects of the Lax operator $L$, associated to \GP, in the setting of \NZBC,
with the goal of rigorously defining the transmission coefficient $a(\lambda)$ and proving that its logarithm has asymptotic expansions in powers of the spectral parameters $\lambda$ and $z = \sqrt{\lambda^2 - 1}$ at infinity.
We define Hamiltonians $\mc{H}^\NLS_n$ and $\mc{H}^\GP_n$ as the expansion coefficients (see \eqref{eqn:trans-expand-GP}--\eqref{eqn:trans-expand-NLS}):
        \begin{align*} 
           \log a(\lambda) \sim i \sum_{n=0}^\infty \frac{\mc{H}^\NLS_n}{(2 \lambda)^{n+1}} \, 
            \qquad\qquad  \log a(\lambda) \sim i \sum_{n=0}^\infty \frac{\mc{H}^\GP_n}{(2 z)^{n+1}} \, 
        \end{align*}   and in the process we derive a recurrence relation for their densities (see \eqref{eqn:sig-3} below), which is the basis of our calculations in Section \ref{section:3}.
This recurrence relation is well-known and also appears in the seminal paper \cite{ZakharovShabat} by V. E. Zakharov and A. B. Shabat, which initiated the study of the IST method for \GP with \NZBC.
We refer to \cite{BiondiniKovacic,CuccagnaJenkins,DemontisVanDerMeeCornelis2023,DemontisPrinariMeeVitale,DemontisPrinariVanDerMeeVitale2014,VanDerMee2015} and the review \cite{Prinari} for further development of the IST method for \NLS-type integrable equations with nonzero boundary data.

Note that most of these references are concerned with the formal and rigorous establishment of the inverse scattering transform, 
but not the rigorous asymptotic expansion of the logarithm of the transmission coefficient.
The standard reference for this is the monograph \cite[I.§3-§8]{FaddeevTakhtajan}, where first the periodic case on an interval of length $L$ is considered, and then the limit $L \rightarrow \infty$ is taken.
More recently in \cite{KochLiao,KochLiao2022}, H. Koch and X. Liao have defined and rigorously expanded the transmission coefficient in the very general setting $E^s(q) < \infty$ for $s \geq 0$, by finding a change of variables for the Zakharov-Shabat scattering problem that essentially causes only the quantities $|q|^2 - 1$ and $q_x$ to appear.

We believe there is expository value in a direct path from the Zakharov-Shabat scattering problem to the rigorous expansion of the transmission coefficient, using the classical approach described in \cite{FaddeevTakhtajan}, but not working with the periodic case first.
Specifically, we use the integral representation Ansatz \eqref{eqn:int-4}, which to the best of our knowledge is not present in the literature.
In order to explain, let us assume as given the Jost solution $\Phi^\pm_1(x, \lambda)$ with $x \in \R$ and $\lambda \in \C$ to the Zakharov-Shabat scattering problem
\begin{align*} \label{eqn:ZS} \tag{\ZS}
    L \Phi^\pm_1(x, \lambda) = \lambda \Phi^\pm_1(x, \lambda) \,, \quad \text{or equivalently} \quad \del_x \Phi^\pm_1(x, \lambda) &= \begin{pmatrix}
        - i \lambda & q(x)
        \\ \conj{q}(x) & i \lambda
    \end{pmatrix} \Phi^\pm_1(x, \lambda) \,,
\end{align*}
for a potential of interest $q$. In addition, we assume that for a reference potential $q_\ast$ the corresponding Jost solution $\Phi^\pm_{\ast,1}$ is known.
We now make the triangular representation Ansatz
\begin{align} \label{eqn:int-demo}
    \Phi^\pm_1(x, \lambda) &= \Phi^\pm_{\ast,1}(x, \lambda) + \int_{\pm \infty}^x \Gamma^\pm(x, y) \Phi^\pm_{\ast,1}(y, \lambda) \dd y \,,
\end{align}
where $\Gamma^\pm(x, y)$ is a kernel that is independent of $\lambda$. 
Such triangular representations are well-known in the case $q_\ast = q_\pm$ (see e.g. \cite{DemontisPrinariVanDerMeeVitale2014,FaddeevTakhtajan,VanDerMee2021}).
Since the kernel $\Gamma^\pm$ is independent of $\lambda$, these integral representations are suitable for asymptotic analysis in the parameter $\lambda$ of the Jost solutions via integration by parts.
In fact, the density of the logarithm of the transmission coefficient $\log a(\lambda)$ is $\del_x \log \Psi^-_{1,1}$, i.e. the logarithmic derivative of the first component of the modified Jost solution $\Psi^-_{1,1}$ (see \eqref{eqn:trans-5} below).
Therefore, we can use \eqref{eqn:int-demo} to derive the desired asymptotic expansion for $\log a(\lambda)$, as long as the error terms can be controlled in $L^1_x$.
An essential source of integrability when working with \eqref{eqn:int-demo} is the difference $q - q_\ast$, but the choice $q_\ast = q_\pm$ does not yield any integrability at $\mp \infty$.
We therefore propose to choose the reference potential $q_\ast$ as the dark soliton profile with boundary data $q_\pm$ instead (see \eqref{eqn:dark-soliton} below), so that $q - q_\ast$ can have decay at both infinities.
Conveniently, explicit formulas for the Jost solutions of the dark soliton are given in \cite{ChenChenHuang}.

It is important for $\Phi^\pm_{\ast,1}$ to be a solution of the Zakharov-Shabat scattering problem \eqref{eqn:ZS} in order for $\Gamma^\pm$ to solve a boundary value problem which does not involve $\lambda$. This is the system \eqref{eqn:603}--\eqref{eqn:605}, which generalizes \cite[Chapter 1, (8.18)--(8.19)]{FaddeevTakhtajan}.
We subsequently state and prove a well-posedness result for $\Gamma^\pm$ in Lemma \ref{lem:678}, which we believe to be of independent interest. Note that we work in the setting where $q$ and $q_\ast$ are smooth, and $q - q_\ast$ is Schwartz. 
The reader interested in weaker assumptions may adapt the proof for this purpose.
The estimates we obtain allow us to deduce the desired expansion of the transmission coefficient in Lemma \ref{lem:1}.

We see potential for the triangular representation \eqref{eqn:int-demo} to be useful in other aspects of the IST method, such as WKB expansions of various scattering data.
Furthermore, it may be of use in the setting with asymmetric boundary conditions (i.e. $|q_-| \neq |q_+|$), when a Jost solution $\Phi^-_{\ast,1}$ for a reference potential $q_\ast$ is known.

\subsubsection{Section 3: Explicit formulas for parts of the hierarchies and the perturbative setting} \label{subsection:1.1.3} 
In Section \ref{section:3} we determine the coefficients of all terms in the \NLS hierarchy which have at most one factor that is $q_x$, $\conj{q}_x$, or a derivative thereof.
When using a perturbative Ansatz $q(t, x) = q_\ast(x) + p(t, x)$, these are precisely the terms which we need to know in order to determine the linear part of the PDE for $p$. 

Let $[u^j]$ denote the functional that extracts the coefficient in front of $u^j$ from a formal power series in the symbol $u$, e.g. $[u^2] (1 + u)^3 = 3$.
With this notation, we define the coefficients
\begin{align}
    \label{eqn:coeffs} J_{n,j} &= \begin{cases} 
        0 &, n \text{ even} \\
        [u^j] 2 (1 + 4 u)^{\frac{n-2}{2}} &, n \text{ odd}
    \end{cases}
    & K_{n,j} &= \begin{cases}
        [u^j] (1 + 4 u)^{\frac{n-1}{2}} &, n \text{ even} \\
        [u^j] (- 1 - 2 u) (1 + 4 u)^{\frac{n-2}{2}} &, n \text{ odd}
    \end{cases} \,.
\end{align}
\begin{theorem}[Structure of the \NLS and \GP hierarchies] \label{thm:2}
    The following hold true.
    \begin{enumerate}[(i)]
        \item Let $m = \floor{\frac{n}{2}} \geq 1$. The equations of the \NLS hierarchy have the form
        \begin{align}
            \label{eqn:struct-NLS-1} i \del_{t_{2m}} q 
            &= \sum_{j=0}^{m-2} J_{2m+1,j} q^{j+2} \conj{q}^j (i \del_x)^{2m - 2 - 2 j} \conj{q} 
            - \sum_{j=0}^{m-1} K_{2m+1,j} q^j \conj{q}^j (i \del_x)^{2 m - 2 j} q
            \\ \nonumber &+ (m+1) C_m \conj{q}^m q^{m+1} + \mc{O}^{2,n-2}_q(q_x)
            \\ \label{eqn:struct-NLS-2} i \del_{t_{2m+1}} q &= \sum_{j=0}^m K_{2m+2,j} q^j \conj{q}^j (i \del_x)^{2 m + 1 - 2 j} q + \mc{O}^{2,n-2}_q(q_x) \,.
        \end{align}
        Here $\mc{O}^{2,n-2}_q(q_x)$ refers to the class of polynomial expressions in $q$, $\conj{q}$, and their derivatives, for which
        each monomial has at most $n - 2$ derivatives in total, and at least two factors are $q_x$, $\conj{q}_x$, or derivatives thereof. 
        This notation is defined fully in \eqref{eqn:O-def} below. The coefficients $C_m$ are the Catalan numbers (see \eqref{eqn:catalan}).

        \item The following is a formal statement. Setting $\rho = q \conj{q}$, we have as $|\lambda| \rightarrow \infty$ the asymptotic expansion
        \begin{align*}
            & \frac{\delta \log a(\lambda)}{\delta \conj{q}} 
            \sim \sum_{n=0}^\infty \frac{\mc{O}^{2,n}_q(q_x)}{(2 i \lambda)^{n+1}}
            \\ &+ \frac{i}{2 (\lambda^2 - \rho)^{\frac12}} \Bigg( 
                q + \frac{2 \lambda + i \del_x}{4 (\lambda^2 - \rho) + \del_x^2} [i q_x]
                + \frac{1}{2 (\lambda^2 - \rho) (4 (\lambda^2 - \rho) + \del_x^2)} [q \del_x^2 (\lambda^2 - \rho)]
            \Bigg) \,.
        \end{align*}
        This is to be understood in the sense of expansion in $\lambda$ by the geometric series and subsequent comparison of orders.
        Note that $\del_x^2$ and the multiplication operator $\rho$ commute up to terms in $\mc{O}^{1,2}_q(q_x)$.
        The square brackets denote operator application.

        \item For $m \geq 2$ the equations of the \GP hierarchy have the form
        \begin{align*}
            \label{eqn:struct-GP-3} \tag{$\GP_{2m}$} i \del_{t_{2m}} q 
            &= \left( (i \del_x)^{2m} + 2 (i \del_x)^{2m-2} \right) q + 2 q^2 (i \del_x)^{2m-2} \conj{q} 
            + \mc{O}^{2,n-2}_{q, |q|^2 - 1}(q_x, |q|^2 - 1)
            \\ \label{eqn:struct-GP-4} \tag{$\GP_{2m+1}$} i \del_{t_{2m+1}} q &= \sum_{j=0}^m 4^{m-j} \binom{\frac12}{m-j} (i \del_x)^{2j+1} q 
            + \mc{O}^{2,n-2}_{q, |q|^2 - 1}(q_x, |q|^2 - 1) \,.
        \end{align*}
        Here $\mc{O}^{2,n-2}_{q, |q|^2 - 1}(q_x, |q|^2 - 1)$ refers to the class of polynomial expressions in $q$, $\conj{q}$, $|q|^2 - 1$, and their derivatives, for which
        each monomial has at most $n - 2$ derivatives in total, and at least two factors are $q_x$, $\conj{q}_x$, $|q|^2 - 1$, or derivatives thereof.
        We refer again to \eqref{eqn:O-def} for the full definition of this notation.
    \end{enumerate}
\end{theorem}
\begin{remark}
    To the best of our knowledge these coefficients are not known in the literature. 
    For further efforts to calculate explicit coefficients in integrable hierarchies, we refer to \cite{Adams2025}, where the coefficients of all cubic terms in the \dNLS hierarchy are calculated, and also \cite{AvramidiSchimming} for the \KdV hierarchy.
    The general theory in \cite{Dickey1991} may also be of use.
\end{remark}
\begin{remark}
    We do not have an application for (ii) and only state it because we consider it an elegant reformulation of (i).
    We would like to point out that it might ``by chance'' describe more terms in the hierarchy than \eqref{eqn:struct-NLS-1}--\eqref{eqn:struct-NLS-2} do.
    To explain, observe that different arrangements of the order of $\del_x^2$ and multiplication by $\rho$ in the geometric series expansion of the resolvent
    yield different residual terms in $(2 i \lambda)^{-n-1} \mc{O}^{2,n}_q(q_x)$. 
    One may hope that by choosing the correct ordering the residual terms can vanish, i.e. (ii) would become an exact formula.
    We have concluded, after testing various natural arrangements, that this does not seem to be the case.
\end{remark}
Given a choice of $q_\ast$, we can now write down a dispersive PDE for $p$.
\begin{proposition}[Perturbative formulation of the \GP hierarchy] \label{prop:2}
    Let $q(t, x) = q_\ast(x) + p(t, x)$, where $\lim_{x \rightarrow \pm \infty} q_\ast(x) = q_\pm \in \SSS^1$.
    Let $n \geq 1$. We consider the extended system for $\bm{p} = (\bm{p}_j)_{1 \leq j \leq 4} = (p, q_\ast^2 \conj{p}, \conj{p}, \conj{q}_\ast^2 p)$. 
    Then $q$ solves \eqref{eqn:GP-n} if and only if $\bm{p}$ solves
    \begin{align} \label{eqn:struct-GP-5}
        i \del_{t_n} \bm{p} &= \mc{L}^n[\bm{p}] + \mc{N}^n[\bm{p}] \,,
    \end{align}
    where with $D_x = - i \del_x$ we have
    \begin{align}
        \label{eqn:lin-even} \mc{L}^{2m} &= D_x^{2m-2} \begin{pmatrix}
        D_x^2 + 2 & 2 & 0 & 0
        \\ - 2 & - D_x^2 - 2 & 0 & 0
        \\ 0 & 0 & - D_x^2 - 2 & - 2
        \\ 0 & 0 & 2 & D_x^2 + 2
    \end{pmatrix}
        \\ \label{eqn:lin-odd} \sum_{k=0}^m \binom{-\frac12}{m-k} 4^{m-k} \mc{L}^{2k+1} &= - D_x^{2m+1} \,,
    \end{align}
    and for $d \in \{1, 2, 3, 4\}$ we have
    \begin{align} \label{eqn:nonlinear}
        (\mathcal{N}^n[p])_d &= \mathds{1}_{\{n=2\}} \mc{O}^{1,0}_{q_\ast, |q_\ast|^2 - 1} (|q_\ast|^2 - 1) + \mc{O}^{1,n}_{q_\ast}((q_\ast)_x) + \mc{O}^{2,n-2}_{p, q_\ast, |q_\ast|^2 - 1}(p, (q_\ast)_x, |q_\ast|^2 - 1) \,.
    \end{align}
    This notation is defined in \eqref{eqn:O-def}.
\end{proposition}
The key points here are that the linear part is benign and that each term in the nonlinear part has sufficient integrability,
in the sense that it has either coefficients which provide integrability, or two factors which are derivatives of $p$. 
\begin{remark}
    We switch to the vector variable $\bm{p}$ in order to eliminate the non-constant coefficient $q^2$ that can be seen in \eqref{eqn:struct-GP-3} and would otherwise appear in the linear part of the equation for $p$.
    This is not necessary for the odd flows, but we treat them the same way for the sake of uniformity.
\end{remark}

\subsubsection{Section 4: Well-posedness of a large class of nonlinear dispersive systems including the \GP hierarchy} \label{subsection:1.1.4} 
In Section \ref{section:4} we construct a local well-posedness theory for a large class of dispersive nonlinear systems of PDEs that includes the perturbative formulation of the \GP hierarchy in Proposition \ref{prop:2}.
We then deduce global well-posedness, proving Theorem \ref{thm:1} by using the conserved energies constructed by H. Koch and X. Liao in \cite{KochLiao,KochLiao2022}.

Before we explain our argument in detail, let us review existing well-posedness results for integrable hierarchies and general dispersive nonlinear equations.
Due to their nature, when working with integrable hierarchies, one expects to have available a plethora of conserved quantities. 
Although the construction of \emph{useful} conserved quantities is not a trivial matter (see for example \cite{KochLiao,KochLiao2022,KochTataru2018} and \cite{KillipVisanHarropGriffiths2024,KillipVisanZhang2018}),
it is greatly simplified when working with high, integer regularity, as we do. 
Therefore, local well-posedness is the essential problem, and here it is not easy to derive a benefit from the integrability of the equation under consideration.
As such, the local well-posedness of large, not necessarily integrable, classes of dispersive nonlinear equations may be studied, and in return integrable hierarchies serve as natural applications for such theories.

We start with the \KdV hierarchy, which we define for a function $u = u(t_1, t_2, \dots, x): \R^{(\N)} \times \R \rightarrow \R$ as the infinite sequence of PDEs 
\begin{align*} \tag{$\KdV_n$} \label{eqn:KdV-n}
    \del_{t_n} u &= \frac12 \del_x \frac{\delta \mc{E}^\KdV_n(u)}{\delta u} 
    \qquad \mc{E}^\KdV_n(u) = \int_{\R} \sigma^\KdV_n(u) \dd x 
    & \\ \nonumber \sigma^\KdV_0(u) &= 0 
    \qquad \sigma^\KdV_1(u) = u 
    \qquad \sigma^\KdV_{n+1}(u) = \del_x \sigma^\KdV_n(u) + \sum_{k=0}^n \sigma^\KdV_k(u) \sigma^\KdV_{n-k}(u)\quad \forall n\geq 1.
\end{align*}
J.-C. Saut first proved in \cite{Saut1979} the existence of global distributional solutions to \eqref{eqn:KdV-n} for initial data in Sobolev spaces. 
Subsequently, M. Schwarz in \cite{Schwarz1984} showed global existence and uniqueness of solutions in Sobolev spaces on the torus, again at high regularity.

We are most interested in a theory developed in the 1990s by C. E. Kenig, G. Ponce and L. Vega \cite{KenigPonceVega-BenjaminOno,KenigPonceVega1991-1,KenigPonceVega-SmallSolutions,KenigPonceVega1993,KenigPonceVega1991-2}.
Here local smoothing and maximal function estimates for the linear evolution are proven and used to derive local and global well-posedness results for dispersive PDEs for large classes of linear and nonlinear parts.
We focus in particular on the papers \cite{KenigPonceVega1994-1,KenigPonceVega1994-2}, where local (and in some cases global) well-posedness results in weighted Sobolev spaces are proven for a class of equations that includes the \KdV hierarchy.

We use their methods to set up a local well-posedness theory for a large class of dispersive nonlinear systems of PDEs that covers the case of Proposition \ref{prop:2}.
Before we provide the details, let us discuss more recent methods and results and why we do not use them here. 

With the aim of finding an $L^2$-based well-posedness theory, D. Pilod \cite{Pilod} studied a family of equations similar to the \KdV hierarchy, but with only quadratic nonlinearities, and proved local well-posedness for small initial data in $H^s(\R)$.
However, they also proved failure of the solution map to be $C^2$ at zero.
Subsequently in \cite{KenigPilod}, D. Pilod and C. E. Kenig proved global well-posedness in $H^s(\R)$ for initial data of arbitrary size for a class of equations that includes the \KdV hierarchy, using subtle energy estimates and parabolic regularization. Both of these papers need to assume smallness of the initial data, which is lifted in the second paper only by using the scaling symmetry.
We want to avoid smallness assumptions and have no scaling symmetry available in the perturbative setting of Proposition \ref{prop:2}.

Based on their work on the Fourier restriction norm method and the usage of Fourier-Lebesgue spaces for dispersive nonlinear equations (see e.g. \cite{Grünrock2004}), A. Grünrock proved in \cite{Grünrock2010} well- and ill-posedness results in the aforementioned spaces for the \mKdV and \KdV hierarchies at low regularity.
Recently, J. Adams continued work in this direction \cite{Adams2025,Adams2024}, proving well- and ill-posedness results for the \NLS- and \dNLS hierarchies in Fourier-Lebesgue and modulation spaces at low regularity.

We certainly expect the Fourier restriction norm method to be applicable to our setting, but not without considerable work. 
Besides the fact that \eqref{eqn:struct-GP-5} contains coefficients in the nonlinearity, we face the following complications: in Proposition \ref{prop:2} we are dealing with a genuine \emph{system}, which is only diagonalizable with a singular change of variables (see \eqref{eqn:diagonalization}); has nonpolynomial dispersion relations in diagonal form; and lastly contains quadratic and even linear terms (with benign coefficients) in the nonlinearity.

We mention also the method of commuting flows, which was developed and used by R. Killip and M. Vișan to prove global well-posedness of \KdV in $H^{-1}(\R)$ in \cite{KillipVisan2019}, and subsequently applied to other equations such as Benjamin-Ono \cite{KillipLaurensVisan2024} and \NLS \cite{KillipVisanHarropGriffiths2024}.
We would like to mention that in \cite{KillipVisan2019} the flow of \KdV is approximated by a flow whose Hamiltonian involves the logarithm of the transmission coefficient. One may understand this as a modulation of the solution along the higher flows, giving an example of higher symmetries aiding in the understanding of an equation in an integrable hierarchy.
In the paper \cite{BringmannKillipVisan2021} by B. Bringmann and the aforementioned authors, this method was applied to the fifth-order equation in the \KdV hierarchy, and finally in the work \cite{KlausKochLiu2023} by F. Klaus, H. Koch and B. Liu to the whole hierarchy, proving global well-posedness in $H^{-1}(\R)$.
We again believe that it is in principle possible to apply this method to the \NLS hierarchy with nonzero boundary data, but recognize that it would take considerable work.

For these reasons, we focus on the classical theory developed by C. E. Kenig, G. Ponce and L. Vega, working with high regularity and the weighted Sobolev spaces $(H^s \cap H^{s',1})(\R)$.
We attempt to push their arguments for local smoothing and maximal function estimates to the limit in terms of the generality of the dispersion relation.
Although the linear estimates we obtain are, to the best of our knowledge, not in the literature, we have relegated them to Appendix \ref{appendix:linear-estimates} due to their independence from the rest of the work.
Our linear estimates admit dispersion relations $\phi$ which are truly ``not like a derivative'' in the sense that critical points $\phi'(\xi) = 0$ with $\xi \neq 0$ are possible.

We use these linear estimates in an adaptation of the method in \cite{KenigPonceVega1994-2} to obtain a local well-posedness result for a large class of dispersive nonlinear systems of PDEs. This is the content of Theorem \ref{thmlwp}.
Since we need the smoothing effect to be compatible with Sobolev spaces, we must restrict ourselves here to dispersion relations which are rather ``like a derivative'', in the sense that $\phi'(\xi) = 0 \Longrightarrow \xi = 0$ and $\phi(\xi) \sim \langle \xi \rangle^n$ for large frequencies.
Afterwards, we prove a blow-up alternative (see Lemma \ref{lem:11}) and use Grönwall's inequality to prevent the $H^{s',1}$-norm from blowing up provided the $H^s$-norm does not (see Lemma \ref{lem:12}). This allows us to prove Theorem \ref{thm:1}.
Because of their technical nature, we refer the reader to Section \ref{section:4} for further details.

\subsection{Notations and definitions}

When we write $A \lesssim_{\text{things}} B$ we mean $A \leq C(\text{things}) B$, where the constant is a function of the objects in the parentheses,
and must be continuous in the real or complex parameters.

Whenever we consider a function space without explicit domain, the domain is implicitly assumed to be $\R$. 
When functions in the variables $(t, x)$ are considered, we use subscripts to denote with respect to which variable the function space should be considered.

Set $D_x = - i \del_x$ and note that $\ha{D_x f} = \xi \ha{f}$ and $\ha{x f} = - D_\xi \ha{f}$ for $f \in \mc{D}(\R; \C)$.

Let $m, K, L \in \N$ and $A, B$ be sets of formal complex-valued functions. We define the sets of formal functions
\begin{align*}
    A \cup \conj{A} &= \{f: f \in A \text{ or } \conj{f} \in A\} 
    & \del_x (B \cup \conj{B}) &= \{\del_x^n f: f \in B \text{ or } \conj{f} \in B, n \in \N\}
\end{align*} 
and
\begin{align} \label{eqn:O-def}
    \mc{O}^{m,L}_A(B) &= \left\{ \sum_{k=m}^K \; \underset{|l| \leq L}{\sum_{l = (l_1, \dots, l_k) \in \N^k}} \; \underset{\del_x^{l_1} h_1, \dots \del_x^{l_m} h_m \in \del_x (B \cup \conj{B})}{\sum_{h_1, \dots, h_k \in A \cup \conj{A}}} c_{h_1, \dots, h_k}^{k,l} \prod_{j=1}^k \del_x^{l_j} h_j: K \in \N, c_{h_1, \dots, h_k}^{k,l} \in \C \right\} \,.
\end{align}
We write $f = \mc{O}^{m,L}_A(B)$ if $f \in \mc{O}^{m,L}_A(B)$. If any of the parameters of $\mc{O}$ are missing, we set by default $m$, $L$ to zero and $A$, $B$ to the empty set.
This class contains all complex polynomial expressions in functions from $A \cup \conj{A}$ and derivatives thereof, for which each monomial has at most $L$ derivatives in total, and contains at least $m$ factors which are derivatives of a function in $B \cup \conj{B}$.

If $u(t)$ is a formal power series in $t$, then $[t^n] u(t)$ denotes the application of the linear functional that extracts the coefficient in front of $t^n$.

For matrices $M \in \C^{n \times n}$ we denote the $j$-th column by $M_j = (M_{j,k})_{1 \leq k \leq n}$, $1 \leq j \leq n$.
We denote by $|M|$ the operator norm and remark that the choice of underlying vector norm does not matter up to universal constants.
We define the Pauli matrix
\begin{align*}
    \sigma_3 = \begin{pmatrix} 1 & 0 \\ 0 & -1 \end{pmatrix} \,.
\end{align*}
For real numbers $x, y \in \R$, we define $x \land y = \min\{x, y\}$ and $x \lor y = \max\{x, y\}$, as well as the Japanese bracket $\langle x \rangle = \sqrt{x^2 + 4}$.

Lastly, we write $\mc{Q}_j$, $j \in \{1, 2, 3, 4\}$ for the four open quadrants of the complex plane. 

\textbf{Acknowledgements.} \\
We thank our colleagues Julia Henniger and Sebastian Ohrem for many insightful discussions. \\
This work was funded by the Deutsche Forschungsgemeinschaft (DFG, German Research Foundation) – Project-ID 258734477 – SFB 1173.

\section{Asymptotic expansions of direct scattering data and definitions of \texorpdfstring{$\mc{H}^\NLS_n$}{Hᴺᴸˢₙ} and \texorpdfstring{$\mc{H}^\GP_n$}{Hᴳᴾₙ}}
\label{section:2}

\subsection{A Riemann surface}

\begin{figure}[htbp]
    \caption{Mapping properties of $\lambda \leftrightarrow z$. 
    An arrow is bold if and only if the corresponding map is continuous from closed quadrant to closed quadrant.}

    \centering

    \newcommand*\xgap{3cm}   
    \newcommand*\ygap{1.8cm}   

    \definecolor{oiGrey}   {HTML}{444444} 
    \definecolor{oiBlue}   {HTML}{0072B2} 
    \definecolor{oiGreen}  {HTML}{009E73} 
    \definecolor{oiOrange} {HTML}{D55E00} 
    \definecolor{oiPurple} {HTML}{CC79A7} 

    \newcommand*\normalarrowwidth{.8pt}   
    \newcommand*\strongarrowwidth{1.6pt}  

    \newcommand*\LR[4]{
    \pgfmathsetmacro\outangle{#4}%
    \pgfmathsetmacro\inangle {180-#4}%
    \draw[#1,out=\outangle,in=\inangle] (#2) to (#3);%
    }
    \newcommand*\RL[4]{
    \pgfmathsetmacro\outangle{180-#4}%
    \pgfmathsetmacro\inangle {#4}%
    \draw[#1,out=\outangle,in=\inangle] (#2) to (#3);%
    }

    \newcommand*\SLoop[4]{%
    \pgfmathsetmacro\inangle{#3+30}
    \path  (#2) edge[#1,
                    loop,
                    out=#3, in=\inangle,
                    looseness=#4] (#2);
    }

    \begin{tikzpicture}[
        >=Stealth,
        base/.style    ={-{Latex[length=1.6mm,width=1mm]}, line width=\normalarrowwidth, oiGrey},
        f/.style    ={-{Latex[length=1.6mm,width=1mm]}, line width=\normalarrowwidth, oiBlue},
        finv/.style ={-{Latex[length=1.6mm,width=1mm]}, line width=\normalarrowwidth, oiGreen},
        g/.style    ={-{Latex[length=1.6mm,width=1mm]}, line width=\normalarrowwidth, oiOrange},
        ginv/.style ={-{Latex[length=1.6mm,width=1mm]}, line width=\normalarrowwidth, oiPurple},
        mydash/.style={dash pattern=on 3pt off 1.5pt},
        strong/.style={line width=\strongarrowwidth},
      ]
      
      \node (Q1) at (-\xgap,0)   {$\mc{Q}_1$};
      \node (Q3) at (-\xgap/2,0) {$\mc{Q}_3$};
      \node (Q2) at ( \xgap/2,0) {$\mc{Q}_2$};
      \node (Q4) at ( \xgap,0)   {$\mc{Q}_4$};
  
      \SLoop{g}              {Q1}{240}{7}
      \SLoop{ginv}           {Q1}{190}{7}
      \SLoop{f,strong}       {Q1}{140}{7}
      \SLoop{finv,strong}    {Q1}{ 90}{7}
  
      \SLoop{g,mydash}           {Q3}{60}{7}
      \SLoop{ginv,mydash}        {Q3}{ 10}{7}
      \SLoop{f,mydash,strong}    {Q3}{-40}{7}
      \SLoop{finv,mydash,strong} {Q3}{-90}{7}
  
      \SLoop{f,mydash}       {Q2}{240}{7}
      \SLoop{finv,mydash}    {Q2}{190}{7}
      \SLoop{g,strong}       {Q2}{140}{7}
      \SLoop{ginv,strong}    {Q2}{ 90}{7}
  
      \SLoop{f}                  {Q4}{60}{7}
      \SLoop{finv}               {Q4}{ 10}{7}
      \SLoop{g,mydash,strong}    {Q4}{-40}{7}
      \SLoop{ginv,mydash,strong} {Q4}{-90}{7}
  
      \RL{g}             {Q3}{Q1}{12}
      \RL{ginv}          {Q3}{Q1}{24}
      \RL{f,strong}      {Q3}{Q1}{36}
      \RL{finv,strong}   {Q3}{Q1}{48}
  
      \LR{g,mydash}             {Q1}{Q3}{-12}
      \LR{ginv,mydash}          {Q1}{Q3}{-24}
      \LR{f,mydash,strong}      {Q1}{Q3}{-36}
      \LR{finv,mydash,strong}   {Q1}{Q3}{-48}
  
      \RL{f,mydash}       {Q4}{Q2}{-12}
      \RL{finv,mydash}    {Q4}{Q2}{-24}
      \RL{g,strong}       {Q4}{Q2}{36}
      \RL{ginv,strong}    {Q4}{Q2}{48}
  
      \LR{f}                   {Q2}{Q4}{12}
      \LR{finv}                {Q2}{Q4}{24}
      \LR{g,mydash,strong}     {Q2}{Q4}{-36}
      \LR{ginv,mydash,strong}  {Q2}{Q4}{-48}
  
      \begin{scope}[
        shift={($(Q4.east)+(2cm,0.45cm)$)},
        font=\scriptsize,
        every node/.style={inner sep=0pt, outer sep=0pt},
      ]
    \matrix[
      matrix of nodes,
      row sep=2pt,
      column sep=0pt,
      nodes in empty cells,
      column 3/.style={nodes={minimum width=16pt}},
      anchor=west,
    ] at (0,-.4) {
      \tikz{\draw[base]        (0,0)--(.8,0);}  & $\in C(\mc{Q}_j; \mc{Q}_k)$ & &
      \tikz{\draw[base,strong] (0,0)--(.8,0);}  & $\in C(\overline{\mc{Q}_j}, \overline{\mc{Q}_k})$ \\
  
      \tikz{\draw[f]           (0,0)--(.8,0);}  & $\sqrt{\lambda^2-1}$  & &
      \tikz{\draw[f,mydash]    (0,0)--(.8,0);}  & $-\sqrt{\lambda^2-1}$ \\
  
      \tikz{\draw[finv]        (0,0)--(.8,0);}  & $\sqrt{z^2+1}$ & &
      \tikz{\draw[finv,mydash] (0,0)--(.8,0);}  & $-\sqrt{z^2+1}$ \\
  
      \tikz{\draw[g]           (0,0)--(.8,0);}  & $\sqrt[\sim]{\lambda^2-1}$ & &
      \tikz{\draw[g,mydash]    (0,0)--(.8,0);}  & $-\sqrt[\sim]{\lambda^2-1}$ \\
  
      \tikz{\draw[ginv]        (0,0)--(.8,0);}  & $\sqrt[\sim]{z^2+1}$ & &
      \tikz{\draw[ginv,mydash] (0,0)--(.8,0);}  & $-\sqrt[\sim]{z^2+1}$ \\
    };
  \end{scope}
  
    \end{tikzpicture}

    \label{fig:1}
\end{figure}
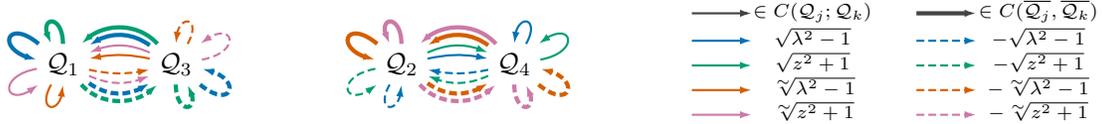

Consider the Riemann surface
\begin{align*}
    \K = \{(\lambda, z) \in \C^2: \lambda^2 - z^2 = 1\}
\end{align*}
and decompose $\K = \conj{\K_+} \cup \conj{\K_-}$ using the two sheets
\begin{align*}
    \K_\pm &= \{(\lambda, z) \in \K: \lambda \in \C \setminus ((-\infty, -1] \cup [1, \infty)), \pm \Imag z > 0\} \,.
\end{align*}
We define two branches of the square root:
\begin{align*}
    \sqrt{r e^{i \theta}} &= \sqrt{r} e^{i \frac{\theta}{2}} \qquad \forall\, \theta \in [0, 2 \pi) \quad \forall\, r > 0
    & \sqrt[\sim]{r e^{i \theta}} &= \sqrt{r} e^{i \frac{\theta}{2}} \qquad \forall\, \theta \in (- \pi, \pi] \quad \forall\, r > 0 \,.
\end{align*}
In figure \ref{fig:1} we depict pictographically the mapping properties between $\lambda$ and $z$ on the four open quadrants $\mc{Q}_j$, $j \in \{1, 2, 3, 4\}$ of the complex plane.
With the help of this pictogram, we can find on each closed quadrant an explicit homeomorphism that maps $\lambda \leftrightarrow z$. There is no single choice which works for all quadrants.
It is convenient to introduce the complex variable $\zeta$, which fulfills
\begin{align} \label{eqn:zeta}
    \zeta &= \lambda + z & \zeta^{-1} &= \lambda - z
    & \lambda &= \frac{\zeta + \zeta^{-1}}{2} & z &= \frac{\zeta - \zeta^{-1}}{2} \,.
\end{align}
When considering $(\lambda, z) \in \K$ and $\zeta \in \C$, we shall freely use the maps depicted in Figure \ref{fig:1}, as well as the relations \eqref{eqn:zeta}, to consider functions of one parameter also as functions of the others.
In particular, we may write $\lambda \in \K_\pm$ or $z \in \K_\pm$ to refer to a pair $(\lambda, z) \in \K_\pm$.
In situations where the details of the mapping matter, we split into cases depending on the quadrant in consideration.
Abusing notation, we define
\begin{align*}
    \gamma_1 &= 1 = \text{``$+$''}
    & \gamma_2 &= -1 = \text{``$-$''}
\end{align*}  
and write
\begin{align*}
    \K_{+\gamma_1} &= \K_{-\gamma_2} = \K_+
    & \K_{-\gamma_1} &= \K_{+\gamma_2} = \K_- \,.
\end{align*}

\subsection{Jost solutions}

Let $q \in C^\infty(\R; \C)$ have boundary values $\lim_{x \rightarrow \pm \infty} q(x) = q_\pm \in \SSS^1$ such that $q - q_\pm \in \mc{S}(\R_\pm; \C)$.
Consider for $j \in \{1, 2\}$ the solutions to the Zakharov-Shabat scattering problem
\begin{align} 
    \nonumber \Phi^\pm_j: \R \times \overline{\K_{\pm \gamma_j}} &\longrightarrow \C^2
    \\ \label{eqn:jost-1} \del_x \Phi^\pm_j(x, \lambda) &= \begin{pmatrix}
        - i \lambda & q(x)
        \\ \conj{q}(x) & i \lambda
    \end{pmatrix} \Phi^\pm_j(x, \lambda) 
    \\ \label{eqn:jost-2} \lim_{x \rightarrow \pm \infty} \Phi^\pm_j(x, \lambda) e^{x i z \gamma_j} &= E^\pm_j(\lambda) 
    \quad \text{ where } \quad E^\pm(\lambda) =
    \begin{pmatrix}
        q_\pm & i (z - \lambda)
        \\ i (\lambda - z) & \conj{q}_\pm
    \end{pmatrix} \,.
\end{align}
Here $E^\pm_1$ is the first and $E^\pm_2$ the second column of $E^\pm$, and the functions $E^\pm_j e^{- x i z \gamma_j}$ are solutions to \eqref{eqn:jost-1} if $q$ is replaced by $q_\pm$.
We call $\Phi^\pm_j$ the Jost solutions for the potential $q$. 
The modified Jost solutions $\Psi^\pm_j$ are defined by
\begin{align*}
    \Psi^\pm_j(x, \lambda) &= \Phi_j^\pm(x, \lambda) e^{x i z \gamma_j} \,,
\end{align*}
or equivalently as the unique solutions to the system
\begin{align} 
    \label{eqn:mjost-1} \del_x \Psi^\pm_j(x, \lambda) &= \begin{pmatrix}
        - i \lambda + i z \gamma_j & q(x)
        \\ \conj{q}(x) & i \lambda + i z \gamma_j
    \end{pmatrix} \Psi^\pm_j(x, \lambda) 
    \\ \label{eqn:mjost-2} \lim_{x \rightarrow \pm \infty} \Psi^\pm_j(x, \lambda) &= E^\pm_j(\lambda) \,.
\end{align}
Let us recall a well-known existence result for the Jost solutions. 
Note that we are interested in the behavior for large $|\lambda|$ and hence do not put emphasis on points of nonanalyticity.
\begin{lemma}[{\cite[Proposition 3]{DemontisPrinariMeeVitale}}]
    If $q - q_\pm \in L^1(\R_\pm)$ then the Jost solutions $\Phi^\pm_j$ and $\Psi^\pm_j$ exist. 
    They are analytic in $\lambda$ except for a finite set of points. They are smooth in $x$ if $q$ is smooth and all derivatives are in $L^1(\R)$.
\end{lemma}
This can be shown by usage of either of the following integral representations for the modified Jost solutions, and the corresponding Neumann series Ansatz:
\begin{align}
    \label{eqn:int-1} \Psi^\pm_j(x, \lambda) &= E^\pm_j(\lambda) + \int_{\pm \infty}^x e^{(x - y) A_\pm(\lambda)} (Q - Q_\pm)(y) \Psi^\pm_j(y, \lambda) e^{(x - y) i z \gamma_j} \dd y
    \\ \label{eqn:int-2} \Psi^\pm_j(x, \lambda) &= E^\pm_j(\lambda) + \int_{\pm \infty}^x e^{x \ti{A}(x, \lambda)} e^{- y \ti{A}(y, \lambda)} (Q - \ti{Q})(y) \Psi^\pm_j(y, \lambda) e^{(x - y) i z \gamma_j} \dd y \,.
\end{align}
Here
\begin{align*}
    Q(x) &= \begin{pmatrix}
        0 & q(x)
        \\ \conj{q}(x) & 0
    \end{pmatrix}
    & Q_\pm &= \begin{pmatrix}
        0 & q_\pm
        \\ \conj{q}_\pm & 0
    \end{pmatrix}
    & A_\pm &= Q_\pm - i \lambda \sigma_3
    \\ \sigma_3 &= \begin{pmatrix} 1 & 0 \\ 0 & - 1 \end{pmatrix} \,.
    & \ti{Q}(x) &= \mathds{1}_{\{x \leq 0\}} Q_- + \mathds{1}_{\{x > 0\}} Q_+ 
    & \ti{A}(x, \lambda) &= \ti{Q}(x) - i \lambda \sigma_3 \,.
\end{align*}
Observe that
\begin{align*}
    e^{x A_\pm(\lambda)} &= E^\pm(\lambda) e^{- x i z \sigma_3} E^\pm(\lambda)^{-1} \,.
\end{align*}
The representation \eqref{eqn:int-2} has the advantage that the integrand contains $Q - \ti{Q}$, which decays at \emph{both} infinities. 
As explained in the introduction, we are interested in integral representations of modified Jost solutions because we wish to study their asymptotic expansion in powers of $\lambda$ and $z$ at infinity.
Since we aim to control these asymptotic expansions at \emph{both} infinities, it is of importance to use an appropriate integral representation.
To obtain such an expansion, we may attempt to repeatedly integrate the factor $e^{(x - y) i z \gamma_j}$ in \eqref{eqn:int-2} using integration by parts.
Since the remaining integrand depends on nontrivially $\lambda$, this is inconvenient. 
We therefore need an integral representation that has an integrand with decay at both infinities, and is also suitable for asymptotic expansion via integration by parts.

\subsection{Triangular representations}

Another integral representation of the modified Jost solutions is given by the triangular representation (see \cite[Chapter 1, §8]{FaddeevTakhtajan}):
\begin{align} \label{eqn:int-3}
    \Psi^\pm_j(x, \lambda) &= E^\pm_j(\lambda) + \int_{\pm \infty}^x \Gamma^\pm(x, y) E^\pm_j(\lambda) e^{(x - y) i z \gamma_j} \dd y \,.
\end{align}
Here the kernel $\Gamma^\pm$ is independent of $\lambda$ and $E^\pm_j(\lambda)$ is bounded and analytic (up to a finite set of points) in $\lambda$.
Therefore, this representation is suitable for the asymptotic expansion of $\Psi^\pm_j(x, \lambda)$ in powers of $z$ at $\infty$.
In \cite{FaddeevTakhtajan} it is shown that $\Gamma^\pm \in \C^\infty(\{(x, y) \in \R^2: \pm (x - y) < 0\}; \C^{2 \times 2})$ exists and fulfills $\Gamma^\pm(x, -) \in \mc{S}(\{y \in \R: \pm (x - y) < 0\}; \C^{2 \times 2})$. 
Subsequently, the authors use this integral representation to derive various asymptotic expansions of the Jost solutions, but not the crucial expansion \eqref{eqn:exp-crucial} in Lemma \ref{lem:1} for the transmission coefficient that we are interested in.
Nevertheless, \eqref{eqn:exp-crucial} is proven in \cite{FaddeevTakhtajan}, but the proof is heavily abridged and works by referring to previous chapters, where the Zakharov-Shabat problem is considered on a Torus of length $L$, and then transferring expansions obtained in this setting through the limit $L \rightarrow \infty$. 
We present here a direct approach to the expansion of the Jost solutions and the transmission coefficient that uses the triangular representation \eqref{eqn:int-4} below.

Instead of using $q - \ti{q}$ as a source of decay at both infinities, we compare $q$ to a reference profile $q_\ast$ that assumes our boundary data $q_\pm$ at infinity.
We choose $q_\ast$ to be the dark soliton profile
\begin{align} \label{eqn:dark-soliton}
    q_\ast(x) &= q_+ \conj{\zeta_+} \left( \Real[\zeta_+] + i \Imag[\zeta_+] \tanh(\Imag[\zeta_+] x) \right)
    \\ &= q_- \conj{\zeta_-} \left( \Real[\zeta_-] + i \Imag[\zeta_-] \tanh(\Imag[\zeta_-] x) \right) \,.
\end{align}
Here $\zeta_+ \in e^{i [0, \pi)}$ and $\zeta_- = \conj{\zeta_+} \in e^{i [\pi, 2 \pi)}$ fulfill $\zeta_+^2 = \frac{q_+}{q_-}$ and $\zeta_-^2 = \frac{q_-}{q_+}$.
Most conveniently, explicit formulas are given in \cite{ChenChenHuang} for the modified Jost solutions of the dark soliton. They are
\begin{align}
    \label{eqn:1000} \Psi^-_{\ast,1}(x, \lambda) &= \begin{pmatrix}
        q_- \frac{\zeta - \zeta_+}{\zeta - \zeta_-} + i q_- \frac{1}{\zeta - \zeta_-} \frac{2 \Imag[\zeta_+]}{e^{2 \Imag[\zeta_+] x} + 1}
        \\ i \zeta^{-1} \frac{\zeta - \zeta_+}{\zeta - \zeta_-} \zeta_-^2 - \frac{1}{\zeta - \zeta_-} \zeta_- \frac{2 \Imag[\zeta_+]}{e^{2 \Imag[\zeta_+] x} + 1}
    \end{pmatrix}
    \;\, \Psi^-_{\ast,2}(x, \lambda) = \begin{pmatrix}
        - i \zeta^{-1} - \frac{1}{\zeta - \zeta_+} \zeta_+ \frac{2 \Imag[\zeta_-]}{e^{2 \Imag[\zeta_-] x} + 1}
        \\ \conj{q}_- - i \conj{q}_- \frac{1}{\zeta - \zeta_+} \frac{2 \Imag[\zeta_-]}{e^{2 \Imag[\zeta_-] x} + 1}
    \end{pmatrix}
    \\ \label{eqn:1002} \Psi^+_{\ast,1}(x, \lambda) &= \begin{pmatrix}
        q_+ \frac{\zeta - \zeta_-}{\zeta - \zeta_+} + i q_+ \frac{1}{\zeta - \zeta_+} \frac{2 \Imag[\zeta_-]}{e^{2 \Imag[\zeta_-] x} + 1}
        \\ i \zeta^{-1} \frac{\zeta - \zeta_-}{\zeta - \zeta_+} \zeta_+^2 - \frac{1}{\zeta - \zeta_+} \zeta_+ \frac{2 \Imag[\zeta_-]}{e^{2 \Imag[\zeta_-] x} + 1}
    \end{pmatrix}
    \;\, \Psi^+_{\ast,2}(x, \lambda) = \begin{pmatrix}
        - i \zeta^{-1} - \frac{1}{\zeta - \zeta_-} \zeta_- \frac{2 \Imag[\zeta_+]}{e^{2 \Imag[\zeta_+] x} + 1}
        \\ \conj{q}_+ - i \conj{q}_+ \frac{1}{\zeta - \zeta_-} \frac{2 \Imag[\zeta_+]}{e^{2 \Imag[\zeta_+] x} + 1}
    \end{pmatrix}
\end{align}
in our notation.
The subscript $\ast$ always denotes that the potential in the Zakharov-Shabat problem is $q_\ast$ instead of $q$.
We may now study $\Psi^\pm_j$ as perturbations of $\Psi^\pm_{\ast,j}$ by use of the integral representation formula
\begin{align} \label{eqn:int-4}
    \Psi^\pm_j(x, \lambda) &= \Psi^\pm_{\ast,j}(x, \lambda) + \int_{\pm \infty}^x \Gamma^\pm(x, y) \Psi^\pm_{\ast,j}(y, \lambda) e^{(x - y) i z \gamma_j} \dd y \,,
\end{align}
where crucially the kernel $\Gamma^\pm(x, y)$ is independent of $\lambda$, and $\Psi^\pm_{\ast,j}(y, \lambda)$ is bounded and analytic (up to a finite set of points) in $\lambda$, i.e. a benign factor.

For a matrix $M \in \C^{2 \times 2}$ we write $\overarcarrow{M}$ to denote the operator of multiplication from the right by $M$, and treat it as if it was a matrix in our use of language and notation.
Define 
\begin{align*}
    Q_\ast(x) &= \begin{pmatrix}
        0 & q_\ast(x) \\
        \conj{q}_\ast(x) & 0
    \end{pmatrix}
    & A_\ast(x, \lambda) &= Q_\ast(x) - i \lambda \sigma_3 \,,
\end{align*}
and recall that $\Psi^\pm_{\ast,j}$ solves
\begin{align} 
    \label{eqn:mjost-infty-1} \del_x \Psi^\pm_{\ast,j}(x, \lambda) &= (A_\ast(x, \lambda) + i z \gamma_j) \Psi^\pm_{\ast,j}(x, \lambda)
    \\ \label{eqn:mjost-infty-2} \lim_{x \rightarrow \pm \infty} \Psi^\pm_{\ast,j}(x, \lambda) &= E^\pm_j(\lambda) \,.
\end{align}
Substituting \eqref{eqn:int-4} into \eqref{eqn:mjost-1} yields
\begin{align*}
    0 &= (\del_x - A(x, \lambda) - i z \gamma_j) \Psi^\pm_j(x, \lambda)
    \\ &= - (A(x, \lambda) + i z \gamma_j - A_\ast(x, \lambda) - i z \gamma_j) \Psi^\pm_{\ast,j}(x, \lambda) + \Gamma^\pm(x, x) \Psi^\pm_{\ast,j}(x, \lambda)
    \\ &+ \int_{\pm \infty}^x (\del_x - A(x, \lambda)) \Gamma^\pm(x, y) \Psi^\pm_{\ast,j}(y, \lambda) e^{(x - y) i z \gamma_j} \dd y
    \\ &= \left( \Gamma^\pm(x, x) - (A - A_\ast)(x, \lambda) \right) \Psi^\pm_{\ast,j}(x, \lambda) + \int_{\pm \infty}^x (\del_x - A(x, \lambda)) \Gamma^\pm(x, y) \Psi^\pm_{\ast,j}(y, \lambda) e^{(x - y) i z \gamma_j} \dd y \,.
\end{align*}
We write
\begin{align*}
    \Gamma^\pm_\dg &= \frac12 (\Gamma^\pm + \sigma_3 \Gamma^\pm \sigma_3)
    & \Gamma^\pm_\odg &= \frac12 (\Gamma^\pm - \sigma_3 \Gamma^\pm \sigma_3)
    & \Gamma^\pm &= \Gamma^\pm_\dg + \Gamma^\pm_\odg = \sigma_3 \Gamma^\pm \sigma_3 + 2 \Gamma^\pm_\odg
\end{align*}
and prescribe
\begin{align*}
    \Gamma^\pm_\odg(x, x) &= \frac12 (A - A_\ast)_\odg(x, \lambda)  = \frac12 (Q - Q_\ast)(x) \,.
\end{align*}
Note also that
\begin{align*}
    (A - A_\ast)_\dg &= (Q - Q_\ast)_\dg = 0 \,.
\end{align*}
We require that
\begin{align*}
    0 &= \sigma_3 \Gamma^\pm(x, x) \sigma_3 \Psi^\pm_{\ast,j}(x, \lambda) + \int_{\pm \infty}^x (\del_x - A(x, \lambda)) \Gamma^\pm(x, y) \Psi^\pm_{\ast,j}(y, \lambda) e^{(x - y) i z \gamma_j} \dd y
\end{align*}
and consider this a special case of $F(x, y) = 0$, where
\begin{align*}
    F(x, y) &= \sigma_3 \Gamma^\pm(x, y) \sigma_3 \Phi^\pm_{\ast,j}(y, \lambda) + \int_{\pm \infty}^y (\del_x - A(x, \lambda)) \Gamma^\pm(x, s) \Phi^\pm_{\ast,j}(s, \lambda) \dd s \,.
\end{align*}
Observe that 
\begin{align*}
    \del_y F(x, y) &= \left( \sigma_3 \del_y \Gamma^\pm(x, y) \sigma_3 + \sigma_3 \Gamma^\pm(x, y) \sigma_3 A_\ast(y, \lambda) + (\del_x - A(x, \lambda)) \Gamma^\pm(x, y) \right) \Phi^\pm_{\ast,j}(y, \lambda)
    \\ \lim_{y \rightarrow \pm \infty} F(x, y) &= \lim_{y \rightarrow \pm \infty} \sigma_3 \Gamma^\pm(x, y) \sigma_3 \Phi^\pm_{\ast,j}(y, \lambda) \,.
\end{align*}
We see that $F(x, y) = 0$ is fulfilled if $\Gamma^\pm$ solves the boundary value problem
\begin{align*}
    \sigma_3 \del_y \Gamma^\pm(x, y) \sigma_3 = - \sigma_3 \Gamma^\pm(x, y) \sigma_3 A_\ast(y, \lambda) - (\del_x - A(x, \lambda)) \Gamma^\pm(x, y)
    \\ \lim_{y \rightarrow \pm \infty} \sigma_3 \Gamma^\pm(x, y) \sigma_3 \Phi^\pm_{\ast,j}(y, \lambda) = 0
    \qquad\qquad \Gamma^\pm_\odg(x, x) = \frac12 (A - A_\ast)_\odg(x, \lambda) \,.
\end{align*}
It is necessary at this point that $\Phi^\pm_{\ast,j}(y, \lambda)$ remains bounded as $y \rightarrow \pm \infty$.
Then a solution to the boundary value problem
\begin{align*}
    \sigma_3 \del_y \Gamma^\pm(x, y) \sigma_3 = - \sigma_3 \Gamma^\pm(x, y) \sigma_3 Q_\ast(y) - (\del_x - Q(x)) \Gamma^\pm(x, y)
    \\ \lim_{y \rightarrow \pm \infty} \Gamma^\pm(x, y) = 0 \quad\qquad \Gamma^\pm_\odg(x, x) = \frac12 (Q - Q_\ast)(x)
\end{align*}
indeed yields a solution to \eqref{eqn:mjost-1}, as long as \eqref{eqn:int-4} is well-defined.
In summary, we aim to construct $\Gamma^\pm_\dg$ and $\Gamma^\pm_\odg$ which solve the boundary value problem
\begin{align}
    \label{eqn:603} (\del_x + \del_y) \Gamma^\pm_\dg(x, y) &= \Big(Q(x) + \overarcarrow{Q_\ast(y)}\Big) \Gamma^\pm_\odg(x, y)
    \\ \label{eqn:604} (\del_x - \del_y) \Gamma^\pm_\odg(x, y) &= \Big(Q(x) - \overarcarrow{Q_\ast(y)}\Big) \Gamma^\pm_\dg(x, y)
    \\ \label{eqn:605} \lim_{y \rightarrow \pm \infty} \Gamma^\pm(x, y) &= 0 \qquad \qquad \Gamma^\pm_\odg(x, x) = \frac12 (Q - Q_\ast)_\odg(x) \,.
\end{align}
If the reference profile $q_\ast$ is chosen as $q_\ast = q_-$, then this system matches \cite[Chapter 1, (8.18)--(8.19)]{FaddeevTakhtajan}.

We state now our well-posedness result for $\Gamma^\pm$.
It provides all the estimates necessary for the aforementioned asymptotic expansion of the Jost solutions and the transmission coefficient in powers of $\lambda$ and $z$.
\begin{lemma} \label{lem:678}
    Let $q \in q_\ast + \mc{S}(\R; \C)$. There exist smooth solutions $\Gamma^\pm$ to \eqref{eqn:603}--\eqref{eqn:605}.
    For every $k, m \in \N$ there exist bounded, monotonic control functions $c_\pm \in C_b(\R; \R_+)$, decreasing faster at $\pm \infty$ than any power of $\frac{1}{x}$, such that
    \begin{align}
        \label{eqn:621} \|((\del_x + \del_y) \del_x^m \del_y^k \Gamma^-)(s + y, s - y)\|_{(L^1 \cap L^\infty)_s((- \infty, x])} &\lesssim c_-(x) e^{\int_{x-y}^{x+y} c_-(s) \dd s}
        \\ \label{eqn:621-alt} \|((\del_x + \del_y) \del_x^m \del_y^k \Gamma^+)(s - y, s + y)\|_{(L^1 \cap L^\infty)_s([x, \infty))} &\lesssim c_+(x) e^{\int_{x-y}^{x+y} c_+(s) \dd s}
    \end{align}
    for all $x \in \R$ and $y \geq 0$. Integrating along the direction $\del_x + \del_y$ yields
    \begin{align}
        \label{eqn:620} |\del_x^m \del_y^k \Gamma^-(x, y)| &\lesssim c_-\left( \frac{x+y}{2} \right) e^{\int_y^x c_-(s) \dd s}
        \\ \label{eqn:620-alt} |\del_x^m \del_y^k \Gamma^+(x, y)| &\lesssim c_+\left( \frac{x+y}{2} \right) e^{\int_x^y c_+(s) \dd s} \,.
    \end{align}
    On the diagonal $y = x$, we have 
    \begin{align}
        \label{eqn:619} (\del_y^k \Gamma^\pm)\vert_{y=x} &\in C_b(\R; \C^{2 \times 2}) \cap \mc{S}(\R_\pm; \C^{2 \times 2})
        & ((\del_x + \del_y) \del_y^k \Gamma^\pm)\vert_{y=x} &\in \mc{S}(\R; \C^{2 \times 2}) \,.
    \end{align}
\end{lemma}
\begin{proof}
    Appendix \ref{appendix:678} contains the proof of \eqref{eqn:621} and \eqref{eqn:621-alt}, as well as the statements on the diagonal (see Claim \ref{claim:1}).
\end{proof}

\subsection{The transmission coefficient, asymptotic expansions, and definitions of \texorpdfstring{$\mc{H}^\NLS_n$}{Hᴺᴸˢₙ} and \texorpdfstring{$\mc{H}^\GP_n$}{Hᴳᴾₙ}}
We focus only on the sheet $\conj{\K_+}$ here for simplicity, but every result in this section has an analogous statement and proof on $\conj{\K_-}$.

If $\Imag z = 0$ then $\Phi^- = (\Phi^-_1, \Phi^-_2)$ and $\Phi^+ = (\Phi^+_1, \Phi^+_2)$ are two fundamental solution matrices to \eqref{eqn:jost-1}, 
hence there exist $a(\lambda), b(\lambda) \in \C$ such that
\begin{align}
    \label{eqn:trans-1} \Phi^-_1(x, \lambda) &= a(\lambda) \Phi^+_1(x, \lambda) + b(\lambda) \Phi^+_2(x, \lambda) \,.
    \\ \label{eqn:trans-2} \Psi^-_1(x, \lambda) &= a(\lambda) \Psi^+_1(x, \lambda) + b(\lambda) e^{2 x i z} \Psi^+_2(x, \lambda) \,.
\end{align}
We call $a(\lambda)$ the Transmission coefficient of $q$. It fulfills
\begin{align} \label{eqn:trans-3}
    a(\lambda) &= \frac{\Psi^-_{1,1} \Psi^+_{2,2} - \Psi^-_{1,2} \Psi^+_{2,1}}{\Psi^+_{1,1} \Psi^+_{2,2} - \Psi^+_{1,2} \Psi^+_{2,1}} = \frac{\Psi^-_{1,1} \Psi^+_{2,2} - \Psi^-_{1,2} \Psi^+_{2,1}}{2 z (\lambda - z)} \,.
\end{align}
Using this formula, we extend $a(\lambda)$ analytically to $\lambda \in \K_+$ up to a finite set of points.
For the dark soliton $q_\ast$ we can explicitly determine $a_\ast(\lambda) = \frac{q_-}{q_+} \frac{\zeta - \zeta_+}{\zeta - \zeta_-}$ and $b_\ast(\lambda) = 0$.
Suppose for now that
\begin{align}
    \label{eqn:ass-limit} \text{the limits } \lim_{x \rightarrow \infty} \Psi^-_1(x, \lambda), \lim_{x \rightarrow - \infty} \Psi^+_2(x, \lambda) \text{ and } \lim_{x \rightarrow \infty} \del_x \Psi^-_1(x, \lambda), \lim_{x \rightarrow - \infty} \del_x \Psi^+_2(x, \lambda) \text{ exist.}
\end{align}
In particular, the latter two limits must be zero. 
Using \eqref{eqn:mjost-1} to substitute $\Psi^-_{1,2}$ in \eqref{eqn:trans-3} and taking the limit $x \rightarrow \infty$, or respectively substituting $\Psi^+_{2,1}$ and taking the limit $x \rightarrow - \infty$, we find that
\begin{align*}
    \lim_{x \rightarrow \infty} \Psi^-_{1,1}(x, \lambda) &= a(\lambda) q_+
    & \lim_{x \rightarrow - \infty} \Psi^+_{2,2}(x, \lambda) &= a(\lambda) \conj{q}_- \,.
\end{align*}
Together with the trivial limit \eqref{eqn:mjost-2}, this implies 
\begin{align} \label{eqn:trans-5}
    \log a(\lambda) &= \log q_- - \log q_+ + \int_{\R} \sigma^\GP(x, \lambda) \dd x 
    = \log q_- - \log q_+ - \int_{\R} \ti{\sigma}^\GP(x, \lambda) \dd x \,,
\end{align}
where
\begin{align}
    \label{eqn:sig-0} \sigma^\GP(x, \lambda) &= \frac{\del_x \Psi^-_{1,1}(x, \lambda)}{\Psi^-_{1,1}(x, \lambda)} 
    & \ti{\sigma}^\GP(x, \lambda) &= \frac{\del_x \Psi^+_{2,2}(x, \lambda)}{\Psi^+_{2,2}(x, \lambda)} \,.
\end{align}
In Lemma \ref{lem:1} below we establish that these densities are indeed integrable.

\subsubsection{Asymptotic expansions for \texorpdfstring{$\Psi^-_1$}{Ψ₁⁻}, \texorpdfstring{$\Psi^+_2$}{Ψ₂⁺}, 
\texorpdfstring{$\sigma^\GP$}{σᴳᴾ}, \texorpdfstring{$\tilde{\sigma}^\GP$}{σ~ᴳᴾ}, and \texorpdfstring{$a(\lambda)$}{a(λ)}}
\begin{definition}
    Let $D \subset \C$, $d \in \N$ and $1 \leq p \leq \infty$.
    \begin{enumerate}[(i)]
        \item We say that a function $f = f(z) \in \C^\infty(D; \C^d)$ has an 
        \textbf{asymptotic expansion in powers of $2 i z$ at infinity on $D$} if there exist $(f_n)_{n \in \N} \subset \C^d$ such that
        \begin{align*}
            \forall\, N \in \N \quad
            \lim_{|z| \rightarrow \infty} (2 i z)^N \left( f(z) - \sum_{n=0}^N \frac{f_n}{(2 i z)^n} \right) = 0 \,.
        \end{align*}
        We call $f_n$ the \textbf{expansion coefficients} of $f$ and note that they are unique.
        If $f$ and $g$ have such an expansion then $2 i z f$ and $f g$ do as well. 
        If $g(z) \neq 0$ for $|z|$ sufficiently large, then $f g^{-1}$ also has such an expansion.
        
        \item We say that a function $f = f(x, z) \in \C^\infty(\R \times D; \C^d)$ has an 
        \textbf{$L^p$-smooth asymptotic expansion in powers of $2 i z$ at infinity on $D$} if there exist $(f_n)_{n \in \N} \subset (C^\infty \cap L^p)(\R; \C^d)$ such that
        \begin{align*}
            \forall\, x \in \R \quad \forall\, k, N \in \N \quad
            \lim_{|z| \rightarrow \infty} \left\| (2 i z)^N \left( \del_x^k f(x, z) - \sum_{n=0}^N \frac{\del_x^k f_n(x)}{(2 i z)^n} \right) \right\|_{L^p_x(\R)} = 0 \,.
        \end{align*}
        In this case $\del_x f$ has such an expansion as well. 
        If $g$ has such an $L^\infty$-smooth expansion and $\inf_{x \in \R} |g(x, z)| > 0$ for $|z|$ sufficiently large, then $f g^{-1}$ also has an $L^p$-smooth asymptotic expansion.
        If $p = 1$ then $\int_{-\infty}^\infty f(x, z) \dd x$ has an asymptotic expansion in powers of $2 i z$ at infinity on $D$ with expansion coefficients $\int_{-\infty}^\infty f_n(x, z) \dd x$.
    \end{enumerate}
    When $f$ has an asymptotic expansion of some kind with expansion coefficients $(f_n)_{n \in \N}$, we write
    \begin{align*}
        f &\sim \sum_{n=0}^\infty \frac{f_n}{(2 i z)^n} \,.
    \end{align*}
\end{definition}
\begin{lemma} \label{lem:1}
    Assume $q \in C^\infty(\R; \C)$ with $q - q_\pm \in \mc{S}(\R_\pm; \C)$.
    Then $\Psi^-_1$, $\Psi^+_2$, $\del_x \Psi^-_1$, $\del_x \Psi^+_2$, $\sigma^\GP$, and $\ti{\sigma}^\GP$
    have $L^\infty$-smooth asymptotic expansions in powers of $2 i z$ at infinity on $\K_+ \cap \{\Imag z > c\}$ for some $c = c(q) > 0$.
    The asymptotic expansions for $\del_x \Psi^-_1$, $\del_x \Psi^+_2$, $\sigma^\GP$, and $\ti{\sigma}^\GP$ are $L^1$-smooth, and as a result $\log a(\lambda)$
    has an asymptotic expansion of the form
    \begin{align} \label{eqn:exp-crucial}
        \log a(\lambda) + \log q_+ - \log q_- &\sim \sum_{n=0}^\infty \frac{\int \sigma^\GP_n(x) \dd x}{(2 i z)^n}
        = - \sum_{n=0}^\infty \frac{\int \ti{\sigma}^\GP_n(x) \dd x}{(2 i z)^n} \,.
    \end{align}
    Furthermore, our assumption \eqref{eqn:ass-limit} is true.
\end{lemma}
\begin{proof}
    We give the proof only for the case of $\Psi^-_1$ and $\sigma^\GP$, as it is analogous for $\Psi^+_2$ and $\ti{\sigma}^\GP$.
    Recall the integral representation formula \eqref{eqn:int-4}. For $k \in \N$ it implies
    \begin{align}
        \label{eqn:int-4-deriv} \del_x^k \Psi^-_1(x, \lambda) 
        &= \del_x^k \Psi^-_{\ast,1}(x, \lambda) - \int_0^\infty \del_x^k \left( \Gamma^-(x, x - y) \Psi^-_{\ast,1}(x - y, \lambda) \right) e^{y i z} \dd y \,.
    \end{align}
    Integrating by parts $N$ times, we have
    \begin{align*}
        \del_x^k \Psi^-_1(x, \lambda) 
         &= \del_x^k \Psi^-_{\ast,1}(x, \lambda) - \sum_{n=1}^N \del_x^k \left( \del_y^{n-1}\big(\Gamma^-(x, y) \Psi^-_{\ast,1}(y, \lambda)\big)\big\vert_{y=x} \right) (i z)^{-n}
        \\ &- \int_0^\infty \del_x^k (- \del_y)^N \left( \Gamma^-(x, x - y) \Psi^-_{\ast,1}(x - y, \lambda) \right) e^{y i z} (i z)^{-N} \dd y \,.
    \end{align*}
    From the explicit formula \eqref{eqn:1000}, we know that $\Psi^-_{\ast,1}$ has an $L^\infty$-smooth and $\del_x \Psi^-_{\ast,1}$ an $L^1$-smooth asymptotic expansion in powers of $2 i z$ at infinity on $\K_+ \cap \{\Imag z > c\}$ for any $c = c(q) > 0$.
    We denote by $\Psi^-_{\ast,1,n}$ the expansion coefficients.
    Then we can write
    \begin{align*}
        \del_x^k \Psi^-_1(x, \lambda) 
        &= \sum_{n=0}^N (2 i z)^{-n} \left( \del_x^k \Psi^-_{\ast,1,n}(x) - \sum_{m=0}^{n-1} \del_x^k \left( 2^{n-m} \del_y^{n-1-m}\big(\Gamma^-(x, y) \Psi^-_{\ast,1,m}(y)\big)\big\vert_{y=x} \right) \right)
        \\  &+ \del_x^k \Psi^-_{\ast,1}(x, \lambda) - \sum_{n=0}^N \frac{\del_x^k \Psi^-_{\ast,1,n}(x)}{(2 i z)^n}
        \\  &- \sum_{n=1}^N (2 i z)^{-n} \del_x^k \left( 2^n \del_y^{n-1} \left( \Gamma^-(x, y) \left( \Psi^-_{\ast,1}(y, \lambda) - \sum_{m=0}^{N-n} \frac{\Psi^-_{\ast,1,m}(y)}{(2 i z)^m} \right) \right) \Bigg\vert_{y=x} \right)
        \\  &- \int_0^\infty \del_x^k (- \del_y)^N \left( \Gamma^-(x, x - y) \Psi^-_{\ast,1}(x - y, \lambda) \right) e^{y i z} (i z)^{-N} \dd y \,.
    \end{align*}
    By combining our knowledge of $\Gamma^-$ on the diagonal \eqref{eqn:619} with \eqref{eqn:1000}, we find that the intended expansion coefficients in the first line of the above expression are in $L^\infty$ if $k \geq 0$, and in $L^1$ if $k \geq 1$.
    Combining the estimates \eqref{eqn:621}--\eqref{eqn:620} with \eqref{eqn:1000}, we know that for every $k \geq 1$ there exists some $c \in L^\infty(\R; \R)$, decreasing faster than any power of $\frac{1}{x}$ at $\pm \infty$, such that
    \begin{align}
        \label{eqn:int-5} \left\| (- \del_y)^N \left( \Gamma^-(x, x - y) \Psi^-_{\ast,1}(x - y, \lambda) \right) e^{y i z} \right\|_{L^\infty_x((\pm \infty, x_0])} &\leq c\left(x_0-\frac{y}{2}\right) e^{y c(x_0) - y \Imag z}
        \\ \label{eqn:int-6} \left\| \del_x^k (- \del_y)^N \left( \Gamma^-(x, x - y) \Psi^-_{\ast,1}(x - y, \lambda) \right) e^{y i z} \right\|_{(L^1 \cap L^\infty)_x((\pm \infty, x_0])} &\leq c\left(x_0-\frac{y}{2}\right) e^{y c(x_0) - y \Imag z}
    \end{align}
    for all $y \geq 0$ and $\lambda$ with $\Imag z > \|c\|_{L^\infty}$.
    We find that $\Psi^-_1$ has an $L^\infty$-smooth and $\del_x \Psi^-_1$ a both $L^\infty$- and $L^1$-smooth asymptotic expansion in powers of $2 i z$ at infinity on $\K_+ \cap \{\Imag z > \|c\|_{L^\infty}\}$.
    Since $\inf_{x \in \R} |\Psi^-_{\ast,1}(x, \lambda)| > 0$ when $\Imag z$ is sufficiently large, we know that also $\sigma^\GP(x, \lambda)$
    has an $L^\infty$- and $L^1$-smooth asymptotic expansion in powers of $2 i z$ at infinity on $\K_+ \cap \{\Imag z > \|c\|_{L^\infty}\}$.
    It remains to verify \eqref{eqn:ass-limit}, which can be seen by applying dominated convergence to \eqref{eqn:int-4-deriv} with the estimates \eqref{eqn:int-5}--\eqref{eqn:int-6}.
\end{proof}
Note that from \eqref{eqn:int-4} and the same dominated convergence argument we obtain the integral representation formula
\begin{align*}
    a(\lambda) &= a_\ast(\lambda) \left( 1 - \int_0^\infty \lim_{x \rightarrow \infty} \Gamma^-(x, x - y) e^{y i z} \dd y \right) \,,
\end{align*}
which we present only for curiosity.
Applying \eqref{eqn:mjost-1} to \eqref{eqn:sig-0} yields the Riccati equations
\begin{align}
    \label{eqn:sig-1} (2 i z) \sigma^\GP - \sigma^\GP_x - (\sigma^\GP)^2 + q \conj{q} - 1 + q_x \frac{\sigma^\GP + i (\lambda - z)}{q} &= 0
    \\ \label{eqn:sigt-1} - (2 i z) \ti{\sigma}^\GP - \ti{\sigma}^\GP_x - (\ti{\sigma}^\GP)^2 + q \conj{q} - 1 + \conj{q}_x \frac{\ti{\sigma}^\GP - i (\lambda - z)}{\conj{q}} &= 0 \,.
\end{align} 
Since they are equivalent when $q$ and $\conj{q}$ are swapped and $i$ is replaced by $- i$, we have $\sigma^\GP_n = (-1)^n \conj{\ti{\sigma}^\GP_n}$. 
In particular
\begin{align} \label{eqn:sig-sigt-relation}
    \int \sigma_n^\GP(x, \lambda) \dd x &= (-1)^{n+1} \int \conj{\sigma_n^\GP(x, \lambda)} \dd x \,.
\end{align}

\subsubsection{Recurrence relations for the expansion coefficients.}
We can obtain a recurrence relation for the expansion coefficients $\sigma^\GP_n$ directly, by substituting $\sigma^\GP$ in \eqref{eqn:sig-1} with the formal power series $\sum_{n=0}^\infty \frac{\sigma^\GP_n}{(2 i z)^n}$ and comparing coefficients.
Due to the presence of $\lambda$ on the right this yields an awkward recurrence relation, so we consider instead the quantity
\begin{align*}
    \sigma^\NLS(x, \lambda) &= \sigma^\GP(x, \lambda) + i (\lambda - z) \,.
\end{align*}
Note that $\sigma^\NLS$ is not integrable in $x$. In the \eqref{eqn:NLS}--\eqref{eqn:ZBC} setting $\sigma^\NLS(x, \lambda)$ is indeed the density of the transmission coefficient, i.e. it fulfills \eqref{eqn:trans-5}, and it is $\sigma^\GP$ which is not integrable.
Equation \eqref{eqn:sig-1} is equivalent to
\begin{align} \label{eqn:sig-2}
    (2 i \lambda) \sigma^\NLS &= \del_x \sigma^\NLS - \frac{q_x}{q} \sigma^\NLS + (\sigma^\NLS)^2 - q \conj{q} \,.
\end{align}
By expanding
\begin{align}
    \label{eqn:sig-nls-expand} \sigma^\NLS(x, \lambda) &\sim \sum_{n=0}^\infty \frac{\sigma^\NLS_n(x)}{(2 i \lambda)^n} \,,
\end{align}
we obtain the recurrence relation
\begin{align}
    \label{eqn:sig-3} \sigma^\NLS_0 &= 0 
    & \sigma^\NLS_1 &= - q \conj{q}
    & \sigma^\NLS_{n+1} &= \del_x \sigma^\NLS_n - \frac{q_x}{q} \sigma^\NLS_n + \sum_{k=0}^n \sigma^\NLS_k \sigma^\NLS_{n-k} \,.
\end{align}

We can then derive the expansion coefficients $\sigma^\GP_n$ from $\sigma^\NLS_n$ and vice versa. 
This requires the use of a suitable map from Figure \ref{fig:1}, depending on the quadrants that $\lambda$ and $z$ are in.
For the quadrant $\conj{\mc{Q}_1}$, we choose $z = \sqrt{\lambda^2 - 1}$ and $\lambda = \sqrt{z^2 + 1}$, using the principal square root, and obtain the relations
\begin{align}
    \label{eqn:NLS-GP-1} \sigma^\NLS_{2m} &= \sum_{k=0}^m \binom{m-1}{m-k} (-4)^{m-k} \sigma^\GP_{2k}
    & \sigma^\NLS_{2m+1} &= \sum_{k=0}^m \binom{m-\tfrac{1}{2}}{m-k} (-4)^{m-k} \sigma^\GP_{2k+1} + (-1)^{m+1} C_m
    \\ \label{eqn:NLS-GP-2} \sigma^\GP_{2m} &= \sum_{k=0}^m \binom{m-1}{m-k} 4^{m-k} \sigma^\NLS_{2k}
    & \sigma^\GP_{2m+1} &= \sum_{k=0}^m \binom{m-\tfrac{1}{2}}{m-k} 4^{m-k} \sigma^\NLS_{2k+1} + C_m \,.
\end{align}
Here $C_n$ are the Catalan numbers. They can be defined either by the recurrence relation they solve, or by their generating function:
\begin{align}
    \label{eqn:catalan} C_0 &= 1 & C_{n+1} &= \sum_{k=0}^n C_k C_{n-k} 
    & \sum_{n=0}^\infty X^n C_n &= \frac{1 - \sqrt{1 - 4 X}}{2 X} \,.
\end{align}
On other quadrants, either the relations \eqref{eqn:NLS-GP-1}--\eqref{eqn:NLS-GP-2}, or alternatively the expansion coefficients in \eqref{eqn:sig-nls-expand}, must be corrected by sign changes.

\subsubsection{The Hamiltonians \texorpdfstring{$\mc{H}^\NLS_n$}{Hᴺᴸˢₙ} and \texorpdfstring{$\mc{H}^\GP_n$}{Hᴳᴾₙ}}

We define for $n \geq 0$ the Hamiltonians
\begin{align*}
    \mc{H}^\NLS_n &= - (- i)^n \int_{\R} \sigma^\NLS_{n+1}(x) \dd x
    & \mc{H}^\GP_n &= - (- i)^n \int _{\R}\sigma^\GP_{n+1}(x) \dd x \,.
\end{align*}

\begin{lemma} \label{lem:trans-comm}
    Assume the setting of Lemma \ref{lem:1}
    \begin{enumerate}[(i)]
        \item The Hamiltonians $\mc{H}^\NLS_n$ and $\mc{H}^\GP_n$ are real-valued functionals.
        \item For $n, m \in \N$ we have
        \begin{align*}
            \left\{\mc{H}^\GP_n, \mc{H}^\GP_m \right\} &= 0 \,.
            & \left\{\mc{H}^\NLS_n, \mc{H}^\NLS_m \right\} &= 0 \,.
        \end{align*}
        More generally, for all $\lambda_1, \lambda_2 \in \K$ we have $\{a(\lambda_1), a(\lambda_2)\} = 0$.
        \item We use the bijections $z = \sqrt{\lambda^2 - 1}$ and $\lambda = \sqrt{z^2 + 1}$ on the closed first quadrant $\conj{\mc{Q}_1}$.
        There exists some $c = c(q) > 0$ for which the functional $\log a(\lambda)$ has an asymptotic expansion on $\K_+ \cap \conj{\mc{Q}_1} \cap \{\Imag z > c\}$ in powers of $2 i z$ at infinity of the form
        \begin{align} \label{eqn:trans-expand-GP}
            \log a(\lambda) &\sim i \sum_{n=0}^\infty \frac{\mc{H}^\GP_n}{(2 z)^{n+1}} \,.
        \end{align}
        If in Lemma \ref{lem:1} and the surrounding theory the dark soliton $q_\ast$ is replaced by the trivial solution $q_\ast = 0$, i.e. we assume \ZBC, then we have instead the asymptotic expansion
        \begin{align} \label{eqn:trans-expand-NLS}
            \log a(\lambda) &\sim i \sum_{n=0}^\infty \frac{\mc{H}^\NLS_n}{(2 \lambda)^{n+1}} \,.
        \end{align}
    \end{enumerate}
\end{lemma}
\begin{proof}
    \begin{enumerate}[(i)]
        \item This is a consequence of \eqref{eqn:sig-sigt-relation}.
        \item We refer to \cite[III.§2]{FaddeevTakhtajan} and \cite[Theorem B.7]{KlausKochLiu2023}.
        \item This follows from Lemma \ref{lem:1}, subsequent elaboration, and the definition of the Hamiltonians.
    \end{enumerate}
\end{proof}

\section{Analysis of the structure of the \texorpdfstring{\NLS}{NLS} and \texorpdfstring{\GP}{GP} Hierarchies}
\label{section:3}

This section is concerned with extract structure from the recurrence relation \eqref{eqn:sig-3} in the form of explicit coefficients. 
For a function $F = F(q, \conj{q})$, we use the shorthand notation
\begin{align*}
    \delta F = \delta F(q, \conj{q}) = \frac{\delta}{\delta \conj{q}} \int F(q(x), \conj{q}(x)) \dd x
\end{align*}
for the functional derivative with respect to $\conj{q}$. With this notation, we write \eqref{eqn:NLS-n} as
\begin{align*}
    \del_{t_n} q &= \frac{\delta \mc{H}^\NLS_n}{\delta \conj{q}} = - (- i)^n \delta \sigma^\NLS_{n+1} \,.
\end{align*}
From the rest of this section we set $\sigma = \sigma^\NLS$.
Given a polynomial $Q$ of $q, \conj{q}$ and their derivatives, we write $\pi_k Q$ for the sum of all monomials in $Q$ which have exactly $k$ factors with derivatives.
We want to obtain an explicit formula for $\pi_1 \delta \sigma_n$. 
The reason is that in our well-posedness theory we make a perturbative Ansatz $q(x, t) = q_\ast(x) + p(x, t)$ and the non-trivial linear part of the equation for $p$ depends on $\delta \pi_0 \sigma_n$ and $\delta \pi_1 \sigma_n$.
The objective of our analysis is to find an explicit formula for this linear part.
Note that
\begin{align*}
    \delta \sigma_n &= \pi_0 \delta \sigma_n + \pi_1 \delta \sigma_n + \mc{O}^{2,n-3}_q(q_x) \,.
\end{align*}
The only information we need and use is the recurrence relation \eqref{eqn:sig-3}.
We start by noting that it implies
\begin{align}
    \nonumber \pi_0 \sigma_1 &= - |q|^2 \\
    \label{eqn:pi-0} \pi_0 \sigma_{n+1} &= \sum_{k=1}^n \pi_0 \sigma_k \pi_0 \sigma_{n-k} \,.
\end{align}
This allows us to find an explicit formula for $\pi_0 \sigma_n$.
Next, we observe that \eqref{eqn:sig-3} implies the following recurrence relation for $\pi_1 \sigma_n$, involving $\pi_0 \sigma_n$:
\begin{align}
    \nonumber \pi_1 \sigma_1 &= 0 \\
    \label{eqn:pi-1} \pi_1 \sigma_{n+1} &= \del_x \pi_0 \sigma_n - \frac{q_x}{q} \pi_0 \sigma_n + \pi_1 \del_x \pi_1 \sigma_n + \sum_{k=1}^n 2 \pi_1 \sigma_k \pi_0 \sigma_{n-k} \,.
\end{align}
This allows us to find an explicit formula for $\pi_1 \sigma_n$. 
Unfortunately, $\pi_1 \delta \sigma_n$ depends not just on $\delta \pi_1 \sigma_n$, but also certain terms from $\pi_2 \sigma_n$.
Let us carefully perform this analysis, starting with the trivial observation
\begin{align*}
    \pi_1 \delta \sigma_n &= \sum_{k=0}^\infty \pi_1 \delta \pi_k \sigma_n \,.
\end{align*}
Clearly, $\delta \pi_0 \sigma_n$ has no derivatives and is therefore in the kernel of $\pi_1$. 
Similarly, $\delta \pi_k \sigma_n$ for $k \geq 3$ always have at least $2$ factors in each monomial which have derivatives, so it is also in the kernel of $\pi_1$. 
We suppose here that we already have an explicit formula for $\pi_1 \delta \pi_1 \sigma_n$, so it remains to study $\pi_1 \delta \pi_2 \sigma_n$.
In fact, we have $\pi_1 \delta \pi_2 \sigma_n = \pi_1 \delta \ti{\pi}_2 \sigma_n$, where $\ti{\pi}_2$ projects onto sums of monomials which have exactly two factors with derivatives, and one of those factors is of the form $\del_x^k \conj{q}$. 
Similarly, we define $\ti{\pi}_1$ to project onto sums of monomials which have only one factor with derivatives, and it is of the form $\del_x^k \conj{q}$.
We obtain the decomposition
\begin{align}
    \label{eqn:pi-1-delta} \pi_1 \delta \sigma_n &= \pi_1 \delta \pi_1 \sigma_n + \pi_1 \delta \ti{\pi}_2 \sigma_n \,.
\end{align}
With our explicit formula for $\pi_1 \sigma_n$, we can compute the first term, so it remains to find an explicit formula for $\ti{\pi}_2 \sigma_n$. 
The corresponding reccurence relation is
\begin{align}
    \nonumber \ti{\pi}_2 \sigma_1 &= 0 \\
    \label{eqn:pi-2} \ti{\pi}_2 \sigma_{n+1} &= \ti{\pi}_2 \del_x \ti{\pi}_2 \sigma_n - \frac{q_x}{q} \ti{\pi}_1 \sigma_n + \ti{\pi}_2 \del_x \pi_1 \sigma_n + 2 \sum_{k=1}^{n-1} (\ti{\pi}_2 \sigma_k \pi_0 \sigma_{n-k} + \ti{\pi}_1 \sigma_k \pi_1 \sigma_{n-k}) \,.
\end{align}

\subsection{Formulas for \texorpdfstring{$\pi_0 \sigma_n$}{π₀σₙ} and \texorpdfstring{$\pi_1 \sigma_n$}{π₁σₙ}}

\begin{lemma}
    We have
    \begin{align}
        \label{eqn:pi-0-formula} \pi_0 \sigma_n &= \mathds{1}_{\{n \text{ odd}\}} C_{\frac{n-1}{2}} \conj{q}^{\frac{n+1}{2}} (- q)^{\frac{n+1}{2}} \,,
    \end{align}
    where $C_n$ are the Catalan numbers defined in \eqref{eqn:catalan}.
\end{lemma}
\begin{proof}
    We prove by induction. The base case is trivial, so we assume the formula holds for $n$. Then
    \begin{align*}
        \pi_0 \sigma_{n+1} &= \sum_{k=1}^n \pi_0 \sigma_k \pi_0 \sigma_{n-k} \\
        &= \sum_{k=1}^n \mathds{1}_{\{k \text{ odd}\}} C_{\frac{k-1}{2}} \conj{q}^{\frac{k+1}{2}} (- q)^{\frac{k+1}{2}}  \mathds{1}_{\{n-k \text{ odd}\}} C_{\frac{n-k-1}{2}} \conj{q}^{\frac{n-k+1}{2}} (- q)^{\frac{n-k+1}{2}}  \\
        &= \mathds{1}_{\{n+1 \text{ odd}\}} \conj{q}^{\frac{n+2}{2}} (- q)^{\frac{n+2}{2}} \sum_{k=1}^n \mathds{1}_{\{k \text{ odd}\}} C_{\frac{k-1}{2}} C_{\frac{n-k-1}{2}} \,,
    \end{align*}
    so by the definition of the Catalan numbers \eqref{eqn:catalan}, we are done.
\end{proof}
Using \eqref{eqn:pi-0-formula}, the iteration \eqref{eqn:pi-1} simplifies to
\begin{align*}
    \pi_1 \sigma_{n+1} &= \del_x (\mathds{1}_{\{n \text{ odd}\}} C_{\frac{n-1}{2}} \conj{q}^{\frac{n+1}{2}} (- q)^{\frac{n+1}{2}}) + \pi_1 \del_x \pi_1 \sigma_n + q_x \mathds{1}_{\{n \text{ odd}\}} C_{\frac{n-1}{2}} \conj{q}^{\frac{n+1}{2}} (- q)^{\frac{n-1}{2}} 
    \\ &+ \sum_{k=1}^n 2 \pi_1 \sigma_k \mathds{1}_{\{n-k \text{ odd}\}} C_{\frac{n-k-1}{2}} \conj{q}^{\frac{n-k+1}{2}} (- q)^{\frac{n-k+1}{2}} \,.
\end{align*}
\begin{lemma}
    We have the formula
    \begin{align}
        \label{eqn:pi-1-formula} \pi_1 \sigma_n &= \sum_{j=0}^{\floor{\frac{n}{2}} - 1} D_{n,j} (- q)^{j+1} \conj{q}^j \del_x^{n - 1 - 2 j} \conj{q} + \sum_{j=0}^{\floor{\frac{n}{2}} - 2} E_{n,j} (- q)^{j+1} \conj{q}^{j+2} \del_x^{n - 3 - 2 j} q \,,
    \end{align}
    where
    \begin{align}
        \label{eqn:D-formula} D_{n,j} &= 4^j \begin{pmatrix}
            \frac{n}{2} - 1 \\ j
        \end{pmatrix} \\
        \label{eqn:E-formula} E_{n,j} &= \sum_{l=0}^j (-1)^l (C_{l+1} - 2 C_l) 4^{j-l} \begin{pmatrix}
            \frac{n - 1}{2} - 1 \\ j - l
        \end{pmatrix} \,.
    \end{align}
\end{lemma}

\begin{proof}
    Recall that
    \begin{align*}
        \pi_1 \sigma_{n+1} &= \del_x (\mathds{1}_{\{n \text{ odd}\}} C_{\frac{n-1}{2}} \conj{q}^{\frac{n+1}{2}} (- q)^{\frac{n+1}{2}}) + \pi_1 \del_x \pi_1 \sigma_n + \sum_{k=1}^n 2 \pi_1 \sigma_k \mathds{1}_{\{n-k \text{ odd}\}} C_{\frac{n-k-1}{2}} \conj{q}^{\frac{n-k+1}{2}} (- q)^{\frac{n-k+1}{2}} \\
        &+ q_x \mathds{1}_{\{n \text{ odd}\}} C_{\frac{n-1}{2}} \conj{q}^{\frac{n+1}{2}} (- q)^{\frac{n-1}{2}} \,.
    \end{align*}
    We plug \eqref{eqn:pi-1-formula} as an Ansatz into this recurrence relation. 
    After a lengthy calculation, we obtain the following recurrence relation for the coefficients:
    \begin{align*}
        D_{2,0} &= 1 &  
        D_{n+1,j} &= 2 \sum_{k=0}^{2j-1} \mathds{1}_{\{k \text{ odd}\}} C_{\frac{k-1}{2}} D_{n-k,j-\frac{k+1}{2}} + \begin{cases}
            C_{\frac{n-1}{2}} \frac{n+1}{2} &, j = \frac{n-1}{2} \\
            D_{n,j} &, \text{ else }
        \end{cases} & \forall\, 0 \leq j \leq \left\lfloor \frac{n+1}{2} \right\rfloor - 1 \\
        E_{4,0} &= - 1 & 
        E_{n+1,j} &= 2 \sum_{k=0}^{2j-1} \mathds{1}_{\{k \text{ odd}\}} C_{\frac{k-1}{2}} E_{n-k,j-\frac{k+1}{2}} + \begin{cases}
            - C_{\frac{n-1}{2}} \frac{n-1}{2} &, j = \frac{n-3}{2} \\
            E_{n,j} &, \text{ else }
        \end{cases} & \forall\, 0 \leq j \leq \left\lfloor \frac{n+1}{2} \right\rfloor - 2
    \end{align*}
    In particular,
    \begin{align*}
        D_{2,0} &= 1 &  
        D_{n+1,j} &= 2 \sum_{k=0}^{j-1} C_k D_{n-2k-1,j-1-k} + \begin{cases}
            C_{\frac{n-1}{2}} \frac{n+1}{2} &, j = \frac{n-1}{2} \\
            D_{n,j} &, \text{ else }
        \end{cases} & \forall\, 0 \leq j \leq \left\lfloor \frac{n+1}{2} \right\rfloor - 1 \\
        E_{4,0} &= - 1 & E_{n+1,j} &= 2 \sum_{k=0}^{j-1} C_k E_{n-2k-1k,j-1-k} + \begin{cases}
            - C_{\frac{n-1}{2}} \frac{n-1}{2} &, j = \frac{n-3}{2} \\
            E_{n,j} &, \text{ else }
        \end{cases} & \forall\, 0 \leq j \leq \left\lfloor \frac{n+1}{2} \right\rfloor - 2 \,.
    \end{align*}
    When writing down such recurrence relations, we always set the coefficients which are not explicitly defined to zero.
    We now verify that \eqref{eqn:D-formula}--\eqref{eqn:E-formula} solve this recurrence relation.
    Note that
    \begin{align*}
        C_m (m + 1) &= \begin{pmatrix}
            2 m \\ m
        \end{pmatrix} = (- 4)^m \begin{pmatrix}
            - \frac12 \\ m
        \end{pmatrix} = 4^m \begin{pmatrix}
            \frac12 + m - 1 \\ m
        \end{pmatrix} \,.
    \end{align*}
    Since solutions to the recurrence are unique, it suffices to show that for $0 \leq j \leq \floor{\frac{n+1}{2}} - 1$ we have
    \begin{align*}
        4^j \begin{pmatrix}
            \frac{n+1}{2} - 1 \\ j
        \end{pmatrix} &= 2 \sum_{k=0}^{j-1} C_k 4^{j - 1 - k} \begin{pmatrix}
            \frac{n-1}{2} - 1 - k \\ j - 1 - k
        \end{pmatrix} + \begin{cases}
            4^j \begin{pmatrix}
                \frac{n}{2} - 1 \\ j
            \end{pmatrix} &, j = \frac{n-1}{2} \\
            4^j \begin{pmatrix}
            \frac{n}{2} - 1 \\ j
        \end{pmatrix} &, 0 \leq j \leq \floor{\frac{n}{2}} - 1 
        \end{cases} \\
        \iff 4^j \begin{pmatrix}
            \frac{n+1}{2} - 1 \\ j
        \end{pmatrix} &= 2 \sum_{k=0}^{j-1} C_k 4^{j - 1 - k} \begin{pmatrix}
            \frac{n-1}{2} - 1 - k \\ j - 1 - k
        \end{pmatrix} + 4^j \begin{pmatrix}
            \frac{n}{2} - 1 \\ j
        \end{pmatrix} \\
        \iff \begin{pmatrix}
            \frac{n+1}{2} - 1 \\ j + 1
        \end{pmatrix} - \begin{pmatrix}
            \frac{n}{2} - 1 \\ j + 1
        \end{pmatrix} &= \frac12 \sum_{k=0}^j 4^{-k} C_k \begin{pmatrix}
            \frac{n-1}{2} - 1 - k \\ j - k
        \end{pmatrix} \qquad \qquad \forall\, -1 \leq j \leq \floor{\frac{n+1}{2}} \\
        \iff \begin{pmatrix}
            \frac{s}{2} + 1 \\ j + 1
        \end{pmatrix} - \begin{pmatrix}
            \frac{s}{2} + \frac12 \\ j + 1
        \end{pmatrix} &= \frac12 \sum_{k=0}^j 4^{-k} C_k \begin{pmatrix}
            \frac{s}{2} - k \\ j - k
        \end{pmatrix} \qquad \qquad \text{ where } s = n-3 \,.
    \end{align*}
    This identity is verified by tedious calculation and subsequent comparison of the generating functions
    \begin{align*}
        \sum_{s, j \geq 0} X^s Y^j \begin{pmatrix}
            \frac{s}{2} + 1 \\ j + 1
        \end{pmatrix} - \sum_{s, j \geq 0} X^s Y^j \begin{pmatrix}
            \frac{s}{2} + \frac12 \\ j + 1
        \end{pmatrix} &= \frac{1 + Y - \sqrt{1 + Y}}{Y} \frac{1}{1 - X \sqrt{1 + Y}}
    \end{align*}
    and
    \begin{align*}
        \sum_{s, j \geq 0} X^s Y^j \frac12 \sum_{k=0}^j 4^{-k} C_k \begin{pmatrix}
            \frac{s}{2} - k \\ j - k
        \end{pmatrix} &= \frac{1 + Y - \sqrt{1 + Y}}{Y} \frac{1}{1 - X \sqrt{1 + Y}} \,.
    \end{align*}
    We proceed similarly for the coefficients $E_{n,j}$. 
    Here it suffices to show that for $0 \leq j \leq \floor{\frac{n+1}{2}} - 2$, we have
    \begin{align*}
        \sum_{l=0}^j (-1)^l (C_{l+1} - 2 C_l) 4^{-l} \begin{pmatrix}
            \frac{n}{2} - 1 \\ j - l
        \end{pmatrix} &= 2 \sum_{k=0}^{j-1} C_k \sum_{l=0}^{j-1-k} (-1)^l (C_{l+1} - 2 C_l) 4^{-1-k-l} \begin{pmatrix}
            \frac{n}{2} - 2 - k \\ j - 1 - k - l
        \end{pmatrix} \\
        &+ \begin{cases}
                - 4^{-j} C_{\frac{n-1}{2}} \frac{n-1}{2} &, j = \frac{n-3}{2} \\
                \sum_{l=0}^j (-1)^l (C_{l+1} - 2 C_l) 4^{-l} \begin{pmatrix}
                    \frac{n-1}{2} - 1 \\ j - l
                \end{pmatrix} &, 0 \leq j \leq \floor{\frac{n}{2}} - 2 
            \end{cases} \,.
    \end{align*}
    Further tedious calculations reveal
    \begin{align*}
        \sum_{m=0}^\infty Y^m \sum_{l=0}^m (-1)^l (C_{l+1} - 2 C_l) 4^{-l} 
        \begin{pmatrix}
            m \\ m - l
        \end{pmatrix} &= \frac{8}{Y^2} - \frac{4 \sqrt{1 - Y}}{Y^2} - \frac{4}{Y^2 \sqrt{1 - Y}}
    \end{align*}
    and
    \begin{align*}
        - \sum_{m=0}^\infty 4^{-m} Y^m C_{m+1} (m+1) &= \frac{8}{Y^2} - 4 \frac{\sqrt{1 - Y}}{Y^2} - \frac{4}{Y^2 \sqrt{1 - Y}} \,,
    \end{align*}
    so it remains to show that
    \begin{align*}
        & \sum_{l=0}^j (-1)^l (C_{l+1} - 2 C_l) 4^{-l} \left( \begin{pmatrix}
            \frac{n}{2} - 1 \\ j - l
        \end{pmatrix}  - \begin{pmatrix}
            \frac{n-1}{2} - 1 \\ j - l
        \end{pmatrix} \right)
        \\ &= 2 \sum_{k=0}^{j-1} C_k \sum_{l=0}^{j-1-k} (-1)^l (C_{l+1} - 2 C_l) 4^{-1-k-l} \begin{pmatrix}
            \frac{n}{2} - 2 - k \\ j - 1 - k - l
        \end{pmatrix} \,.
    \end{align*}
    This follows by yet another lengthy calculation of the generating functions
    \begin{align*}
        & \sum_{n, j \geq 0} X^n Y^j \sum_{l=0}^j (C_{l+1} - 2 C_l) 4^{-l} \left( \begin{pmatrix}
            \frac{n}{2} - 1 \\ j - l
        \end{pmatrix} - \begin{pmatrix}
            \frac{n-1}{2} - 1 \\ j - l
        \end{pmatrix} \right) \\
        &= \left( 2 \left( \frac{4}{Y} - 2 \right) \frac{1 - \sqrt{1 - Y}}{Y} - \frac{4}{Y} \right) \frac{1}{1 + Y} \frac{1}{1 - X \sqrt{1 + Y}} \left( 1 - \frac{1}{\sqrt{1 + Y}} \right) 
    \end{align*}
    and
    \begin{align*}
        & \sum_{n, j \geq 0} X^n Y^j 2 \sum_{k=0}^{j-1} C_k \sum_{l=0}^{j-1-k} (C_{l+1} - 2 C_l) 4^{-1-k-l} \begin{pmatrix}
            \frac{n}{2} - 2 - k \\ j - 1 - k - l
        \end{pmatrix} \\
        &= \left( 2 \left( \frac{4}{Y} - 2 \right) \frac{1 - \sqrt{1 - Y}}{Y} - \frac{4}{Y} \right) \frac{1}{1 + Y} \frac{1}{1 - X \sqrt{1 + Y}} \left( 1 - \frac{1}{\sqrt{1 + Y}} \right) \,.
    \end{align*}
\end{proof}

\subsection{Formulas for \texorpdfstring{$\ti{\pi}_2 \sigma_n$}{π̃₂ σₙ} and \texorpdfstring{$\pi_1 \delta \ti{\pi}_2 \sigma_n$}{π₁ δ π̃₂ σₙ}}

\begin{lemma}
    We have the formula
    \begin{align}
        \label{eqn:pi-2-formula} \pi_1 \delta \ti{\pi}_2 \sigma_n &= \sum_{j=0}^{\floor{\frac{n-1}{2}} - 2} (- q)^{j+2} \conj{q}^j \del_x^{n-3-2j} \conj{q} \ti{F}_{n,j} + \sum_{j=0}^{\floor{\frac{n-1}{2}} - 1} (- q)^j \conj{q}^j \del_x^{n-1-2j} q \ti{G}_{n,j} \,,
    \end{align}
    where
    \begin{align}
        \label{F-formula-1}\ti{F}_{2m} &= 0 \\
        \label{F-formula-2}\ti{F}_{2m+1,j} &= - (8 m - 8 j - 6) 4^j \binom{m - \frac12}{j} \\
        \label{G-formula-1}\ti{G}_{2m,j} &= \sum_{k=0}^{j-1} C_k 4^{j-k} (m-1-j) \binom{m-1-k-\frac12}{j-1-k} \\
        \label{G-formula-2}\ti{G}_{2m+1,j} &= 4^j (2j+1) \binom{m-\frac12}{j} - 4^j \frac12 \binom{m-\frac12}{j-1} + \sum_{k=0}^j C_k 4^{j-k} (j+1) \binom{m-1-k}{j-k} \,.
    \end{align}
\end{lemma}

\begin{proof}
    Recall the recurrence relation that $\ti{\pi}_2$ solves:
    \begin{align*}
        \ti{\pi}_2 \sigma_1 &= 0 \\
        \ti{\pi}_2 \sigma_{n+1} &= \ti{\pi}_2 \del_x \ti{\pi}_2 \sigma_n - \frac{q_x}{q} \ti{\pi}_1 \sigma_n  + \ti{\pi}_2 \del_x \pi_1 \sigma_n + 2 \sum_{k=1}^{n-1} (\ti{\pi}_2 \sigma_k \pi_0 \sigma_{n-k} + \ti{\pi}_1 \sigma_k \pi_1 \sigma_{n-k}) \,.
    \end{align*}
    We make the Ansatz
    \begin{align} \label{eqn:tipi-2-sig}
        \ti{\pi}_2 \sigma_n &= \sum_{j=1}^{\floor{\frac{n-1}{2}} - 1} (- q)^j \conj{q}^j \sum_{t=1}^{n-2-2j} G_{n,j,t} \del_x^{n-1-2j-t} q \del_x^t \conj{q} 
        \\ &+ \sum_{j=0}^{\floor{\frac{n-1}{2}} - 2} (- q)^{j+2} \conj{q}^j \sum_{t=1}^{\floor{\frac{n-1}{2}}-1-j} F_{n,j,t} \del_x^{n-3-2j-t} \conj{q} \del_x^t \conj{q} \,.
    \end{align}
    Unfortunately, we must now plug this Ansatz, as well as the formulas for $\pi_0 \sigma_n$ and $\pi_1 \sigma_n$ that we have obtained, into the recurrence relation.
    A lengthy calculation, which we spare the reader, yields the following recurrence relations for $G_{n,j,t}$ and $F_{n,j,t}$.
    We set $G_{5,1,1} = -6$ and for $n \geq 0$, $1 \leq j \leq \floor{\frac{n}{2}} - 1$, $1 \leq t \leq n - 1 -2 j$ we have
    \begin{align*}
        G_{n+1,j,t} &= 2 \sum_{k=0}^{j-2} G_{n-2k-1,j-k-1,t} C_k + 2 \sum_{k=0}^{j-2} D_{t+1+2k,k} E_{n-1-2k-t,j-2-k} \\
        &+ \mathds{1}_{\{j \neq  \frac{n}{2}-1\}} \mathds{1}_{\{t \neq n - 1 - 2j\}} G_{n,j,t} + \mathds{1}_{\{j \neq  \frac{n}{2}-1\}} \mathds{1}_{\{t \neq 1\}} G_{n,j,t-1} \\
        &+ \mathds{1}_{\{t = 1\}} (j+1) E_{n,j-1} - \mathds{1}_{\{t = n - 1 - 2j\}} j D_{n,j} \,.
    \end{align*}
    Similarly, we set $F_{5,0,1} = 5$ and for $n \geq 5$, $0 \leq j \leq \floor{\frac{n}{2}} - 2$, $1 \leq t \leq \floor{\frac{n}{2}} - 1 - j$ we have
    \begin{align*}
        F_{n+1,j,t} &= 2 \sum_{k=0}^{j-1} F_{n-2k-1,j-k-1,t} C_k + (1 + \mathds{1}_{\{t \neq \frac{n}{2}-1-j\}}) \sum_{k=0}^j D_{t+1+2k,k} D_{n-1-2k-t,j-k} \\
        &+ \mathds{1}_{\{t = 1\}} (j + 1) D_{n,j+1} + \mathds{1}_{\{j \neq \frac{n}{2} - 2\}} \Big( \mathds{1}_{\{t \neq \frac{n}{2} - 1 - j\}} F_{n,j,t} + \mathds{1}_{\{t \neq 1\}} F_{n,j,t-1} + \mathds{1}_{\{t = \frac{n-1}{2} - 1 - j\}} F_{n,j,t} \Big) \,.
    \end{align*}
    The Ansatz \eqref{eqn:tipi-2-sig} implies
    \begin{align*}
        \pi_1 \delta \ti{\pi}_2 \sigma_n &= \sum_{j=1}^{\floor{\frac{n-1}{2}} - 1} (- q)^j \conj{q}^j \del_x^{n-1-2j} q \sum_{t=1}^{n-2-2j} (-1)^t G_{n,j,t} 
        \\ &+ \sum_{j=0}^{\floor{\frac{n-1}{2}} - 2} (- q)^{j+2} \conj{q}^j \del_x^{n-3-2j} \conj{q} \sum_{t=1}^{\floor{\frac{n-1}{2}}-1-j} (-1)^t (1 - (-1)^n) F_{n,j,t} \,.
    \end{align*}
    Accordingly, we define
    \begin{align*}
        \ti{G}_{n,j} &= \sum_{t=1}^{n-2-2j} (-1)^t G_{n,j,t} 
        & \ti{F}_{n,j} &= (1 - (-1)^n) \sum_{t=1}^{\floor{\frac{n-1}{2}}-1-j} (-1)^t F_{n,j,t} \,.
    \end{align*}
    Then
    \begin{align*}
        \ti{G}_{n+1,j} &= 2 \sum_{k=0}^{j-2} \sum_{t=1}^{n-1-2j} (-1)^t G_{n-2k-1,j-k-1,t} C_k + 2 \sum_{k=0}^{j-2} \sum_{t=1}^{n-1-2j} (-1)^t D_{t+1+2k,k} E_{n-1-2k-t,j-2-k} \\
        &+ \mathds{1}_{\{j \neq \frac{n}{2} - 1\}} \sum_{t=1}^{n-1-2j} (-1)^t (\mathds{1}_{\{t \neq n - 1 - 2j\}} G_{n,j,t} + \mathds{1}_{\{t \neq 1\}} G_{n,j,t-1}) \\
        &+ \sum_{t=1}^{n-1-2j} (-1)^t (\mathds{1}_{\{t = 1\}} (j+1) E_{n,j-1} - \mathds{1}_{\{t = n - 1 - 2j\}} j D_{n,j}) \\
        &= 2 \sum_{k=0}^{j-1} \ti{G}_{n-2k-1,j-k-1} C_k + 2 \sum_{k=0}^{j-2} \sum_{t=1}^{n-1-2j} (-1)^t D_{t+1+2k,k} E_{n-1-2k-t,j-2-k} \\
        &- (j+1) E_{n,j-1} + (-1)^n j D_{n,j}
    \end{align*}
    and
    \begin{align*}
        \ti{F}_{n+1,j} &= 2 \mathds{1}_{\{n \text{ even}\}} \Bigg( 2 \sum_{k=0}^{j-1} \sum_{t=1}^{\floor{\frac{n}{2}}-1-j} (-1)^t F_{n-2k-1,j-k-1,t} C_k 
        \\ &+ \sum_{k=0}^j \sum_{t=1}^{\floor{\frac{n}{2}}-1-j} (-1)^t (1 + \mathds{1}_{\{t \neq \frac{n}{2}-1-j\}}) D_{t+1+2k,k} D_{n-1-2k-t,j-k}
        + \sum_{t=1}^{\floor{\frac{n}{2}}-1-j} (-1)^t \Big( \mathds{1}_{\{t = 1\}} (j + 1) D_{n,j+1} 
        \\ &+ \mathds{1}_{\{j \neq \frac{n}{2} - 2\}} \Big( \mathds{1}_{\{t \neq \frac{n}{2} - 1 - j\}} F_{n,j,t} + \mathds{1}_{\{t \neq 1\}} F_{n,j,t-1} + \mathds{1}_{\{t = \frac{n-1}{2} - 1 - j\}} F_{n,j,t} \Big) \Big) \Bigg) \\
        &= 2 \mathds{1}_{\{n \text{ even}\}} \Bigg( \sum_{k=0}^{j-1} \ti{F}_{n-2k-1,j-k-1} C_k + \sum_{k=0}^j \sum_{t=1}^{\floor{\frac{n}{2}}-1-j} (-1)^t (1 + \mathds{1}_{\{t \neq \frac{n}{2}-1-j\}}) D_{t+1+2k,k} D_{n-1-2k-t,j-k} \\
        &- (j + 1) D_{n,j+1} + \mathds{1}_{\{n \text{ odd}\}} 2 (-1)^{\frac{n-1}{2} - 1 - j} F_{n,j,\frac{n-1}{2} - 1 - j} \Bigg) \\
        &= \mathds{1}_{\{n \text{ even}\}} \Bigg( 2 \sum_{k=0}^{j-1} \ti{F}_{n-2k-1,j-k-1} C_k + 2 \sum_{k=0}^j \sum_{t=1}^{\floor{\frac{n}{2}}-1-j} (-1)^t (1 + \mathds{1}_{\{t \neq \frac{n}{2}-1-j\}}) D_{t+1+2k,k} D_{n-1-2k-t,j-k} \\
        &- 2 (j + 1) D_{n,j+1} \Bigg) \,.
    \end{align*}
    Recall now that
    \begin{align*}
        D_{n,j} &= 4^j \binom{\frac{n}{2} - 1}{j}
        & E_{n,j} &= \sum_{l=0}^j (-1)^l (C_{l+1} - 2 C_l) 4^{j-l} \begin{pmatrix}
            \frac{n - 1}{2} - 1 \\ j - l
        \end{pmatrix} \,.
    \end{align*}
    A calculation with binomial identities (or alternatively verified manually using generating functions) yields
    \begin{align*}
        \sum_{k=0}^{j-2} \sum_{t=1}^{n-1-2j} (-1)^t D_{t+1+2k,k} E_{n-1-2k-t,j-2-k}
        &= - \mathds{1}_{\{n \text{ even}\}} E_{n,j-2} \,.
    \end{align*}
    Therefore
    \begin{align*}
        \ti{G}_{n+1, j} &= 2 \sum_{k=0}^{j-2} \ti{G}_{n-2k-1,j-k-1} C_k - 2 \mathds{1}_{\{n \text{ even}\}} \sum_{l=0}^{j-2} (C_{l+1} - 2 C_l) (-1)^l j 4^{j-2-l} \binom{\frac{n-1}{2}-1}{j-2-l} \\
        &- (j+1) \sum_{l=0}^{j-1} (C_{l+1} - 2 C_l) (-1)^l 4^{j-1-l} \binom{\frac{n-1}{2}-1}{j-1-l} + (-1)^n 4^j \binom{\frac{n}{2}-1}{j} \\
        &= 2 \sum_{k=0}^{j-2} \ti{G}_{n-2k-1,j-k-1} C_k + (-1)^n j 4^j \binom{\frac{n}{2}-1}{j} \\
        &- \sum_{l=0}^{j-1} (C_{l+1} - 2 C_l) (-1)^l \left( (j+1) 4^{j-1-l} \binom{\frac{n-1}{2}-1}{j-1-l} + 2 \mathds{1}_{\{n \text{ even}\}} 4^{j-2-l} \binom{\frac{n-1}{2}-1}{j-2-l} \right) \,.
    \end{align*}
    We claim that
    \begin{align*}
        \ti{G}_{n,j} &= (-1)^{n+1} 4\left(\floor*{\frac{n-1}{2}}-j\right) \sum_{k=0}^{j-1} C_k 4^{j-1-k} \binom{\floor*{\frac{n-1}{2}}-\frac12-k}{j-1-k} + 2 \mathds{1}_{\{n \text{ odd}\}} 4^{j-1} \binom{\frac{n}{2}-1}{j-1} \,.
    \end{align*}
    Again, a set of tedious generating function calculations needs to be performed. We calculate
    \begin{align*}
        \sum_{n=0,j=0}^\infty X^n Y^j \ti{G}_{n+1,j}
        &= \frac{- 4 Y}{\sqrt{1 + 4 Y}} \frac{- X}{1 - X^2 (1 + 4 Y)} + \frac{- 2 Y}{\sqrt{1 + 4 Y}} \frac{1}{1 - X^2 (1 + 4 Y)}
        \\ &+ \left(1 - \frac{1}{\sqrt{1 + 4 Y}} \right) 2 \sqrt{1 + 4 Y} \frac{X^2 (1 - X)}{(1 - X^2 (1 + 4 Y))^2} \,.
    \end{align*}
    Next, we calculate
    \begin{align*}
        \sum_{n=0,j=0}^\infty X^n Y^j 2 \sum_{k=0}^{j-2} C_k \ti{G}_{n-2k-1,j-k-1}
        &= \left( 1 - \frac{1}{\sqrt{1 + 4 Y}} \right)^2 2 \sqrt{1 + 4 Y} \frac{X^2 (1 - X)}{(1 - X^2 (1 + 4 Y))^2}
        \\ &+ \left( 1 - \frac{1}{\sqrt{1 + 4 Y}} \right) \frac{2 Y}{\sqrt{1 + 4 Y}} \left( - 2 - \frac{2}{\sqrt{1 + 4 Y}} \right) \frac{- X}{1 - X^2 (1 + 4 Y)}
        \\ &+ \left( 1 - \frac{1}{\sqrt{1 + 4 Y}} \right) \frac{2 Y}{\sqrt{1 + 4 Y}} \left( - 1 - \frac{2}{\sqrt{1 + 4 Y}} \right) \frac{1}{1 - X^2 (1 + 4 Y)} \,.
    \end{align*}
    Lastly, we calculate
    \begin{align*}
        & \sum_{n=0,j=0}^\infty X^n Y^j (-1)^n j 4^j \binom{\frac{n}{2} - 1}{j} \\
        &= \left( \frac{1}{\sqrt{1 + 4 Y}} - X \right) \frac{4 Y}{\sqrt{1 + 4 Y}} \frac{X^2}{(1 - X^2 (1 + 4 Y))^2}
        + \left( \frac12 X - \frac{1}{\sqrt{1 + 4 Y}} \right) \frac{4 Y}{\sqrt{1 + 4 Y}^3} \frac{1}{1 - X^2 (1 + 4 Y)}
    \end{align*}
    and
    \begin{align*}
        &- \sum_{n=0, j=0}^\infty X^n Y^j \sum_{l=0}^{j-1} (C_{l+1} - 2 C_l) (-1)^l \left( (j+1) 4^{j-1-l} \binom{\frac{n-1}{2}-1}{j-1-l} + 2 \mathds{1}_{\{n \text{ even}\}} 4^{j-2-l} \binom{\frac{n-1}{2}-1}{j-2-l} \right) \\
        &= \left( \sqrt{1 + 4 Y} + \frac{1}{\sqrt{1 + 4 Y}} - 2 \right) \left( X + \frac{1}{\sqrt{1 + 4 Y}} \right) \frac{X^2}{(1 - X^2 (1 + 4 Y))^2} \\
        &+ \frac{- 2 Y}{\sqrt{1 + 4 Y}^3} \frac{- X}{1 - X^2 (1 + 4 Y)} + \left( \frac12 - \frac{1}{\sqrt{1 + 4 Y}} + \frac12 \frac{1}{1 + 4 Y} + \frac{1}{\sqrt{1 + 4 Y}^3} - \frac{1}{(1 + 4 Y)^2} \right) \frac{1}{1 - X^2 (1 + 4 Y)} \,.
    \end{align*}
    Summing up all the contributions verifies that the claimed formula for $\ti{G}_{n,j}$ satisfies the given recurrence relation. 
    We move on to $\ti{F}_{n,j}$.
    For even $n = 2 m$ we calculate
    \begin{align*}
        \sum_{k=0}^j \sum_{t=1}^{\floor{\frac{n}{2}}-1-j} (-1)^t (1 + \mathds{1}_{\{t \neq \frac{n}{2}-1-j\}}) D_{t+1+2k,k} D_{n-1-2k-t,j-k}
        &= - \mathds{1}_{\{j \leq m - 2\}} (- 4)^j \binom{j-m}{j} \,.
    \end{align*}
    Therefore,
    \begin{align*}
        \ti{F}_{2m+1,j} &= 2 \sum_{k=0}^{j-1} \ti{F}_{2m-2k-1,j-k-1} C_k - 2 (- 4)^j \binom{j-m}{j} - 2 (j + 1) 4^{j+1} \binom{m-1}{j+1} \,.
    \end{align*}
    We next claim that
    \begin{align*}
        \ti{F}_{n,j} &= - \mathds{1}_{\{n \text{ odd}\}} \left(8 \floor{\frac{n-1}{2}} - 8 j - 6\right) 4^j \binom{\frac{n}{2} - 1}{j}
        &= \begin{cases}
            - (8 m - 8 j - 6) 4^j \binom{m - \frac12}{j} &, n = 2m + 1\\
            0 &, n = 2m
        \end{cases} \,.
    \end{align*}
    Plugging everything into the recurrence relation yields the trivial identity for odd $n$ and for even $n = 2 m$ the binomial identity
    \begin{align*}
        (8 m - 8 j - 6) \binom{m-\frac12}{j} &= 2 \sum_{k=0}^{j-1} (8 m - 8 j - 6) 4^{-k-1} \binom{m-k-\frac32}{j-k-1} C_k \,.
        + 2 (- 1)^j \binom{j-m}{j} + 8 (j + 1) \binom{m-1}{j+1}
    \end{align*}
    It remains to prove this for all $(m, j) \in \Z^2$.
    We verify this by once more calculating and comparing the generating functions
    \begin{align*}
        & \sum_{m=0}^\infty \sum_{j=0}^\infty X^m Y^j 2 \sum_{k=0}^{j-1} (8 (m-j) - 6) 4^{-k-1} \binom{m-k-1-\frac12}{j-k-1} C_k \\
        &= \left( - 8 \frac{Y}{\sqrt{1 + Y}^3} + 8 \frac{1}{\sqrt{1 + Y}} \right) \frac{X (1 + Y)}{(1 - X (1 + Y))^2} + \left( 4 \frac{Y}{\sqrt{1 + Y}^3} - 6 \frac{1}{\sqrt{1 + Y}} \right) \frac{1}{1 - X (1 + Y)} \\
        &- \frac{8}{(1 + Y)^2} \frac{X (1 + Y)}{(1 - X (1 + Y))^2} + \frac{6 - 2 Y}{(1 + Y)^2} \frac{1}{1 - X (1 + Y)}
    \end{align*}
    and
    \begin{align*}
        & \sum_{m=0}^\infty \sum_{j=0}^\infty X^m Y^j 2 (-1)^j \binom{- m + j}{j} + \sum_{m=0}^\infty \sum_{j=0}^\infty X^m Y^j 8 (j+1) \binom{m - 1}{j + 1} \\
        &= \frac{2 Y - 6}{(1 + Y)^2} \frac{1}{1 - X (1 + Y)} + \frac{8}{(1 + Y)^2} \frac{X (1 + Y)}{(1 - X (1 + Y))^2}
    \end{align*}
    with 
    \begin{align*}
        & \sum_{m=0}^\infty \sum_{j=0}^\infty X^m Y^j (8 (m-j) - 6) \binom{m-\frac12}{j} \\
        &= \left( - 8 \frac{Y}{\sqrt{1 + Y}^3} + 8 \frac{1}{\sqrt{1 + Y}} \right) \frac{X (1 + Y)}{(1 - X (1 + Y))^2} + \left( 4 \frac{Y}{\sqrt{1 + Y}^3} - 6 \frac{1}{\sqrt{1 + Y}} \right) \frac{1}{1 - X (1 + Y)} \,.
    \end{align*}
\end{proof}

\subsection{Formulas for \texorpdfstring{$\pi_0 \delta \sigma_n$}{π₀ δ σₙ} and \texorpdfstring{$\pi_1 \delta \sigma_n$}{π₁ δ σₙ}}

\begin{lemma} \label{lem:20}
    We have
    \begin{align*}
        \pi_0 \delta \sigma_n &= \mathds{1}_{\{n \text{ odd}\}} \frac{n+1}{2} C_{\frac{n-1}{2}} \conj{q}^{\frac{n-1}{2}} (- q)^{\frac{n+1}{2}} \\
        \pi_1 \delta \sigma_n &= \sum_{j=0}^{\floor{\frac{n}{2}} - 2} J_{n,j} (- q)^{j+2} \conj{q}^j \del_x^{n - 3 - 2 j} \conj{q} 
        + \sum_{j=0}^{\floor{\frac{n}{2}} - 1} K_{n,j} (- q)^j \conj{q}^j \del_x^{n - 1 - 2 j} q \,,
    \end{align*}
    where $J_{n,j}$ and $K_{n,j}$ are defined in \eqref{eqn:coeffs}.
\end{lemma}
\begin{proof}
    Recall that
    \begin{align*}
        \pi_0 \sigma_n &= \mathds{1}_{\{n \text{ odd}\}} C_{\frac{n-1}{2}} \conj{q}^{\frac{n+1}{2}} (- q)^{\frac{n+1}{2}} \\
        \pi_1 \sigma_n &= \sum_{j=0}^{\floor{\frac{n}{2}} - 1} D_{n,j} (- q)^{j+1} \conj{q}^j \del_x^{n - 1 - 2 j} \conj{q} + \sum_{j=0}^{\floor{\frac{n}{2}} - 2} E_{n,j} (- q)^{j+1} \conj{q}^{j+2} \del_x^{n - 3 - 2 j} q \\
        \pi_1 \delta \ti{\pi}_2 \sigma_n &= \sum_{j=0}^{\floor{\frac{n-1}{2}} - 2} (- q)^{j+2} \conj{q}^j \del_x^{n-3-2j} \conj{q} \ti{F}_{n,j} + \sum_{j=1}^{\floor{\frac{n-1}{2}} - 1} (- q)^j \conj{q}^j \del_x^{n-1-2j} q \ti{G}_{n,j} \,.
    \end{align*}
    We easily obtain
    \begin{align*}
        \pi_0 \delta \sigma_n &= \pi_0 \delta \pi_0 \sigma_n 
        = \mathds{1}_{\{n \text{ odd}\}} \frac{n+1}{2} C_{\frac{n-1}{2}} \conj{q}^{\frac{n-1}{2}} (- q)^{\frac{n+1}{2}} \,.
    \end{align*}
    Furthermore,
    \begin{align*}
        \pi_1 \delta \pi_1 \sigma_n &= \sum_{j=1}^{\floor{\frac{n}{2}} - 1} j D_{n,j} (- q)^{j+1} \conj{q}^{j-1} \del_x^{n - 1 - 2 j} \conj{q} + \sum_{j=0}^{\floor{\frac{n}{2}} - 2} (j+2) E_{n,j} (- q)^{j+1} \conj{q}^{j+1} \del_x^{n - 3 - 2 j} q \\
        &+ \pi_1 \sum_{j=0}^{\floor{\frac{n}{2}} - 1} (- 1)^{n - 1 - 2 j} D_{n,j} \sum_{k=0}^{n-1-2j} \binom{n-1-2j}{k} \del_x^k((- q)^{j+1}) \del_x^{n-1-2j-k}(\conj{q}^j) \\
        &= \sum_{j=0}^{\floor{\frac{n}{2}} - 2} (j+1) D_{n,j+1} (- q)^{j+2} \conj{q}^j \del_x^{n - 3 - 2 j} \conj{q} + \sum_{j=1}^{\floor{\frac{n}{2}} - 1} (j+1) E_{n,j-1} (- q)^j \conj{q}^j \del_x^{n - 1 - 2 j} q \\
        &+ \sum_{j=1}^{\floor{\frac{n}{2}} - 1} j (- 1)^{n - 1} D_{n,j} (- q)^{j+1} \conj{q}^{j-1} \del_x^{n-1-2j} \conj{q} + \sum_{j=0}^{\floor{\frac{n}{2}} - 1} (j+1) (- 1)^n D_{n,j} (- q)^j \conj{q}^j \del_x^{n-1-2j} q \\
        &= (1 - (-1)^n) \sum_{j=0}^{\floor{\frac{n}{2}} - 2} (j+1) D_{n,j+1} (- q)^{j+2} \conj{q}^j \del_x^{n - 3 - 2 j} \conj{q} 
        \\ &+ \sum_{j=0}^{\floor{\frac{n}{2}} - 1} (j+1) (\mathds{1}_{\{j \geq 1\}} E_{n,j-1} + (-1)^n D_{n,j}) (- q)^j \conj{q}^j \del_x^{n - 1 - 2 j} q \,.
    \end{align*}
    In total,
    \begin{align*}
        \pi_1 \delta \sigma_n &= \sum_{j=0}^{\floor{\frac{n}{2}} - 2} \big( \mathds{1}_{\{0 \leq j \leq \floor{\frac{n-1}{2}} - 2\}} \ti{F}_{n,j} + \mathds{1}_{\{n \text{ odd}\}} 2 (j+1) D_{n,j+1} \big) (- q)^{j+2} \conj{q}^j \del_x^{n - 3 - 2 j} \conj{q} \\
        &+ \sum_{j=0}^{\floor{\frac{n}{2}} - 1} \big( \mathds{1}_{\{1 \leq j \leq \floor{\frac{n-1}{2}} - 1\}} \ti{G}_{n,j} + (j+1) (\mathds{1}_{\{j \geq 1\}} E_{n,j-1} + (-1)^n D_{n,j}) \big) (- q)^j \conj{q}^j \del_x^{n - 1 - 2 j} q \\
        &= \sum_{j=0}^{\floor{\frac{n}{2}} - 2} J_{n,j} (- q)^{j+2} \conj{q}^j \del_x^{n - 3 - 2 j} \conj{q} + \sum_{j=0}^{\floor{\frac{n}{2}} - 1} K_{n,j} (- q)^j \conj{q}^j \del_x^{n - 1 - 2 j} q
    \end{align*}
    for certain coefficients $J_{n,j}$ and $K_{n,j}$, whose claimed formulas we need to verify by simplifying the expressions above. 
    Recall the definitions
    \begin{align*}
        D_{n,j} &= 4^j \begin{pmatrix}
            \frac{n}{2} - 1 \\ j
        \end{pmatrix} \\
        E_{n,j} &= \sum_{l=0}^j (-1)^l (C_{l+1} - 2 C_l) 4^{j-l} \begin{pmatrix}
            \frac{n - 1}{2} - 1 \\ j - l
        \end{pmatrix} \\
        \ti{F}_{n,j} &= - \mathds{1}_{\{n \text{ odd}\}} \left(8 \floor{\frac{n-1}{2}} - 8 j - 6\right) 4^j \binom{\frac{n}{2} - 1}{j} \\
        \ti{G}_{n,j} &= (-1)^{n+1} 4\left(\floor*{\frac{n-1}{2}}-j\right) \sum_{k=0}^{j-1} C_k 4^{j-1-k} \binom{\floor*{\frac{n-1}{2}}-\frac12-k}{j-1-k} + 2 \mathds{1}_{\{n \text{ odd}\}} 4^{j-1} \binom{\frac{n}{2}-1}{j-1} \,.
        \end{align*}
    Let $0 \leq j \leq \floor{\frac{n}{2}} - 2$ and $n = 2 m$ or $n = 2 m + 1$. We have
    \begin{align*}
        J_{2m,j} &= \mathds{1}_{\{j \leq \floor{\frac{2m-1}{2}} - 2\}} \ti{F}_{2m,j} + \mathds{1}_{\{2m \text{ odd}\}} 2 (j+1) D_{2m,j+1} = 0
    \end{align*}
    and
    \begin{align*}
        J_{2m+1,j} &= \mathds{1}_{\{j \leq \floor{\frac{2m+1-1}{2}} - 2\}} \ti{F}_{2m+1,j} + \mathds{1}_{\{2m+1 \text{ odd}\}} 2 (j+1) D_{2m+1,j+1} \\
        &= - \mathds{1}_{\{j \leq m - 2\}} (8 m - 8 j - 6) 4^j \binom{m - \frac12}{j} + 2 (j+1) 4^{j+1} \binom{m - \frac12}{j+1} \\
        &= (8 j + 8) 4^j \binom{m + \frac12}{j + 1} - (8 m + 2) 4^j \binom{m - \frac12}{j} \,.
    \end{align*}
    A computation shows that
    \begin{align*}
        \sum_{m=0}^\infty \sum_{j=0}^\infty X^m Y^j (8 j + 8) 4^j \binom{m + \frac12}{j + 1} - \sum_{m=0}^\infty \sum_{j=0}^\infty X^m Y^j (8 m + 2) 4^j \binom{m - \frac12}{j}
        &= \frac{2}{\sqrt{1 + 4 Y}} \frac{1}{1 - X (1 + 4 Y)}
        \\ &= \sum_{m=0}^\infty \sum_{j=0}^\infty X^m Y^j 4^j 2 \binom{m-\frac12}{j} \,.
    \end{align*}
    Hence
    \begin{align*}
        J_{2m+1,j} &= 4^j 2 \binom{m-\frac12}{j} \qquad \text{ and } \qquad J_{n,j} = 2 \mathds{1}_{\{n \text{ odd}\}} 4^j \binom{\frac{n}{2} - 1}{j} \,.
    \end{align*}
    Using the notation $[u^j]$ for the extraction of the coefficient in front of $u^j$ from a formal power series in the symbol $u$, we can also write this as
    \begin{align*}
        J_{n,j} &= \mathds{1}_{\{n \text{ odd}\}} [u^j] 2 (1 + 4 u)^{\frac{n-2}{2}} \,.
    \end{align*}
    Now let $0 \leq j \leq \floor{\frac{n}{2}} - 1$. We have
    \begin{align*}
        K_{n,j} &= \mathds{1}_{\{1 \leq j \leq \floor{\frac{n-1}{2}} - 1\}} \ti{G}_{n,j} + (j+1) (\mathds{1}_{\{j \geq 1\}} E_{n,j-1} + (-1)^n D_{n,j}) \\
        &= \ti{G}_{n,j} + (j+1) (E_{n,j-1} + (-1)^n D_{n,j}) \\
        &= \left( (-1)^{n+1} 4\left(\floor*{\frac{n-1}{2}}-j\right) \sum_{k=0}^{j-1} C_k 4^{j-1-k} \binom{\floor*{\frac{n-1}{2}}-\frac12-k}{j-1-k} + 2 \mathds{1}_{\{n \text{ odd}\}} 4^{j-1} \binom{\frac{n}{2}-1}{j-1} \right) \\
        &+ (j+1) \left( \mathds{1}_{\{j \geq 1\}} \sum_{l=0}^{j-1} (-1)^l (C_{l+1} - 2 C_l) 4^{j-1-l} \begin{pmatrix}
            \frac{n - 1}{2} - 1 \\ j - 1 - l
        \end{pmatrix} + (-1)^n 4^j \begin{pmatrix}
            \frac{n}{2} - 1 \\ j
        \end{pmatrix} \right) \,.
    \end{align*} 
    We claim that
    \begin{align}
        \nonumber K_{n,j} &= (-1)^n D_{n+1,j} + 2 \mathds{1}_{\{n \text{ odd}\}} E_{n,j-2} \\
        \\ \label{eqn:K-1} &= (-1)^n 4^j \binom{\frac{n+1}{2} - 1}{j} + 2 \mathds{1}_{\{n \text{ odd}\}} \sum_{l=0}^{j-2} (-1)^l (C_{l+1} - 2 C_l) 4^{j-2-l} \binom{\frac{n-1}{2}-1}{j-2-l} \,.
        \\ \label{eqn:K-2} &= \begin{cases}
            [u^j] (1 + 4 u)^{\frac{n-1}{2}} &, n \text{ even} \\
            [u^j] (- 1 - 2 u) (1 + 4 u)^{\frac{n-2}{2}} &, n \text{ odd}
        \end{cases} \,.
    \end{align}
    Note that
    \begin{align*}
        & \sum_{n=0,j=0}^\infty X^n Y^j (-1)^n 4^j \binom{\frac{n+1}{2} - 1}{j} + \sum_{n=0,j=0}^\infty X^n Y^j 2 \mathds{1}_{\{n \text{ odd}\}} \sum_{l=0}^{j-2} (-1)^l (C_{l+1} - 2 C_l) 4^{j-2-l} \binom{\frac{n-1}{2}-1}{j-2-l} \\
        &= \frac{1}{\sqrt{1 + 4 Y}} \frac{1}{1 - X^2 (1 + 4 Y)} + \left( - \frac12 \frac{1}{\sqrt{1 + 4 Y}} - \frac12 \sqrt{1 + 4 Y} \right) \frac{X}{1 - X^2 (1 + 4 Y)} \,.
        \\ &= \frac{1}{\sqrt{1 + 4 Y}} \frac{1 - X (1 + 2 Y)}{1 - X^2 (1 + 4 Y)}
        \\ &= \sum_{m=0,j=0}^\infty X^{2m} Y^j [u^j] (1 + 4 u)^{\frac{2m-1}{2}}
        + \sum_{m=0,j=0}^\infty X^{2m+1} Y^j [u^j] (- 1 - 2 u) (1 + 4 u)^{\frac{2m-1}{2}} \,.
    \end{align*}
    This shows that \eqref{eqn:K-1} and \eqref{eqn:K-2} agree.
    To show that these are indeed formulas for $K_{n,j}$, we have to compute the following generating function, which mostly reduces to previous calculations:
    \begin{align*}
        & \sum_{n=0,j=0}^\infty X^n Y^j K_{n,j} \\
        &= X \sum_{n=0,j=0}^\infty X^n Y^j \ti{G}_{n+1,j} 
        + \sum_{j=0}^\infty Y^j \ti{G}_{0,j}
        + \sum_{n=0,j=0}^\infty X^n Y^j (j+1) E_{n,j-1}
        + \sum_{n=0,j=0}^\infty X^n Y^j (j+1) (-1)^n D_{n,j} \,.
    \end{align*}
    First, note that $\ti{G}_{0,j} = 0$. Second, recall that
    \begin{align*}
        X \sum_{n=0,j=0}^\infty X^n Y^j \ti{G}_{n+1,j}
        &= X\left( \frac{- 4 Y}{\sqrt{1 + 4 Y}} \frac{- X}{1 - X^2 (1 + 4 Y)} + \frac{- 2 Y}{\sqrt{1 + 4 Y}} \frac{1}{1 - X^2 (1 + 4 Y)} \right.
        \\ &+ \left. \left(1 - \frac{1}{\sqrt{1 + 4 Y}} \right) 2 \sqrt{1 + 4 Y} \frac{X^2 (1 - X)}{(1 - X^2 (1 + 4 Y))^2} \right) \,.
    \end{align*}
    For the new term involving $\ti{G}_{0,j}$, we have
    \begin{align*}
        \sum_{j=0}^\infty Y^j \ti{G}_{0,j}
        &= \frac12 \left( 1 - \frac{1}{\sqrt{1 + 4 Y}} \right) \left( - \frac12 \frac{16 Y}{\sqrt{1 + 4 Y}^3} + \frac{4}{\sqrt{1 + 4 Y}} \right) 
        + \frac{Y}{\sqrt{1 + 4 Y}} \left( \frac{- 16 Y}{\sqrt{1 + 4 Y}^3} + \frac{4}{\sqrt{1 + 4 Y}} \right) \\
        &= \frac{1}{\sqrt{1 + 4 Y}} + \frac{1}{\sqrt{1 + 4 Y}^3} + \frac{2}{(1 + 4 Y)^2} \,.
    \end{align*}
    Lastly, we calculate
    \begin{align*}
        \sum_{n=0,j=0}^\infty X^n Y^j (j+1) (-1)^n D_{n,j}
        &= \left( \frac{1}{\sqrt{1 + 4 Y}} - \sqrt{1 + 4 Y} \right) \frac{X^3}{(1 - X^2 (1 + 4 Y))^2} + \frac{4 Y}{1 + 4 Y}  \frac{X^2}{(1 - X^2 (1 + 4 Y))^2}
        \\ &- \frac12 \frac{1}{\sqrt{1 + 4 Y}} \left( 1 + \frac{1}{1 + 4 Y} \right) \frac{X}{1 - X^2 (1 + 4 Y)} + \frac{1}{(1 + 4 Y)^2} \frac{1}{1 - X^2 (1 + 4 Y)}
    \end{align*}
    and
    \begin{align*}
        & \sum_{n=0, j=0}^\infty X^n Y^j (j+1) E_{n,j-1}
        \\ &= \left( - \sqrt{1 + 4 Y} - \frac{1}{\sqrt{1 + 4 Y}} + 2 \right) \frac{X^3}{(1 - X^2 (1 + 4 Y))^2}
        + \left( - 1 - \frac{1}{1 + 4 Y} + \frac{2}{\sqrt{1 + 4 Y}} \right) \frac{X^2}{(1 - X^2 (1 + 4 Y))^2}
        \\ &- \frac{2 Y}{\sqrt{1 + 4 Y}^3} \frac{X}{1 - X^2 (1 + 4 Y)} + \left( \frac{1}{(1 + 4 Y)^2} - \frac{1}{\sqrt{1 + 4 Y}^3} \right) \frac{1}{1 - X^2 (1 + 4 Y)} \,.
    \end{align*}
    Adding up all the contributions, we find that the generating functions of both ofour formulas for $K_{n,j}$ agree.
\end{proof}

\subsection{Proofs of Theorem \ref{thm:2} and Proposition \ref{prop:2}}

\begin{proof}[Proof of (i) from Theorem \ref{thm:2}]
    We can write \eqref{eqn:NLS-n} as
    \begin{align*}
        \del_{t_n} q &= \delta \mc{H}^\NLS_n = - (- i)^n \delta \sigma_{n+1} \,,
    \end{align*}
    so Lemma \ref{lem:20} directly implies (i).
\end{proof}

\begin{proof}[Proof of (ii) from Theorem \ref{thm:2}]
    In order to proceed, we must use the Lagrange inversion formula. We refer the reader to \cite{LagrangeInversion} for a detailed exposition and adopt the notations from this reference.
    If $u(t)$ is a formal power series in $t$, then $[t^n] u(t)$ denotes the application of the linear functional of extracting the coefficient in front of $t^n$ (e.g. $[t^2] (1 + u)^3 = 3$).
    With this notation, a version of the Lagrange inversion formula reads as follows.
    Let $F(t)$, $u(t)$, and $\phi(t)$ be formal power series. Assume that $[t^0] \phi(t) = 0$, which implies that
    there exists a formal power series solution $w(t)$ to the implicit equation
    \begin{align*}
        w(t) &= t \phi(w(t)) \,.
    \end{align*}
    Then
    \begin{align*}
        [t^n] F(w(t)) &= [t^n] F(t) \phi(t)^{n-1} (\phi(t) - t \phi'(t)) \,.
    \end{align*}
    In the subsequent calculations, the relation $\approx$ shall refer to equality up to terms in $\mc{O}^{2,n-2}_q(q_x)$.
    Square brackets after a differential operator denote operator application.
    We use the Lagrange inversion formula to calculate
    \begin{align*}
        \frac{\delta \mc{H}_\NLS^{2m+1}}{\delta \conj{q}} &\approx \sum_{j=0}^m [t^j] \rho^j (1 + 4 t)^{\frac{2m+1}{2}} (- \del_x^2)^{m-j} [i \del_x q]
        \approx [t^m] \frac{(1 + 4 \rho t)^{\frac32}}{1 + t \del_x^2} (1 + 4 \rho t)^{m-1} [i \del_x q]
        \\ &\approx [t^m] \frac{(1 + 4 \rho \frac{t}{1 - 4 \rho t})^{\frac32}}{1 + \frac{t}{1 - 4 \rho t} \del_x^2} [i \del_x q]
        \approx [t^m] \frac{1}{(1 - 4 \rho t)^{\frac12} (1 - 4 \rho t + t \del_x^2)} [i \del_x q]
    \end{align*}
    and
    \begin{align*}
        \frac{\delta \mc{H}_\NLS^{2m}}{\delta \conj{q}} 
        &\approx [t^m] (1 - 4 t)^{-\frac12} \rho^m q
        + \sum_{j=0}^{m-2} [t^j] \rho^j 2 q^2 (1 + 4 t)^{\frac{2m-1}{2}} (- \del_x^2)^{m-2-j} [- \del_x^2 \conj{q}]
        \\ &+ \sum_{j=0}^{m-1} [t^j] \rho^j (1 + 2 t) (1 + 4 t)^{\frac{2m-1}{2}} (- \del_x^2)^{m-1-j} [- \del_x^2 \conj{q}]
        \\ &\approx [t^m] (1 - 4 \rho t)^{-\frac12} q
        + [t^m] \frac{2 q^2 t^2 (1 + 4 \rho t)^\frac12}{1 + t \del_x^2} (1 + 4 \rho t)^{m-1} [- \del_x^2 \conj{q}]
        \\ &+ [t^m] \frac{t (1 + 2 \rho t) (1 + 4 \rho t)^{\frac12}}{{1 + t \del_x^2}} (1 + 4 \rho t)^{m-1} [- \del_x^2 q]
        \\ &\approx [t^m] (1 - 4 \rho t)^{-\frac12} q
        + [t^m] \frac{2 q^2 \frac{t^2}{(1 - 4 \rho t)^2} (1 + 4 \rho \frac{t}{1 - 4 \rho t})^\frac12}{1 + \frac{t}{1 - 4 \rho t} \del_x^2} [- \del_x^2 \conj{q}]
        \\ &+ [t^m] \frac{\frac{t}{1 - 4 \rho t} (1 + 2 \rho \frac{t}{1 - 4 \rho t}) (1 + 4 \rho \frac{t}{1 - 4 \rho t})^{\frac12}}{{1 + \frac{t}{1 - 4 \rho t} \del_x^2}} [- \del_x^2 q]
        \\ &\approx [t^m] (1 - 4 \rho t)^{-\frac12} q
        + [t^m] \frac{1}{(1 - 4 \rho t)^{\frac32} (1 - 4 \rho t + t \del_x^2)} [- 2 t^2 \del_x^2 (\rho q) + (6 t^2 \rho - t) \del_x^2 q] \,.
    \end{align*}
    Here we have used
    \begin{align*}
        2 t^2 q^2 \del_x^2 \conj{q} + t (1 - 2 \rho t) \del_x^2 q 
        &\approx t \del_x^2 q + 2 t^2 (q^2 \del_x^2 \conj{q} - \rho \del_x^2 q)
        \approx t \del_x^2 q + 2 t^2 (\del_x^2 (\rho q) - 3 \rho \del_x^2 q) \,.
    \end{align*}
    Recall that for $\lambda \in \mc{Q}_1$ with $z = \sqrt{\lambda^2 - 1}$ and $\lambda = \sqrt{z^2 + 1}$, we have from \eqref{eqn:trans-expand-NLS} the asymptotic expansion
    \begin{align*}
        - 2 i \lambda \frac{\delta \log a(\lambda)}{\delta \conj{q}} &= \sum_{n=0}^\infty (2 \lambda)^{-n} \frac{\delta \mc{H}_\NLS^n}{\delta \conj{q}}
        \sim \sum_{m=0}^\infty (4 \lambda^2)^{-m} \left( \frac{\delta \mc{H}_\NLS^{2m}}{\delta \conj{q}} + (2 \lambda)^{-1} \frac{\delta \mc{H}_\NLS^{2m+1}}{\delta \conj{q}} \right) \,.
    \end{align*}
    Substituting the above formulas replaces each $t$ with $(4 \lambda^2)^{-1}$. After further simplification (up to $\approx$), we obtain the claimed formula.
\end{proof}

\begin{proof}[Proof of (iii) from Theorem \ref{thm:2}]
    Combining \eqref{eqn:struct-NLS-1}--\eqref{eqn:struct-NLS-2} with \eqref{eqn:NLS-GP-3}--\eqref{eqn:NLS-GP-4}, we can write the equations \eqref{eqn:GP-n} of the \GP hierarchy as
    \begin{align}
        \label{eqn:struct-GP-1} i \del_{t_{2m}} q 
        &= \sum_{k=0}^m \binom{m-\frac12}{m-k} (-4)^{m-k} 
        \Bigg( \sum_{j=0}^{k - 2} J_{2k+1,j} q^2 (i \del_x)^{2k - 2 - 2 j} \conj{q} 
        - \sum_{j=0}^{k - 1} K_{2k+1,j} (i \del_x)^{2 k - 2 j} q \Bigg)
        \\ \nonumber &+ \sum_{k=0}^m \binom{m-\frac12}{m-k} (-4)^{m-k} 
        \Bigg( \sum_{j=0}^{k - 2} J_{2k+1,j} (|q|^{2 j} - 1) q^2 (i \del_x)^{2k - 2 - 2 j} \conj{q} 
        \\ \nonumber &- \sum_{j=0}^{k - 1} K_{2k+1,j} (|q|^{2 j} - 1) (i \del_x)^{2 k - 2 j} q 
        + (k+1) C_k |q|^{2k} q \Bigg) + \mc{O}^{2,n-2}_q(q_x)
        \\ \label{eqn:struct-GP-2} i \del_{t_{2m+1}} q 
        &= \sum_{k=0}^m \binom{m-1}{m-k} (-4)^{m-k} 
        \Bigg( \sum_{j=0}^{k-1} K_{2k,j} q^j \conj{q}^j (i \del_x)^{2 k - 1 - 2 j} q \Bigg)
        \\ \nonumber &+ \sum_{k=0}^m \binom{m-1}{m-k} (-4)^{m-k} 
        \Bigg( \sum_{j=0}^{k-1} K_{2k,j} (|q|^{2 j} - 1) (i \del_x)^{2 k - 1 - 2 j} q \Bigg) 
        + \mc{O}^{2,n-2}_q(q_x) \,.
    \end{align}
    Note that all terms with a derivative on $q$ or the factor $|q|^{2j} - 1$ are in $\mc{O}^{2,n-2}_{q, |q|^2 - 1}(q_x, |q|^2 - 1)$.
    We calculate
    \begin{align*}
        & \sum_{m=0}^\infty \sum_{k=0}^m X^m \binom{m-\frac12}{m-k} (-4)^{m-k} (k+1) C_k |q|^{2k} q
        = \left( 1 - 4 X (|q|^2 - 1) \right)^{- \frac12} q
        \\ &= \sum_{m=0}^\infty X^m \binom{2m}{m} (|q|^2 - 1)^m q
        = \sum_{m=0}^\infty X^m \left( \mathds{1}_{\{m=0\}} q + 2 \mathds{1}_{\{m=1\}} (q^2 \conj{q} - q) + \mc{O}^{2,0}_{|q|^2 - 1}(|q|^2 - 1) \right) \,.
    \end{align*}
    It remains to simplify the sums in the first lines of \eqref{eqn:struct-GP-1} and \eqref{eqn:struct-GP-2}. Swapping the order of summation, the goal is to evaluate
    \begin{align*}
        &\sum_{j=1}^{m-1} q^2 (i \del_x)^{2j} \conj{q} \sum_{k=j+1}^m \binom{m-\frac12}{m-k} (-4)^{m-k} J_{2k+1,k-1-j} 
        - \sum_{j=1}^m (i \del_x)^{2j} q \sum_{k=j}^m \binom{m-\frac12}{m-k} (-4)^{m-k} K_{2k+1,k-j} \\
        \text{ and } &\sum_{j=0}^{m-1} (i \del_x)^{2j+1} q \sum_{k=j+1}^m \binom{m-1}{m-k} (-4)^{m-k} K_{2k,k-1-j} \,.
    \end{align*}
    We calculate
    \begin{align*}
        \sum_{k=j+1}^m \binom{m-\frac12}{m-k} (-4)^{m-k} J_{2k+1,k-1-j}
        &= 2 \mathds{1}_{\{j = m-1\}}
    \end{align*}
    and
    \begin{align*}
        - \sum_{k=j}^m \binom{m-\frac12}{m-k} (-4)^{m-k} K_{2k+1,k-j}
        &= 4^{m-j} \binom{\frac12}{m-j}
        - 2 \sum_{l=0}^{m-2-j} (-1)^l (C_{l+1} - 2 C_l) 4^{m-2-j-l} \binom{-\frac12}{m-2-j-l} \,.
    \end{align*}
    Studying the generating function
    \begin{align*}
        \sum_{k=-2}^\infty Y^k 4^{k+2} \binom{\frac12}{k+2} 
        - 2 \sum_{k=-2}^\infty Y^k \sum_{l=0}^k (-1)^l (C_{l+1} - 2 C_l) 4^{k-l} \binom{-\frac12}{k-l}
        &= \frac{2}{Y} + \frac{1}{Y^2} \,,
    \end{align*}
    we find that
    \begin{align*}
        - \sum_{k=j}^m \binom{m-\frac12}{m-k} (-4)^{m-k} K_{2k+1,k-j}
        &= \mathds{1}_{\{j = m\}} + 2 \mathds{1}_{\{j = m-1\}} \,.
    \end{align*}
    Lastly, we calculate
    \begin{align*}
        \sum_{k=j+1}^m \binom{m-1}{m-k} (-4)^{m-k} K_{2k,k-1-j}
        &= 4^{m-1-j} \binom{\frac12}{m-1-j} \,.
    \end{align*}
    The results of our calculations now imply \eqref{eqn:struct-GP-3} and \eqref{eqn:struct-GP-4}.
\end{proof}

\begin{proof}[Proof of Proposition \ref{prop:2}]
    We plug the Ansatz $q(t, x) = q_\ast(x) + p(t, x)$ into \eqref{eqn:GP} for $m = 1$ and \eqref{eqn:struct-GP-3}--\eqref{eqn:struct-GP-4} for $m \geq 2$. 
    First, observe that
    \begin{align*}
        \mc{O}^{2,n-2}_{q, |q|^2 - 1}(q_x, |q|^2 - 1)
        &\subseteq \mc{O}^{2,n-2}_{p, q_\ast, |q_\ast|^2 - 1}(p, (q_\ast)_x, |q_\ast|^2 - 1) \,,
    \end{align*}
    and define the shorthand notation
    \begin{align*}
        \bm{O}_n &= \mc{O}^{1,n}_{q_\ast}((q_\ast)_x) + \mc{O}^{2,n-2}_{p, q_\ast, |q_\ast|^2 - 1}(p, (q_\ast)_x, |q_\ast|^2 - 1) \,.
    \end{align*}
    We can already deduce \eqref{eqn:lin-odd} from \eqref{eqn:struct-GP-4} by inverting the sequence convolution, using the Chu-Vandermonde identity.
    Observe now that for all $m \geq 1$ we have
    \begin{align*}
        (i \del_x)^2 q + 2 (|q|^2 - 1) q &= (D_x^2 + 2) p + 2 q_\ast^2 \conj{p} + \mc{O}^{1,0}_{q_\ast, |q_\ast|^2 - 1} (|q_\ast|^2 - 1) + \bm{O}_2
        \\ ((i \del_x)^{2m} + 2 (i \del_x)^{2m-2}) q + 2 q^2 (i \del_x)^{2m-2} \conj{q} &= \left( D_x^{2m} + 2 D_x^{2m-2} \right) p + 2 q_\ast^2 D_x^{2m-2} \conj{p} + \bm{O}_n \,.
    \end{align*}
    To obtain, \eqref{eqn:lin-even} we calculate
    \begin{align*}
        i \del_{t_{2m}} \begin{pmatrix} 
            p \\ q_\ast^2 \conj{p} \\ \conj{p} \\ \conj{q}_\ast^2 p
        \end{pmatrix}
        &= \begin{pmatrix}
            \left( D_x^{2m} + 2 D_x^{2m-2} \right) p + 2 q_\ast^2 D_x^{2m-2} \conj{p}
            \\ - q_\ast^2 \left( D_x^{2m} + 2 D_x^{2m-2} \right) \conj{p} - 2 |q_\ast|^2 D_x^{2m-2} p
            \\ - \left( D_x^{2m} + 2 D_x^{2m-2} \right) \conj{p} - 2 \conj{q}_\ast^2 D_x^{2m-2} p
            \\ \conj{q}_\ast^2 \left( D_x^{2m} + 2 D_x^{2m-2} \right) p + 2 |q_\ast|^2 D_x^{2m-2} \conj{p}
        \end{pmatrix} + \bm{O}_n
        \\ &= \begin{pmatrix}
            \left( D_x^{2m} + 2 D_x^{2m-2} \right) p + 2 D_x^{2m-2} (q_\ast^2 \conj{p})
            \\ - \left( D_x^{2m} + 2 D_x^{2m-2} \right) (q_\ast^2 \conj{p}) - 2 D_x^{2m-2} p
            \\ - \left( D_x^{2m} + 2 D_x^{2m-2} \right) \conj{p} - 2 D_x^{2m-2} (\conj{q}_\ast^2 p)
            \\ \left( D_x^{2m} + 2 D_x^{2m-2} \right) (\conj{q}_\ast^2 p) + 2 D_x^{2m-2} \conj{p}
        \end{pmatrix} + \bm{O}_n \,.
    \end{align*}
\end{proof}

\section{Well-posedness results}
\label{section:4}

\subsection{Local well-posedness for a large class of dispersive nonlinear systems}
As mentioned in the introduction, the well-posedness theory presented here builds on works by C.E. Kenig, G. Ponce, and L. Vega from the 1990s (see \cite{KenigPonceVega1991-1,KenigPonceVega1991-2,KenigPonceVega1993,KenigPonceVega-SmallSolutions,KenigPonceVega-BenjaminOno}), specifically \cite{KenigPonceVega1994-2}.

Let $n, D, K \in \N$ with $n \geq 2$ and set $m = \frac{n-1}{2} \geq \frac12$. 
For $u = u(t, x) : \R^2 \rightarrow \C^D$, we consider the system
\begin{align} \label{eqn:dispersive}
    \del_t u &= i \Phi(D_x) u + \mc{N}[u] \,.
\end{align}
Here
\begin{align} \label{eqn:dispersive-nonlinear}
    (\mc{N}[u])_d &= g_d + \sum_{k=1}^K \sum_{b \in \{1,\dots,D\}^k} \sum_{l \in \N^k, |l| \leq n - 2} f^{k,l}_{d,b} \mc{N}^{k,l}_b[u] 
    & \mc{N}^{k,l}_b[u] &= \prod_{j=1}^k \del_x^{l_j} u_{b_j}
\end{align}
for $1 \leq d \leq D$, with $g_d, f^{k,l}_{d,b}: \R \rightarrow \C$ being space-dependent coefficients,
and there exist measurable functions $\mc{V}, \mc{D}, \mc{W} = \mc{V}^{-1}: \R \rightarrow \R^{D \times D}$ such that
\begin{align*}
    \Phi = \mc{V} \mc{D} \mc{W} = \begin{pmatrix}
        \mc{V}_{1,1} & \cdots & \mc{V}_{D,1}
        \\ \vdots & \ddots & \vdots
        \\ \mc{V}_{1,D} & \cdots & \mc{V}_{D,D}
    \end{pmatrix} \begin{pmatrix}
        \phi_1 & \cdots & 0
        \\ \vdots & \ddots & \vdots
        \\ 0 & \cdots & \phi_D
    \end{pmatrix} \begin{pmatrix}
        \mc{W}_{1,1} & \cdots & \mc{W}_{D,1}
        \\ \vdots & \ddots & \vdots
        \\ \mc{W}_{1,D} & \cdots & \mc{W}_{D,D}
    \end{pmatrix} \,.
\end{align*}
Instead of diagonalizing $\Phi$ by the change of variables $u \mapsto \mc{V} u$, we shall use the relation
\begin{align} \label{eqn:exp-decomp}
    \left( e^{i t \Phi(D_x)} \right)_{p,q} &= \sum_{d=1}^D \mc{V}_{d,q}(D_x) \mc{W}_{p,d}(D_x) e^{i t \phi_d(D_x)} \,.
\end{align}
\begin{definition}
    Consider an open interval $I \subset \R$ and let $\alpha \in (0, 1)$. For a function $\phi \in C^1(I; \R)$ and a point $\xi_0 \in I$, we say that $\phi$ has a \textbf{critical point with steepness $\alpha$ at $\xi_0$} 
    if $\phi'(\xi_0) = 0$, and there exist $r > 0$ and $h_- \in C^1([\xi_0 - r, \xi_0]; \R)$, $h_+ \in C^1([0, \xi_0 + r]; \R)$ such that $h_\pm(0) \neq 0$ and
    \begin{align*}
        \phi(\xi) &= \phi(\xi_0) + |\xi - \xi_0|^{\frac{1}{\alpha}} h_\pm(\xi) \qquad \qquad \forall\, \xi \in \R \text{ s.t. } \pm (\xi - \xi_0) \in (0, r) \,.
    \end{align*}
    We say that $\phi$ has a \textbf{critical point with positive steepness at $\xi_0$} if this is true for some $\alpha \in (0, 1)$.
\end{definition}
Our results require $\Phi: \R \rightarrow \C^{D \times D}$ to be continuous and satisfy certain additional assumptions, which we list below.
We choose some $\mu \in \N$ with $\mu \leq 3 m$ and define
\begin{align}
    \label{eqn:a1} \ti{a}_{d,p,q} &= \frac{|\xi|^m}{|\phi_d'(\xi)|^{\frac12}} \frac{|\xi|^\mu}{\langle \xi \rangle^\mu} \mc{V}_{d,q} \mc{W}_{p,d}
    \\ \label{eqn:a2} a_{d,p,q} &= \frac{|\xi|^{2 m}}{\phi_d'(\xi)} \frac{|\xi|^\mu}{\langle \xi \rangle^\mu} \mc{V}_{d,q} \mc{W}_{p,d} \,.
\end{align}
For every $(\phi, \ti{a}, a) \in \{(\phi_d, \ti{a}_{d,p,q}, a_{d,p,q}): d,p,q \in \{1, \dots, D\}\}$, we require the following to hold true:
\begin{enumerate}[(P1)]
    \item $\phi \in C^3(\R; \R)$ has finitely many critical points, all of which have positive steepness.
    \item There exist some $R > 0$ and $C, c > 0$ such that for all $k \in \{0, 1, 2, 3\}$ and $\xi \in \R$ with $|\xi| > R$, we have
    \begin{align*} 
        c |\xi|^{n-k} &\leq |\del_\xi^k \phi(\xi)| \leq C |\xi|^{n-k} \,.
    \end{align*}
    \item We have $a \in C^{0,1}(\R; \R)$ and there exist $C, \delta > 0$ such that
    \begin{align*}
        |\ti{a}(\xi)| + |a(\xi)| + \langle \xi \rangle^{1 + \delta} |a'(\xi)| &\leq C \,.
    \end{align*}
    \item There exists $M > 0$ such that for any $\tau \in \R$ we can decompose $\R$ into $N$ intervals on whose interiors the function $\frac{a(\xi)}{\phi(\xi) - \tau}$ is monotonic.
\end{enumerate}
Note that \eqref{eqn:exp-decomp} together with \eqref{eqn:a2}, (P1)--(P3) and continuity of $\Phi$ imply that there exists some $C > 0$ such that
\begin{align} \label{eqn:500}
    |\Phi'(\xi)| &\leq C \langle \xi \rangle^{n-1} \,.
\end{align}
\begin{lemma}
    If (P1)--(P4) are fulfilled, then the following linear estimates are available.
    \begin{enumerate}[(i)]
        \item The local smoothing estimate (compare to \cite[Theorem 4.1]{KenigPonceVega1991-1}, \cite[Theorem 2.1, Corollary 2.2]{KenigPonceVega1994-1}, \cite[Theorem 3.5]{KenigPonceVega1993}).
        Define $\alpha_h = \frac{1 - \theta_h}{2}$ and $\frac{1}{p_h} = \frac{\theta_h}{1} + \frac{1 - \theta_h}{2}$, where 
        \begin{align*}
            \theta_h = \begin{cases}
                \frac{h}{m + \mu} &, 0 \leq h \leq m + \mu
                \\ \frac{h - m - \mu}{m} &, m + \mu \leq h \leq 2 m + \mu
            \end{cases} \,.
        \end{align*}
        Then for all $h \in \N$ with $0 \leq h \leq m + \mu$ we have
        \begin{align}
            \label{eqn:501} \left\| |D_x|^{m+\mu} e^{i t \Phi(D_x)} u_0 \right\|_{L^\infty_x L^2_{t \in [-T, T]}} 
            &\lesssim \|\langle D_x \rangle^\mu u_0\|_{L^2_x}
            \\ \label{eqn:503} \left\| \int_0^t D_x^h e^{i (t - t') \Phi(D_x)} u(t', x) \dd t' \right\|_{L^\infty_{t \in [-T, T]} L^2_x} 
            &\lesssim T^{\alpha_h} \|\langle D_x \rangle^{\theta_h \mu} u\|_{L^{p_h}_x L^2_{t \in [-T, T]}}
            \\ \label{eqn:511} \left\| \int_0^t D_x^h e^{i (t - t') \Phi(D_x)} u(t', x) \dd t' \right\|_{L^\infty_x L^2_{t \in [-T, T]}} 
            &\lesssim T^{\alpha_h} \|\langle D_x \rangle^\mu u\|_{L^{p_h}_x L^2_{t \in [-T, T]}}
        \end{align}
        \item The maximal function estimate (compare to \cite[Theorem 2.5]{KenigPonceVega1991-1}, \cite[Theorem 2.3]{KenigPonceVega1994-1}, \cite[Corollary 2.9]{KenigPonceVega1991-2}).
        Set $r = \frac12 \lor \frac{m}{2}$. We have
        \begin{align}
            \label{eqn:502} \left\| e^{i t \Phi(D_x)} u_0 \right\|_{L^2_x L^\infty_{t \in [-T, T]}} &\lesssim_T \|u_0\|_{H^r_x}
            \\ \label{eqn:504} \left\| \Phi'(D_x) e^{i t \Phi(D_x)} u_0 \right\|_{L^2_x L^\infty_{t \in [-T, T]}} &\lesssim_T \|u_0\|_{H^{r+n-1}_x} \,.
        \end{align}
    \end{enumerate}
\end{lemma}
\begin{proof}
    \begin{enumerate}[(i)]
        \item We want to apply Theorem \ref{thm:thm3.4}. This requires $\ti{a} \in L^\infty$, which we assumed in (P3), and $a$ and $\phi$ to fulfill (H1)--(H5).
        We trivially obtain (H1) from (P1) and adopt the definitions given there. 
        We must show that (H2) also follows from (P1). 
        Here we only write the proof of the required estimates for $\psi_j$ and $\psi_0$.
        Let $\alpha$, $r$ and $h_\pm$ be given by the definition of positive steepness at $\xi_j$, and set $\xi_j = 0$, $\eta_j = 0$, $\sigma_j = 1$ without loss of generality.
        Furthermore, it suffice to consider $\eta \in (0, \phi(r))$. We have
        \begin{align*}
            \eta &= h_+(\psi_j(\eta)) |\psi_j(\eta)|^{\frac{1}{\alpha}}
            \\ 1 &= h_+'(\psi_j(\eta)) |\psi_j(\eta)|^{\frac{1}{\alpha}} \psi_j'(\eta) + h_+(\psi_j(\eta)) \frac{1}{\alpha} |\psi_j(\eta)|^{\frac{1}{\alpha}-1} \sign(\psi_j(\eta)) \psi_j'(\eta) \,.
        \end{align*}
        This implies
        \begin{align*}
            \frac{|\eta|^\alpha}{|\psi_j(\eta)|} &= |h_+(\psi_j(\eta))|^\alpha
            \\ \frac{|\eta|^{\alpha - 1}}{|\psi_j'(\eta)|} 
            &= |\eta|^{\alpha - 1} \left| h_+'(\psi_j(\eta)) |\psi_j(\eta)|^{\frac{1}{\alpha}} + |h_+(\psi_j(\eta)) \frac{1}{\alpha} |\psi_j(\eta)|^{\frac{1}{\alpha}-1} \sign(\psi_j(\eta)) \right| \,.
        \end{align*}
        Since $h_+ \circ \psi_j$ is non-zero and continuous to the right of $\eta = 0$, \eqref{eqn:ass-1} and then also \eqref{eqn:ass-2} follow. 
        Applying the estimates from (P2) to $\xi = \psi_j(\eta)$ where $j \in \{0, N\}$ yields (H3) with $\beta = \frac{1}{n}$.
        Finally (P3)--(P4) are just (H4)--(H5) for a specific choice of $a$.
        
        We obtain \eqref{eqn:501} by diagonalizing the exponential with \eqref{eqn:exp-decomp} and applying \eqref{eqn:smooth-1} to $\langle D_x \rangle^\mu u_0$, with the choice \eqref{eqn:a1} for $\ti{a}$. 
        For \eqref{eqn:503} we first note that
        \begin{align*}
            \left\| \int_0^t |D_x|^h e^{i (t - t') \phi(D_x)} u(t', x) \dd t' \right\|_{L^\infty_{t \in [-T, T]} L^2_x} 
            &\lesssim T^{\alpha_h} \|\langle D_x \rangle^{\theta_h \mu} u\|_{L^{p_h}_x L^2_{t \in [-T, T]}} \,.
        \end{align*}
        This is shown by interpolation\footnote{The necessary interpolation inequality follows from the argument in \cite{Stein1956}, using a version of the ``three-lines lemma''.} of the case $h = 0$, which is a direct consequence of applying Hölder's inequality in $t'$, 
        with the case $h = m + \mu$, which is the dual of \eqref{eqn:501}.
        Since $h \in \N$, we may write $|D_x|^h = D_x^h H^h$, where $H$ is the Hilbert transform. 
        By \cite{Burkholder} the vector-valued Hilbert transform is bounded on $L^{p_h}_x L^2_{[-T,T]}$, so we obtain \eqref{eqn:503}.
        To prove \eqref{eqn:511}, we again first show that
        \begin{align*}
            \left\| \int_0^t |D_x|^h e^{i (t - t') \phi(D_x)} u(t', x) \dd t' \right\|_{L^\infty_x L^2_{t \in [-T, T]}} 
            &\lesssim T^{\alpha_h} \|\langle D_x \rangle^{\theta_h \mu} u\|_{L^{p_h}_x L^2_{t \in [-T, T]}}
        \end{align*}
        by interpolation of the case $h = m + \mu$, which is shown by combining Minkowski's inequality with \eqref{eqn:501}, 
        with the case $h = 2 m + \mu$, which is \eqref{eqn:exp-decomp} combined with \eqref{eqn:smooth-3} and the choice \eqref{eqn:a2} for $a$.
        The claimed estimate follows again from the boundedness of the Hilbert transform.

        \item It suffices to prove \eqref{eqn:502}, because together with \eqref{eqn:500} it implies \eqref{eqn:504}. Let $\psi_0 \in C_c^\infty(\R; \R)$. 
        We estimate the low frequency part $\psi_0(D_x) u_0$ by applying Theorem \ref{thm:cor2.8} with $a = 0$ and $b(t, D_x) = e^{i t \Phi(D_x)} \psi_0(D_x)$.
        For the high frequency part $(1 - \psi_0(D_x)) u_0$ we use \eqref{eqn:exp-decomp} to reduce to the scalar case with $\phi = \phi_d$, and set $a = \mc{V}_{d,q} \mc{W}_{p,d}$, $b = 0$.
        The requirements (J1)--(J6) for $\phi$ are all direct consequences of (P2). The boundedness of $a$ also follows from (P2), together with the boundedness of \eqref{eqn:a2} due to (P3).
    \end{enumerate}
\end{proof}
From now on we assume $D = 1$, which allows us to drop the (multi-)indices $d$ and $b$ from the notation.
Specifically, $u_d$, $\mc{N}^{k,l}_b[u]$, $g_d$ and $f^{k,l}_{d,b}$, and are replaced by $u$, $\mc{N}^{k,l}[u]$, $g$ and $f^{k,l}$.
This is possible because there is no relevant interplay between the components; the presented arguments transfer directly to the case $D > 1$.

For the rest of the section follows, we always consider $h \in \N$ so that $\del_x^h$ is well-defined. 
Let $T > 0$ and consider the (pseudo-)norms
\begin{align*}
    \|u\|_{Y_{1,T}} &= \max_{0 \leq h \leq s} \|\del_x^h u\|_{L^\infty_{t \in [-T, T]} L^2_x}
    & \|u\|_{Y_{2,T}} &= \max_{s \leq h \leq s + m} \|\del_x^h u\|_{L^\infty_x L^2_{t \in [-T, T]}}
    \\ \|u\|_{Y_{3,T}} &= \max_{0 \leq h \leq s - m - r} \|\del_x^h u\|_{L^2_x L^\infty_{t \in [-T, T]}}
    & \|u\|_{Y_{4,T}} &= \max_{0 \leq h \leq s' - m - r} \|x \del_x^h u\|_{L^2_x L^\infty_{t \in [-T, T]}}
    \\ \|u\|_{Y_{5,T}} &= \max_{0 \leq h \leq s'} \|x \del_x^h u\|_{L^\infty_{t \in [-T, T]} L^2_x}
    & \|u\|_{Y_T} &= \sum_{j=1}^5 \|u\|_{Y_{j,T}} \,.
\end{align*}
Our fixed-point argument uses the Banach space
\begin{align*}
    \left( Y_T, \|\cdot\|_{Y_T} \right) = \left( \left\{ u \in C_{t \in [-T, T]} (H^s \cap H^{s',1})_x: \|u\|_{Y_T} < \infty \right\}, \|\cdot\|_{Y_T} \right) \,.
\end{align*}
\begin{definition}[Mild solution]
    For a given $u_0 \in H^s \cap H^{s',1}$ and a time $T > 0$ we say that \textbf{$u \in Y_T$ is a mild solution with initial data $u_0$ to \eqref{eqn:dispersive}}
    if
    \begin{align} \label{eqn:mild}
        \Lambda_{u_0}[u] &= u \qquad \text{ where } \qquad \Lambda_{u_0}[u](t) = e^{i t \Phi(D_x)} u_0 + \int_0^t e^{i (t - t') \Phi(D_x)} \mc{N}[u(t')] \dd t' \,.
    \end{align}
\end{definition}
\begin{restatable}[Local well-posedness]{thm}{thmlwp} \label{thmlwp}
    Let $s, s' \in \N$ fulfill $s \geq 5 m$ and $\frac{s+m}{2} \leq s' \leq s - 2 m$.
    Assume 
    \begin{align} \label{eqn:ass-coeff-1}
        C_{f,g} &= \|g\|_{H^{s+1-\ceil{m},1}} + \sum_{l=0}^{n-2} \|f^{1,l}\|_{H^{s+1-\ceil{m},1}} + \sum_{k=1}^K \sum_{l \in \N^k, |l| \leq n - 2} \|f^{k,l}\|_{W^{s+1-\ceil{m},\infty}} < \infty \,.
    \end{align}
    For any $u_0 \in H^s \cap H^{s',1}$ there exists a time $T = T(\|u_0\|_{H^s \cap H^{s',1}}) > 0$ for which \eqref{eqn:dispersive} has a mild solution in $Y_T$. 
    This solution is unique in the sense that it agrees with any other mild solution in $Y_T$. The data-to-solution map is locally Lipschitz continuous in the sense that for two initial data $u_0, \ti{u}_0 \in H^s \cap H^{s',1}$, with times of existence $T$ and $\ti{T}$ respectively, we have
    \begin{align*}
        \|u - \ti{u}\|_{Y_{\min\{T, \ti{T}\}}} &\leq C(\|u_0\|_{H^s \cap H^{s',1}}, \|\ti{u}_0\|_{H^s \cap H^{s',1}}) \|u_0 - \ti{u}_0\|_{H^s \cap H^{s',1}} \,.
    \end{align*}
\end{restatable}
\begin{proof}
    We decompose
    \begin{align*}
        \Lambda_{u_0}[u] &= e^{i t \Phi(D_x)} u_0
        + \sum_{k=0}^K \sum_{|l| \leq n - 2} \int_0^t e^{i (t - t') \Phi(D_x)} \left( f^{k,l} \mc{N}^{k,l}[u(t')] \right) \dd t'
        \\ &= (I) + (II) \,.
    \end{align*}
    We fix some $T > 0$ to be chosen later and start with the estimates for $(I)$:
    \begin{align*}
        \|e^{i t \Phi(D_x)} u_0\|_{Y_{1,T}} &= \max_{0 \leq h \leq s} \|\del_x^h e^{i t \Phi(D_x)} u_0\|_{L^\infty_{t \in [-T, T]} L^2_x} \leq \max_{0 \leq h \leq s} \|\del_x^h u_0\|_{L^2_x} \lesssim \|u_0\|_{H^s_x}
        \\ \|e^{i t \Phi(D_x)} u_0\|_{Y_{2,T}} &= \max_{s \leq h \leq s + m} \|\del_x^h e^{i t \Phi(D_x)} u_0\|_{L^\infty_x L^2_{t \in [-T, T]}} \overset{\eqref{eqn:501}}{\lesssim} \max_{0 \leq h \leq s} \||D_x|^h u_0\|_{L^2_x} \lesssim \|u_0\|_{H^s_x}
        \\ \|e^{i t \Phi(D_x)} u_0\|_{Y_{3,T}} &= \max_{0 \leq h \leq s - m - r} \|\del_x^h e^{i t \Phi(D_x)} u_0\|_{L^2_x L^\infty_{t \in [-T, T]}} \overset{\eqref{eqn:502}}{\lesssim_T}  \max_{0 \leq h \leq s - m - r} \||D_x|^h u_0\|_{H^r_x} \lesssim \|u_0\|_{H^s_x} \,.
    \end{align*}
    For $Y_{4,T}$ we note that
    \begin{align} \label{eqn:512}
        x e^{i t \Phi(D_x)} u &= \mc{F}^{-1}\left[ i \del_\xi \big( e^{i t \Phi(\xi)} \ha{u} \big) \right] = \mc{F}^{-1}\left[ e^{i t \Phi(\xi)} (i \del_\xi - t \Phi'(\xi)) \ha{u} \right] = e^{i t \Phi(D_x)} (x u - t \Phi'(D_x) u) \,,
    \end{align}
    and hence
    \begin{align*}
        \|e^{i t \Phi(D_x)} u_0\|_{Y_{4,T}} &= \max_{0 \leq h \leq s' - m - r} \|x \del_x^h e^{i t \Phi(D_x)} u_0\|_{L^2_x L^\infty_{t \in [-T, T]}}
        \\ &\leq \max_{0 \leq h \leq s' - m - r} \|\del_x^h e^{i t \Phi(D_x)} (x u_0)\|_{L^2_x L^\infty_{t \in [-T, T]}} 
        + \|\del_x^h e^{i t \Phi(D_x)} t \Phi'(D_x) u_0\|_{L^2_x L^\infty_{t \in [-T, T]}} 
        \\ &\quad + \|e^{i t \Phi(D_x)} u_0\|_{Y_{3,T}} \,.
    \end{align*}
    We estimate
       \begin{align*}
        \max_{0 \leq h \leq s' - m - r} \|\del_x^h e^{i t \Phi(D_x)} t \Phi'(D_x) u_0\|_{L^2_x L^\infty_{t \in [-T, T]}} 
        &\overset{\eqref{eqn:504}}{\lesssim_T} \|u_0\|_{H^{s' + n - 1}_x} \lesssim \|u_0\|_{H^s_x} \,.
        \\ \max_{0 \leq h \leq s' - m - r} \|\del_x^h e^{i t \Phi(D_x)} (x u_0)\|_{L^2_x L^\infty_{t \in [-T, T]}} 
        &\overset{\eqref{eqn:502}}{\lesssim_T} \|x u_0\|_{H^{s'}_x} \,.
    \end{align*}
    For $Y_{5,T}$ we use the same strategy and obtain
    \begin{align*}
        \max_{0 \leq h \leq s'} \|\del_x^h e^{i t \Phi(D_x)} (x u_0)\|_{L^\infty_{t \in [-T, T]} L^2_x} 
        &\leq \|x u_0\|_{H^{s'}} 
        \\ \max_{0 \leq h \leq s'} \|\del_x^h e^{i t \Phi(D_x)} t \Phi'(D_x) u_0\|_{L^\infty_{t \in [-T, T]} L^2_x} 
        &\overset{\eqref{eqn:500}}{\lesssim_T} \|u_0\|_{H^{s' + n - 1}} \,.
    \end{align*}
    In summary,
    \begin{align*}
        \|e^{i t \Phi(D_x)} u_0\|_{Y_T} &\lesssim_T \|u_0\|_{H^s} + \|u_0\|_{H^{s',1}} \,.
    \end{align*}
    We now estimate $(II)$. Here we always assume that $l$ is a multiindex with $l_1 \geq \dots \geq l_k$.
    In the proof below, we write everything as if $k \geq 2$. For the case $k = 1$ the argument is identical, but each $\del_x^{l_j} u$ with $j > k$ must be replaced by $1$.
    Similarly, for the term with $g$ we replace every $\del_x^{l_j} u$ by $1$.
    The factor $\del_x^{l_2} u$ is the only one which we control with the $L^2_x$-norm, i.e. the only one that requires integrability. 
    We shall therefore always group (derivatives of) the coefficients $\del_x^h g$ and $\del_x^h f^{1,l}$ with $\del_x^{l_2} u$. 
    When $k \geq 2$ we place $f^{k,l}$ in an $L^\infty$-norm and keep $\del_x^{l_2} u$ in the $L^2$-norm. 
    
    Let us first state for all $h \leq s - 1$ the inequalities
    \begin{align}
        \label{eqn:520} \|\del_x^h u\|_{L^2_{t \in [-T, T]} L^2_x} \leq T^\frac12 \|\del_x^h u\|_{L^\infty_{t \in [-T, T]} L^2_x} \leq T^{\frac12} \|u\|_{Y_{1,T}}
        \\ \label{eqn:521} \|\del_x^h u\|_{L^\infty_{t \in [-T, T]} L^\infty_x} \lesssim \|\del_x^h u\|_{L^\infty_{t \in [-T, T]} L^2_x} + \|\del_x^{h+1} u\|_{L^\infty_{t \in [-T, T]} L^2_x} \lesssim \|u\|_{Y_{1,T}} \,.
    \end{align}
    \begin{enumerate}[(1)]
        \item Estimates for $Y_{1,T}$. We assume $|l| \leq n - 2$ and $h \leq s$, and estimate
        \begin{align*}
            \left\| \int_0^t e^{i (t - t') \Phi(D_x)} \del_x^h \left( f^{k,l} \mc{N}^{k,l}[u] \right) \dd t' \right\|_{L^\infty_{t \in [-T, T]} L^2_x} 
            &\overset{\eqref{eqn:503}}{\lesssim_T} T^\alpha \left\| \del_x^{h - (\ceil{m+\mu-1} \land h)} \langle D_x \rangle^{\theta \mu} \left( f^{k,l} \mc{N}^{k,l}[u] \right) \right\|_{L^p_x L^2_{t \in [-T, T]}} \,.
        \end{align*}
        Here $\theta = \theta_{\ceil{m+\mu-1} \land h}$, and $\alpha > 0$ and $p \in (1, 2)$ depend on $\ceil{m+\mu-1} \land h$, which we suppress from our notation.
        Note that we must remove an integer number of derivatives strictly less than $m + \mu$ and not exceeding $h$, i.e. $\ceil{m+\mu-1} \land h$, in order to have $\alpha > 0$.
        Since we do not want to use fractional Leibniz inequalities, we shall later replace $\langle D_x \rangle^{\theta \mu}$ by up to $\ceil{\theta \mu}$ derivatives.
        Then the total number of derivatives present is bounded by
        \begin{align*}
            & \sup_{0 \leq h \leq s} h - (\ceil{m + \mu - 1} \land h) + \ceil{\frac{\ceil{m + \mu - 1} \land h}{m + \mu} \mu} + 2 m - 1
            \\ &\leq s - (\ceil{m} + \mu - 1) \land s + 2 m - 1 + \ceil{\frac{(\ceil{m} + \mu - 1) \land s}{m + \mu} \mu}
            \leq s + m \,.
        \end{align*}
        Here we have used that $\ceil{m} + \mu - 1 \leq \ceil{5 m} \leq \ceil{s} = s$.
        Distributing the derivatives, it suffices to assume $|l| \leq s - \ceil{m+\mu-1} + 2 m - 1 \leq s + m - \ceil{\theta \mu}$ and $h \leq s + 1 - m - \ceil{\theta \mu}$ for some $\theta \in [0, 1)$, and find an estimate for $\langle D_x \rangle^{\theta \mu} (\del_x^h f^{k,l} \mc{N}^{k,l}[u])$.
        Applying Hölder's inequality with $\frac{1}{\frac{2}{p}} + \frac{1}{\frac{2}{2-p}} = 1$ yields
        \begin{align*}
            & \left\| \langle D_x \rangle^{\theta \mu} \left( \del_x^h f^{k,l} \mc{N}^{k,l}[u] \right) \right\|_{L^p_x L^2_{t \in [-T, T]}} 
            \\ &\lesssim \left\| \langle D_x \rangle^{\theta \mu} \left( \del_x^h f^{k,l} \mc{N}^{k,l}[u] \right) \right\|_{L^2_x L^2_{t \in [-T, T]}} + \left\| x \langle D_x \rangle^{\theta \mu} (\del_x^h f^{k,l} \mc{N}^{k,l}[u]) \right\|_{L^2_x L^2_{t \in [-T, T]}} \left\| \frac{1}{x} \right\|_{L^{\frac{2p}{2-p}}_{x \not\in [-1,1]} L^\infty_{t \in [-T, T]}}
            \\ &\lesssim \sum_{j=0}^{\ceil{\theta \mu}} \left\| \del_x^j \left( \del_x^h f^{k,l} \mc{N}^{k,l}[u] \right) \right\|_{L^2_x L^2_{t \in [-T, T]}} 
            + \left\| x \del_x^j\left( \del_x^h f^{k,l} \mc{N}^{k,l}[u] \right) \right\|_{L^2_x L^2_{t \in [-T, T]}}
            \\ &= (I) + (II) \,.
        \end{align*}
        It now suffices to assume $|l| + h \leq s + m$ and $h \leq s + 1 - \ceil{m}$, and estimate $\del_x^h f^{k,l} \mc{N}^{k,l}[u]$ and $x \del_x^h f^{k,l} \mc{N}^{k,l}[u]$. We first focus on $(I)$.
        If $l_1 \leq s - 1$ then \eqref{eqn:520} and \eqref{eqn:521} imply
        \begin{align*}
            \left\| \del_x^h f^{k,l} \mc{N}^{k,l}[u] \right\|_{L^2_{t \in [-T, T]} L^2_x} &\lesssim_T C_{f,g} \prod_{j=1}^k \|\del_x^{l_j} u\|_{Y_{1,T}} \,.
        \end{align*}
        If $l_1 \geq s$ then $l_2 \leq m \leq s - m - r$, and hence
        \begin{align*}
             \left\| \del_x^h f^{k,l} \mc{N}^{k,l}[u] \right\|_{L^2_{t \in [-T, T]} L^2_x}
             &\leq \left\| \del_x^{l_1} u \right\|_{L^\infty_x L^2_{t \in [-T, T]}} \left\| \del_x^h f^{k,l} \del_x^{l_2} u \right\|_{L^2_x L^\infty_{t \in [-T, T]}} \prod_{j=3}^k \|\del_x^{l_j} u\|_{L^\infty_x L^\infty_{t \in [-T, T]}}
            \\  &\lesssim_T C_{f,g} \left\| u \right\|_{Y_{2,T}} \left\| u \right\|_{Y_{3,T}} \prod_{j=3}^k \|u\|_{Y_{1,T}} \,.
        \end{align*}
        For $(II)$ we group the weight $x$ with $\del_x^{l_2} u$ and replace usage of $Y_{1,T}$ with $Y_{5,T}$ and $Y_{3,T}$ with $Y_{4,T}$ for that term. 
        Since $l_2 \leq l_1 \leq s + m$ we have $l_2 \leq \frac{s + m}{2} \leq s'$, which allows the usage of $Y_{5,T}$.
        In the case $l_1 \geq s$ it is possible to use $Y_{3,T}$ because $l_2 \leq m \leq s' - m - r$.

        \item Estimates for $Y_{2,T}$. We assume that $|l| \leq n - 2$ and $s \leq h \leq s + m$, and estimate
        \begin{align*}
            \left\| \int_0^t e^{i (t - t') \Phi(D_x)} \del_x^h(f^{k,l} \mc{N}^{k,l}[u]) \dd t' \right\|_{L^\infty_x L^2_{t \in [-T, T]}} 
            &\overset{\eqref{eqn:511}}{\lesssim} T^\alpha \left\| \del_x^{h - \ceil{2 m + \mu - 1}} \langle D_x \rangle^\mu (f^{k,l} \mc{N}^{k,l}[u]) \right\|_{L^p_x L^2_{t \in [-T, T]}} \,.
        \end{align*}
        Here $\theta = \theta_{\ceil{2 m + \mu - 1}}$, and $\alpha > 0$ and $p \in (1, 2)$ depend on $\ceil{2 m + \mu - 1}$.
        Note that $\theta < 1$, and so $\alpha$ is positive. Note furthermore that $\ceil{2 m + \mu - 1} \leq s \leq h$, since $\mu \leq 3 m \leq s - 2 m$, which ensures that the above is well-defined.
        We shall distribute the derivatives and proceed as in (1). This is possible because total number of derivatives we obtain is bounded by
        \begin{align*}
            s + m - (\ceil{2 m} + \mu - 1) + 2 m - 1 + \mu
            &\leq s + m \,.
        \end{align*}
        The maximum number of derivatives that can fall on $f^{k,l}$ is $n - 2 = 2 m - 1$ less this quantity, so $s + 1 - \ceil{m}$, which also matches the case (1).

        \item Estimates for $Y_{3,T}$. We assume that $|l| \leq s - m - r + n - 2$ and $h \leq s - m - r$, and estimate
        \begin{align*}
            \left\| \int_0^t e^{i (t - t') \Phi(D_x)} \del_x^h f^{k,l} \mc{N}^{k,l}[u] \dd t' \right\|_{L^2_x L^\infty_{t \in [-T, T]}} &\leq \left\| e^{i (t - t') \Phi(D_x)} \del_x^h f^{k,l} \mc{N}^{k,l}[u](t', x) \right\|_{L^1_{t' \in [-T, T]} L^2_x L^\infty_{t \in [-T, T]}}
            \\ &\overset{\eqref{eqn:502}}{\lesssim_T} \left\| \del_x^h f^{k,l} \mc{N}^{k,l}[u] \right\|_{L^1_{t \in [-T, T]} H^r_x}
            \\ &\lesssim T^{\frac12} \left\| \del_x^h f^{k,l} \mc{N}^{k,l}[u] \right\|_{L^2_{t \in [-T, T]} H^r_x} \,.
        \end{align*}
        We replace the $H^r$-norm by an $L^2$-norm and distribute the derivatives, leading to a total number of derivatives $|l| \leq s - m + n - 2 \leq s + m$ and also $h \leq s + 1 - \ceil{m}$. Then we proceed as in (1).

        \item Estimates for $Y_{4,T}$. Here we assume that $|l| \leq s' - m - r + n - 2$ and $h \leq s' - m - r$, and estimate
        \begin{align*}
            & \left\| x \int_0^t e^{i (t - t') \Phi(D_x)} \del_x^h f^{k,l} \mc{N}^{k,l}[u] \dd t' \right\|_{L^2_x L^\infty_{t \in [-T, T]}} 
            \\ &\overset{\eqref{eqn:512}}{=} \left\| \int_0^t e^{i (t - t') \Phi(D_x)} (x - (t - t') \Phi'(D_x)) \del_x^h f^{k,l} \mc{N}^{k,l}[u] \dd t' \right\|_{L^2_x L^\infty_{t \in [-T, T]}}
            \\ &\leq \left\| e^{i (t - t') \Phi(D_x)} (x - (t - t') \Phi'(D_x)) \del_x^h f^{k,l} \mc{N}^{k,l}[u](t', x) \right\|_{L^1_{t' \in [-T, T]} L^2_x L^\infty_{t \in [-T, T]}}
            \\ &\overset{\eqref{eqn:502}, \eqref{eqn:504}}{\lesssim_T} T^{\frac12} \left\| x \del_x^h f^{k,l} \mc{N}^{k,l}[u] \right\|_{L^2_{t \in [-T, T]} H^r_x}
            + T^{\frac12} \left\| \del_x^h f^{k,l} \mc{N}^{k,l}[u] \right\|_{L^2_{t \in [-T, T]} H^{r+n-1}_x} \,.
        \end{align*}
        We distribute the derivatives in the Sobolev norm and obtain a number of derivatives $|l| \leq s' - m + n - 2 + n - 1 \leq s + m$ as well as $h \leq s' + 1 - \ceil{m}$, with a possible weight $x$ in the expressions to estimate, so we can proceed as in (1).

        \item Estimates for $Y_{5,T}$. Here we assume that $|l| \leq s' - m + n - 2$ and $h \leq s' - m$, and estimate
        \begin{align*}
            & \left\| x \int_0^t e^{i (t - t') \Phi(D_x)} \del_x^h f^{k,l} \mc{N}^{k,l}[u] \dd t' \right\|_{L^\infty_{t \in [-T, T]} L^2_x}
            \\ &\overset{\eqref{eqn:512}}{=} \left\| \int_0^t e^{i (t - t') \Phi(D_x)} (x - (t - t') \Phi'(D_x)) \del_x^h f^{k,l} \mc{N}^{k,l}[u] \dd t' \right\|_{L^\infty_{t \in [-T, T]} L^2_x}
            \\ &\leq T^{\frac12} \left\| x \del_x^h f^{k,l} \mc{N}^{k,l}[u] \right\|_{L^2_{t \in [-T, T]} L^2_x} + T^{\frac32} \left\| \Phi'(D_x) \del_x^h f^{k,l} \mc{N}^{k,l}[u] \right\|_{L^2_{t \in [-T, T]} L^2_x} \,.
        \end{align*}
        Using \eqref{eqn:500} we can replace $\|\Phi'(D_x) \cdot\|_{L^2_x}$ by $\|\cdot\|_{H^{n-1}_x}$.
        We again distribute all derivatives and obtain a number of derivatives $|l| \leq s + m$ as well as $h \leq s' + 1 - \ceil{m}$, with a possible weight $x$, so we can proceed as in (1).
    \end{enumerate}
    We have shown that there exists some $\alpha > 0$ such that
    \begin{align*}
        \|e^{i t \Phi(D_x)} u_0\|_{Y_T} &\lesssim_T \|u_0\|_{H^s \cap H^{s',1}}
        \\ \|\Lambda_{u_0}[u] - e^{i t \Phi(D_x)} u_0\|_{Y_T} &\lesssim_T T^\alpha C_{f,g} (1 + \|u\|_{Y_T})^K
        \\ &\lesssim_T T^\alpha C_{f,g} \left( \|u - e^{i t \Phi(D_x)} u_0\|_{Y_T} + \|u_0\|_{H^s \cap H^{s',1}} \right)^K \,.
    \end{align*}
    Consequently, for any $R > 0$ there exists some small $T = T(\|u_0\|_{H^s \cap H^{s',1}}, R) > 0$ (depending of course also on $\Phi, s, s', m, C_{f,g}$ etc.) such that
    \begin{align*}
        \|\Lambda_{u_0}[u] - e^{i t \Phi(D_x)} u_0\|_{Y_T} < R \,.
    \end{align*}
    We now show that $\Lambda_{u_0}[u] \in C_{t \in [-T, T]} (H^s \cap H^{s',1})_x$. 
    By dominated convergence and \eqref{eqn:512}, the assumption on the initial data $u_0 \in H^s \cap H^{s',1}$ directly implies $e^{i t \Phi(D_x)} u_0 \in C_{t \in [-T, T]} (H^s \cap H^{s',1})$. 
    It remains to show continuity of the nonlinear part. Let $0 < t' < t$ and decompose
    \begin{align*}
        & \int_0^t e^{i (t - t'') \Phi(D_x)} \mc{N}[u(t'')] \dd t'' - \int_0^{t'} e^{i (t' - t'') \Phi(D_x)} \mc{N}[u(t'')] \dd t''
        \\ &= \int_{t'}^t e^{i (t - t'') \Phi(D_x)} \mc{N}[u(t'')] \dd t'' + \int_0^{t'} e^{i (t' - t'') \Phi(D_x)} \left( 1 - e^{i (t - t') \Phi(D_x)} \right) \mc{N}[u(t'')] \dd t'' \,.
    \end{align*}
    As before, it suffices for estimating the $H^s$-norm to replace $\mc{N}[u]$ by $\del_x^h f^{k,l} \mc{N}^{k,l}[u]$, where $0 \leq k \leq K$, $|l| \leq s + n - 2$, $h \leq s$. 
    Then, proceeding as we did in (1) above, we obtain
    \begin{align*}
        \left\| \int_{t'}^t e^{i (t - t'') \Phi(D_x)} \del_x^h(f^{k,l} \mc{N}^{k,l}[u(t'')]) \dd t'' \right\|_{L^2_x} 
        &\lesssim_T (t - t')^\alpha C_{f,g} \|u(t'+\cdot)\|_{Y_T}^K
    \end{align*}
    for some $\alpha > 0$, and see that this term vanishes as $t' \rightarrow t$.
    For the other term the situation is more difficult, as we want to avoid incurring additional derivatives of the form $\Phi'(D_x)$ and therefore need to use dominated convergence.
    We again proceed as in (1), but stop when the quantity that remains to be estimated is
    \begin{align*}
        & \left\| \int_0^{t'} e^{i (t' - t'') \Phi(D_x)} \left( 1 - e^{i (t - t') \Phi(D_x)} \right) \del_x^h(f^{k,l} \mc{N}^{k,l}[u(t'')]) \dd t'' \right\|_{L^2_x} 
        \\ &\lesssim_{t'} \sum_{j=0}^{\ceil{\theta \mu}} \left\| \left( 1 - e^{i (t - t') \Phi(D_x)} \right) \del_x^{h - (\ceil{m+\mu-1} \land h)} \langle D_x \rangle^{\theta \mu} \left( f^{k,l} \mc{N}^{k,l}[u] \right) \right\|_{L^2_x L^2_{t \in [-T, T]}}
        \\ &+ \sum_{j=0}^{\ceil{\theta \mu}} \left\| x \left( 1 - e^{i (t - t') \Phi(D_x)} \right) \del_x^{h - (\ceil{m+\mu-1} \land h)} \langle D_x \rangle^{\theta \mu} \left( f^{k,l} \mc{N}^{k,l}[u] \right) \right\|_{L^2_x L^2_{t \in [-T, T]}} \,.
    \end{align*}
    We now use the Plancherel theorem in $L^2_x$ and subsequently apply dominated convergence. 
    For the term with a weight $x$, this involves once more the usage of \eqref{eqn:512}.
    We can estimate the $H^{s',1}$-norm in the same way, using (5) instead of (1).
    We conclude that $\Lambda_{u_0}$ maps
    \begin{align*}
        \Lambda_{u_0}: e^{i t \Phi(D_x)} u_0 + \overline{B_R^{Y_T}} \longrightarrow e^{i t \Phi(D_x)} u_0 + \overline{B_R^{Y_T}} \,.
    \end{align*}
    With the same estimates as above, we can show furthermore that if $T$ is sufficiently small, then $\Lambda_{u_0}$ is a contraction mapping on the complete metric space $\big(e^{i t \Phi(D_x)} u_0 + \overline{B_R^{Y_T}}, \|\cdot\|_{Y_T}\big)$ and hence has a fixed point. 
    Specifically, we study the expression
    \begin{align*}
        \Lambda_{u_0}[u] - \Lambda_{u_0}[\ti{u}] &= \sum_{k=1}^K \sum_{|l| \leq n - 2} \int_0^t e^{i (t - t') \Phi(D_x)} f^{k,l} \left( \mc{N}^{k,l}[u] - \mc{N}^{k,l}[\ti{u}] \right) \dd t' \,.
    \end{align*}
    We expand the difference of products into a sum of products where exactly one factor is a difference between a component of $u$ and one of $\ti{u}$. 
    The estimates previously described then directly show that there exists some $\alpha > 0$ such that
    \begin{align*}
        \|\Lambda_{u_0}[u] - \Lambda_{u_0}[\ti{u}]\|_{Y_T} 
        &\lesssim_T T^\alpha C_{f,g} \|u - \ti{u}\|_{Y_T} (1 + \|u\|_{Y_T} + \|\ti{u}\|_{Y_T})^{K - 1}
        \\ &\leq T^\alpha C_{f,g} \|u - \ti{u}\|_{Y_T} (1 + 2 R + 2 \|u_0\|_{H^s \cap H^{s',1}})^{K - 1} \,.
    \end{align*}
    Hence, if $T$ is sufficiently small, there exists a unique mild solution $u \in e^{i t \Phi(D_x)} u_0 + \overline{B_R^{Y_T}}$.
    
    So far we have obtained uniqueness only in $e^{i t \Phi(D_x)} u_0 + \overline{B_R^{Y_T}}$, a restriction that we shall now lift.
    Let $\ti{R}$ be arbitrarily large and $\ti{u} \in e^{i t \Phi(D_x)} u_0 + \overline{B_{\ti{R}}^{Y_T}} \subset Y_T$ be another fixed point of $\Lambda_{u_0}$.
    Then $u, \ti{u} \in C_{t \in [-T, T]} (H^s \cap H^{s',1})_x$. Consider the set of times
    \begin{align*}
        \mc{T} = \left\{ T^\ast \in [-T, T]: u(t) = \ti{u}(t) \quad \forall\, t \in \big[-|T^\ast|, |T^\ast|\big] \right\} \ni 0 \,.
    \end{align*}
    Clearly $T^\ast \in \mc{T}$ implies $\big[- |T^\ast|, |T^\ast|\big] \subseteq \mc{T}$. Furthermore, since $u, \ti{u} \in C_{t \in [-T, T]} (H^s \cap H^{s',1})_x$ the set $\mc{T}$ is closed.
    Let $T^\ast \in \mc{T} \cap (-T, T)$. For small $|t|$ a calculation shows that $\Lambda_{v_0}[u(\cdot + T^\ast)](t) = u(t + T^\ast)$.
    Furthermore, defining $v_0 = u(T^\ast)$, we can find a small $T' = T(R) > 0$ such that by previous estimates
    \begin{align*}
        \|u(t + T^\ast) - e^{i t \Phi(D_x)} v_0\|_{Y_{T'}} &\lesssim (T')^\alpha C(T', R) < 1 \,,
    \end{align*}
    and hence $u(t + T^\ast) \in e^{i t \Phi(D_x)} v_0 + \overline{B_1^{Y_{T'}}}$. If $T' > 0$ is sufficiently small then $\Lambda_{v_0}$ is a contraction on $e^{i t \Phi(D_x)} v_0 + \overline{B_1^{Y_{T'}}}$ and $u(t + T^\ast)$ its unique fixed point. 
    The same holds true for $\ti{u}$ for a small $T' = T'(\|u_0\|_{H^s \cap H^{s',1}}, \ti{R}) > 0$, which implies that $u(t + T^\ast)\big\vert_{t \in [-T', T']} = \ti{u}(t + T^\ast)\big\vert_{t \in [-T', T']}$ as elements of $C_{t \in [-T', T']} H^s_x$. 
    Therefore $\mc{T}$ is open and closed in $[-T, T]$, and hence $\mc{T} = [-T, T]$.
    
    The proof of local Lipschitz continuity of the data-to-solution map is analogous to the proof of the contraction mapping property, requiring again $T$ to be sufficiently small.
\end{proof}
We now prove a blow-up alternative. This result is more naturally formulated using intervals $[0, T]$ and $[-T, 0]$ instead of both-sided intervals $[-T, T]$. 
Accordingly, we define analogous spaces $Y_{[0,T]}$ and $Y_{[-T,0]}$, mild solutions on $[0, T]$ and $[-T, 0]$ and derive the same well-posedness result. 
We only consider the case $[0, T]$ below.
\begin{lemma}[Blow-up alternative] \label{lem:11}
    Assume the setting of Theorem \ref{thmlwp}.
    Let $u \in Y_{[0,T]}$ and $\ti{u} \in Y_{[0,\ti{T}]}$ be mild solutions with initial data $u_0 \in H^s \cap H^{s',1}$ and $\ti{u}_0 = u(T)$.
    They can be concatenated to form a mild solution in $Y_{[0,T + \ti{T}]}$. 
    As a result, exactly one of the following statements holds true:
    \begin{enumerate}[(i)]
        \item There exists a maximal time of existence $T^\ast = T^\ast(u_0) > 0$ and a continuous function $u: [0, T^\ast) \rightarrow H^s \cap H^{s',1}$, which is a solution in the sense that $u\big\vert_{t \in [0, T]}$ is the unique mild solution in $Y_T$ with initial data $u_0$ for every $0 < T < T^\ast$.
        Furthermore,
        \begin{align*}
            \sup_{t \in [0, T^\ast)} \|u(t)\|_{H^s \cap H^{s',1}} = \infty \,.
        \end{align*} 
        \item There exists a global mild solution, i.e. a continuous function $u: [0, \infty) \rightarrow H^s \cap H^{s',1}$ such that $u(0) = u_0$ and $u\big\vert_{t \in [0, T]}$ is the unique mild solution in $Y_T$ with initial data $u_0$ for any $T > 0$.
    \end{enumerate}
\end{lemma}
\begin{proof}
    In the proof of uniqueness above, we have demonstrated how mild solutions can be concatenated. 
    Furthermore, mild solutions can be restricted to smaller time intervals. 
    Therefore (ii) can not be true if (i) is true, and it remains to show that if (ii) fails then (i) holds true. 
    Let $u_0 \in H^s \cap H^{s',1}$ and assume that (ii) is false. We consider the set of times
    \begin{align*}
        \mc{T} = \left\{ T > 0: \text{ there exists a mild solution in } Y_{[0, T]} \text{ with initial data } u_0 \right\} \,.
    \end{align*}
    Since mild solutions can be restricted to smaller time intervals, we have $\mc{T} = [0, T^\ast)$ or $\mc{T} = [0, T^\ast]$ for some $T^\ast > 0$, or $\mc{T} = \R$. 
    We can exclude the latter case, as together with the uniqueness result it contradicts (ii). 
    We know that there exists some unique $u \in C_{t \in [0, T^\ast)} H^s_x$ such that $u\big\vert_{t \in [0, T]}$ is the unique mild solution in $Y_T$ for every $T \in [0, T^\ast)$. Suppose now that
    \begin{align*}
        M = \|u\|_{L^\infty_{t \in [0, T^\ast]} (H^s \cap H^{s',1})_x} < \infty \,.
    \end{align*}
    Then any initial data $v_0 = u(T)$ for which $T \in [0, T^\ast)$ fulfills $\|v_0\|_{H^s \cap H^{s',1}} \leq M$, a property which only fails on a zero set, can be continued to a mild solution with a minimal duration $T' = T'(M) > 0$, bounded uniformly from below. 
    By concatenation of solutions, there exists some $T > T^\ast$ for which a mild solution in $Y_T$ with initial data $u_0$ can be found. 
    This contradicts the definition of $\mc{T}$, and hence it must be the case that $M = \infty$.
\end{proof}

\begin{lemma} \label{lem:284}
    Let $s, s' \in \N$ with $s, s' \geq n$ and $s' \leq s - n + 1$. 
    Let $u \in Y_T$ be a mild solution and assume 
    \begin{align} \label{eqn:ass-coeff-2}
        C_{f,g} &= \|g\|_{H^{s-n+2} \cap H^{s'-n+2,1}} + \sum_{k=1}^K \sum_{l \in \N^k, |l| \leq n - 2} \|f^{k,l}\|_{W^{s-n+2,\infty}} < \infty \,.
    \end{align}
    We have $\Phi(D_x) u \in C_{t \in [-T, T]} (H^{s-n} \cap H^{s'-n,1})_x$ and $\mc{N}[u] \in C_{t \in [-T, T]} (H^{s-n+2} \cap H^{s'-n+2,1})_x$. 
    Furthermore, $u \in C^1_{t \in [-T, T]} (H^{s-n} \cap H^{s'-n,1})_x$ with derivative $\del_t u = i \Phi(D_x) u + \mc{N}[u]$.
    In particular, we have for all $h \leq s - n$ the $L^2_x$-Bochner integral identity
    \begin{align*}
        \del_x^h u(t) &= \del_x^h u(0) + \int_0^t \del_x^h(i \Phi(D_x) u(t') + \mc{N}[u(t')]) \dd t' \,.
    \end{align*}
    Note that this implies $\del_x^h (i \Phi(D_x) u(t, x) + \mc{N}[u](t, x)) \in L^1_{t \in [-T, T]}$ for almost all $(t, x) \in \R^2$, and hence the integral identity above holds pointwise for almost all $x \in \R$.
\end{lemma}
\begin{proof}
    We know that $u \in C_{t \in [-T, T]} (H^s \cap H^{s',1})_x$, which implies $\Phi(D_x) u \in C_{t \in [-T, T]} (H^{s-n} \cap H^{s'-n,1})_x$ due to \eqref{eqn:500}. 
    We want to show that $\mc{N}[u] \in C_{t \in [-T, T]} (H^{s-n+2} \cap H^{s'-n+2,1})_x$. Note that we have assumed $g \in H^{s-n+2} \cap H^{s'-n+2}$. 
    We focus now on the $H^{s-n+2}$-norm of the remaining terms, where it suffices to assume $1 \leq k \leq K$, $|l| \leq s - n + 2 + n - 2 = s$, $h \leq s - n + 2$ and apply Hölder's inequality:
    \begin{align*}
        & \left\| \del_x^h f^{k,l} (\mc{N}^{k,l}[u](t) - \mc{N}^{k,l}[u](t')) \right\|_{L^2} 
        \\ &\leq \|\del_x^h f^{k,l}\|_{L^\infty} \sum_{j=1}^k \left\| \del_x^{l_j} u(t) - \del_x^{l_j} u(t') \right\|_{L^{q_j}} \underset{i \neq j}{\prod_{i=1}^k} \left( \left\| \del_x^{l_i} u_{b_i}(t) \right\|_{L^{q_i}} + \left\| \del_x^{l_i} u_{b_i}(t') \right\|_{L^{q_i}} \right)
        \\ &\lesssim C_{f,g} \|u(t) - u(t')\|_{H^s} \|u\|_{L^\infty_{t \in [-T, T]} H^s_x}^{k-1} \,.
    \end{align*}
    Here we choose $q_j = 2$ for the $j$ which maximizes $|l_j|$, and $q_j = \infty$ for all other $j$, where it must be the case that $|l_j| \leq s - 1$ and hence we can use the Sobolev embedding $H^1 \xhookrightarrow{\quad} L^\infty$.
    We estimate the $H^{s'-n+2,1}$-norm the same way, grouping the weight $x$ with the $j$-th factor.
    Next, we show that $u \in C_{t \in [-T, T]} (H^{s-n} \cap H^{s'-n,1})_x$ is differentiable at $t=0$. Differentiability at any $t_0 \in [-T, T]$ can then be obtained by considering $u(t + t_0)$ as a mild solution in $Y_{T - t_0}$ with initial data $u(t_0)$. 
    Applying \eqref{eqn:mild} yields
    \begin{align*}
        \frac{u(t) - u_0}{t} - i \Phi(D_x) u_0 - \mc{N}[u_0]
        &= \left( \frac{e^{i t \Phi(D_x)} - 1}{t} - i \Phi(D_x) \right) u_0
        \\ &+ \frac{1}{t} \int_0^t (e^{i (t - t') \Phi(D_x)} - 1) \mc{N}[u(t')] + \mc{N}[u(t')] - \mc{N}[u_0] \dd t'
        \\ &= (I) + (II) \,.
    \end{align*}
    For the $H^{s-n}$-norm, we use \eqref{eqn:500} to see that for all $h \leq s - n$ we have $\|\del_x^h (I)\|_{L^2_x} \xrightarrow{t \rightarrow 0} 0$, and that
    \begin{align*}
        \left\| \del_x^h (II) \right\|_{L^2_x} &= \left\| \int_0^1 \xi^h \left( (e^{i t (1 - t') \Phi(\xi)} - 1) \ha{\mc{N}[u]}(t t') + \ha{\mc{N}[u]}(t t') - \ha{\mc{N}[u_0]} \right) \dd t' \right\|_{L^2} \xrightarrow{t \rightarrow 0} 0
    \end{align*}
    by dominated convergence. For the $H^{s'-n}$-norm the same methods works, except some additional terms with $\Phi'(D_x)$ appear that can be treated analogously to the unweighted case, since $s' + n - 1 \leq s$.
\end{proof}
\begin{lemma}[Control of $H^s \cap H^{s',1}$-norm] \label{lem:12}
    Let $s, s' \in \N$ with $s \geq 2 n - 1$ and $\frac{s-1}{2} \leq s' \leq s - n$. Assume
    \begin{align} \label{eqn:ass-coeff-3}
        C_{f,g} &= \|g\|_{H^{s-n+2} \cap H^{s',1}} + \sum_{l=0}^{n-2} \|f^{1,l}\|_{H^{s',1}} + \sum_{k=1}^K \sum_{l \in \N^k, |l| \leq n - 2} \|f^{k,l}\|_{W^{s-n+2,\infty}} < \infty \,.
    \end{align}
    There exists $C > 0$ such that for any mild solution $u \in Y_T$ with initial data $u_0 \in H^s \cap H^{s',1}$, we have
    \begin{align*}
        \|u\|_{L^\infty_{t \in [-T, T]} H^{s',1}_x} &\leq \exp\left( (1 + T) C \left( 1 + C_{f,g} + \|u_0\|_{H^{s',1}} + \|u\|_{L^\infty_{t \in [-T, T]} L^\infty_x} + \|u_x\|_{L^\infty_{t \in [-T, T]} H^{s-1}_x} \right)^{K + 1} \right) \,.
    \end{align*}
\end{lemma}
\begin{proof}
    Consider formally the integral form of a differential equation
    \begin{align} \label{eqn:360}
        v(t) &= v(0) + \int_0^t F(v(t')) \dd t' \,.
    \end{align}
    Applying this to $v(t)$ and $\conj{v}(t)$ yields
    \begin{align*}
        v(t) \conj{v}(t) &= v(0) \conj{v(0)} + \int_0^t \conj{v(t')} F(v(t')) + v(t') \conj{F(v(t'))} \dd t' \,.
    \end{align*}
    From Lemma \ref{lem:284}, we know that the integral formulation of our equation
    \begin{align*}
        \del_x^h u(t) = \del_x^h u(0) + \int_0^t \del_x^h (i \Phi(D_x) u(t') + \mc{N}[u(t')]) \dd t'
    \end{align*}
    holds for all $h \leq s - n$. 
    We shall estimate $\|x \del_x^h u\|_{L^\infty_{t \in [-T, T]} L^2_x}$ for $h \leq s' \leq s - n$ with Grönwall's inequality.
    The remaining terms of the form $\|\del_x^h u\|_{L^\infty_{t \in [-T, T]} L^2_x}$ can be treated analogously.
    By the above reasoning and Fubini's theorem, we have
    \begin{align*}
        \int_{\R} x^2 \del_x^h u \conj{\del_x^h u} \dd x 
        &= \int_{\R} x^2 \del_x^h u_0 \conj{\del_x^h u_0} \dd x
        \\ &+ \int_0^t \int_{\R} x^2 \left( i \Phi(D_x) \del_x^h u(t') \conj{\del_x^h u(t')} + \del_x^h u(t') \conj{i \Phi(D_x) \del_x^h u(t')} \right) \dd x \dd t'
        \\ &+ \int_0^t \int_{\R} x^2 \left( \del_x^h \mc{N}[u(t')] \conj{\del_x^h u(t')} + \del_x^h u(t') \conj{\del_x^h \mc{N}[u(t')]} \right) \dd x \dd t'
        \\ &= (I) + (II) + (III) \,.
    \end{align*}
    Firstly, note the trivial estimate $|(I)| \leq \|u_0\|_{H^{s',1}}$.
    Next, we estimate the nonlinear part $(III)$. 
    We have
    \begin{align*}
        \left| \int_0^t \int_{\R} x^2 \del_x^h g \conj{\del_x^h u(t')} \dd x \dd t' \right| 
        &\leq \int_0^t \|g\|_{H^{s',1}} \|u(t')\|_{H^{s',1}} \dd t' 
        \leq \int_0^t \|u(t')\|_{H^{s',1}}^2 + \|g\|_{H^{s',1}}^2 \dd t'  \,.
    \end{align*}
    For the other term, it suffices to assume $1 \leq k \leq K$, $l_1 \geq l_2 \geq \dots \geq l_k$ and $|l| \leq s - n + n - 2 = s - 2$, $h \leq s'$ and estimate $x^2 \del_x^h f^{k,l} \mc{N}^{k,l}[u] \conj{\del_x^h u}$.
    If $k \geq 2$ then $l_k \leq \frac{s-1}{2} \leq s'$, and we have
    \begin{align*}
        & \left| \int_0^t \int_{\R} x^2 \del_x^h f^{k,l} \mc{N}^{k,l}[u](t') \conj{\del_x^h u(t')} \dd x \dd t' \right| 
        \\ &\leq \int_0^t \|\del_x^h f^{k,l} x \del_x^{l_k} u(t')\|_{L^2} \|x \del_x^h u(t')\|_{L^2} \prod_{j=1}^{k-1} \|\del_x^{l_j} u(t')\|_{L^\infty} \dd t'
        \\ &\leq \int_0^t \|u(t')\|_{H^{s',1}}^2 (1 + \|f^{k,l}\|_{W^{s',\infty}} + \|u(t')\|_{L^\infty} + \|u_x(t')\|_{H^{s-1}})^{K-1} \dd t' \,.
    \end{align*}
    If $k = 1$ then we proceed analogously, but place $\del_x^{l_1} u$ into the $L^\infty$-norm and group the weight $x$ with $\del_x^h f^{k,l}$.
    In addition, we need to apply Young's inequality for products.
    It remains to estimate the linear part $(II)$. With the Plancherel theorem, we can write
    \begin{align*}
        - \frac12 (II) &= \Imag\left[ \int_0^t \int_{\R} x \Phi(D_x) D_x^h u(t') \conj{x D_x^h u(t')} \dd x \dd t' \right]
        = \Imag\left[ \int_0^t \int_{\R} \del_\xi(\Phi(\xi) \xi^h \ha{u}(t')) \conj{\del_\xi(\xi^h \ha{u}(t'))} \dd \xi \dd t' \right]
        \\ &= \Imag\left[ \int_0^t \int_{\R} 
            \del_\xi(\Phi(\xi) \xi^h) \ha{u}(t') \conj{\xi^h \del_\xi \ha{u}(t')}
            + \xi^h \del_\xi \ha{u}(t') \conj{\Phi(\xi) \del_\xi(\xi^h) \ha{u}(t')}
        \dd \xi \dd t' \right]
        \\ &+ \Imag\left[ \int_0^t \int_{\R} 
            \del_\xi(\Phi(\xi) \xi^h) \del_\xi(\xi^h) |\ha{u}(t')|^2
            + \Phi(\xi) \xi^{2h} |\del_\xi \ha{u}(t')|^2
        \dd \xi \dd t' \right] \,.
    \end{align*}
    The second integral is zero because it is the imaginary part of a real expression.
    We estimate the first integral using $h \leq s'$ and \eqref{eqn:500}, and obtain
    \begin{align*}
        |(II)| &\lesssim \int_0^t \|u(t')\|_{H^{s',1}} \|u_x(t')\|_{H^{s-1}} \dd t' \leq \int_0^t \|u(t')\|_{H^{s',1}}^2 + \|u_x(t')\|_{H^{s-1}}^2 \dd t' \,.
    \end{align*}
    Here we have paid attention to ensure that at least one full derivative is applied to $u(t)$ when it appears in the unweighted Sobolev norm.
\end{proof}

\subsection{Global well-posedness of the \texorpdfstring{\GP}{GP} hierarchy}

\begin{proof}[Proof of Theorem \ref{thm:1}] \label{proof:main-thm}
    We use the formulation of \eqref{eqn:GP-n} given in Proposition \ref{prop:2}. 
    For the odd flows, we use the linear combination of the hierarchy described by \eqref{eqn:lin-odd}, i.e. we define a new hierarchy
    \begin{align} \label{eqn:tGP-n} \tag{$\ti{\GP}_n$}
        i \del_{t_n} q &= \frac{\mc{H}^{\ti{\GP}}_n}{\delta \conj{q}}
    \end{align}
    with Hamiltonians
    \begin{align*}
        \mc{H}^{\ti{\GP}}_{2m} &= \mc{H}^{\GP}_{2m}
        & \mc{H}^{\ti{\GP}}_{2m+1} &= \sum_{k=0}^m \binom{-\frac12}{m-k} 4^{m-k} \mc{H}^{\GP}_{2k+1} \,.
    \end{align*}
    The \GP and $\ti{\GP}$ hierarchies are equivalent in the sense of Proposition \ref{prop:hierarchy-equivalence}.
    We are now in the setting of a system of dispersive equations that is admissible for Theorem \ref{thmlwp} Lemma, \ref{lem:284}, and Lemma \ref{lem:12}.
    We first verify this for the nonlinear part. Here \eqref{eqn:nonlinear} implies that the coefficients $g_d$ and $f^{k,l}_{d,b}$ have the form
    \begin{align*}
        g_d &= \mathds{1}_{\{n=2\}} \mc{O}^{1,0}_{q_\ast, |q_\ast|^2 - 1} (|q_\ast|^2 - 1) + \mc{O}^{1,n}_{q_\ast}((q_\ast)_x)
        & f^{k,l}_{d,b} &= \mc{O}^{0,n-2}_{q_\ast, |q_\ast|^2 - 1}((q_\ast)_x, |q_\ast|^2 - 1) \,.
    \end{align*}
    These coefficients need to satisfy \eqref{eqn:ass-coeff-1} for Theorem \ref{thmlwp}, \eqref{eqn:ass-coeff-2} for Lemma \ref{lem:284}, and \eqref{eqn:ass-coeff-3} for Lemma \ref{lem:12}.
    In summary, we need $g_d, f^{1,l}_{d,b} \in H^{s+1-\ceil{m},1}$ and $f^{k,l}_{d,b} \in W^{s+1-\ceil{m},\infty}$ for $k \geq 1$. 
    This is fulfilled if $E^{s+1-\ceil{m}+n,1}(q_\ast) < \infty$.
    Now we verify that the linear part is admissible. Clearly \eqref{eqn:500} holds.
    In the case $n = 2 m - 1$, the linear part is $\Phi(D_x) = - D_x^{2m-1}$, which trivially satisfies (P1)--(P4).
    In the case $n = 2 m$, the linear part is
    \begin{align*}
        \mc{L}^{2m} &= D_x^{2m-2} \begin{pmatrix}
            D_x^2 + 2 & 2 & 0 & 0
            \\ - 2 & - D_x^2 - 2 & 0 & 0
            \\ 0 & 0 & - D_x^2 - 2 & - 2
            \\ 0 & 0 & 2 & D_x^2 + 2
        \end{pmatrix}
        = \mc{V} \mc{D} \mc{W} \,,
    \end{align*}
    where with $\langle D_x \rangle = \sqrt{D_x^2 + 4}$ we can write
    \begin{align} \label{eqn:diagonalization}
        & \mc{D} = D_x^{2m-1} \langle D_x \rangle \begin{pmatrix} 
            1 & 0 & 0 & 0
            \\ 0 & -1 & 0 & 0
            \\ 0 & 0 & 1 & 0
            \\ 0 & 0 & 0 & -1
        \end{pmatrix} 
        \\ & {\begingroup
        \setlength\arraycolsep{1pt}
        \mc{V} = \frac12 \begin{pmatrix}
            \frac{D_x}{\langle D_x \rangle} + 1 & \frac{D_x}{\langle D_x \rangle} - 1 & 0 & 0
            \\ \frac{D_x}{\langle D_x \rangle} - 1 & \frac{D_x}{\langle D_x \rangle} + 1 & 0 & 0
            \\ 0 & 0 & \frac{D_x}{\langle D_x \rangle} - 1 & \frac{D_x}{\langle D_x \rangle} + 1
            \\ 0 & 0 & \frac{D_x}{\langle D_x \rangle} + 1 & \frac{D_x}{\langle D_x \rangle} - 1
        \end{pmatrix}
        \;\; \mc{W} = \frac12 \begin{pmatrix}
            \frac{\langle D_x \rangle}{D_x} + 1 & \frac{\langle D_x \rangle}{D_x} - 1 & 0 & 0
            \\ \frac{\langle D_x \rangle}{D_x} - 1 & \frac{\langle D_x \rangle}{D_x} + 1 & 0 & 0
            \\ 0 & 0 & \frac{\langle D_x \rangle}{D_x} - 1 & \frac{\langle D_x \rangle}{D_x} + 1
            \\ 0 & 0 & \frac{\langle D_x \rangle}{D_x} + 1 & \frac{\langle D_x \rangle}{D_x} - 1
        \end{pmatrix} \,.
        \endgroup}
    \end{align}
    Note that $\mc{W} = \mc{V}^{-1}$. Here (P1)--(P4) are fulfilled with $\mu = 1$.
    Applying Theorem \ref{thmlwp}, we obtain a local mild solution $\bm{p}$ with a time of existence $T > 0$ and initial data $\bm{p}_0$.
    Then $q = q_\ast + p_0$ solves \eqref{eqn:tGP-n} with initial data $q_\ast + p_0$ in the sense of distributions.
    Lemma \ref{lem:284} implies that the solution is strong in $H^{s-n} \cap H^{s'-n,1}$.
    Using Lemmas \ref{lem:11} and \ref{lem:12}, it suffices to show that
    \begin{align*}
        \|p\|_{L^\infty_{t \in [-T, T]} L^\infty_x} + \|p_x\|_{L^\infty_{t \in [-T, T]} H^{s-1}_x} \lesssim_{T,s,q_\ast,p_0} 1 \,.
    \end{align*} 
    H. Koch and X. Liao proved in \cite{KochLiao,KochLiao2022} the global well-posedness of \eqref{eqn:GP} for $s > 0$ in the complete metric space $(X^s, d^s)$, whose definition is
    \begin{align} \label{eqn:def-Xs}
        X^s &= \left\{ q \in H^s_{loc}: E^{s-1}(q) < \infty \right\} \big/ \SSS^1 
        & d^s(q, p) &= \left( \int_{\R} \inf_{\lambda \in \SSS^1} \|\sech(\cdot - y) (\lambda p - q)\|_{H^s}^2 \dd y \right)^{\frac12} \,.
    \end{align}
    By \cite[Lemma 2.1]{KochLiao}, the map $[p \mapsto q_\ast + p]: H^s \rightarrow X^s$ is continuous.
    Consequently, we know that $q \in C_{t \in [-T, T]} X^s$. 
    Furthermore, in \cite{KochLiao,KochLiao2022} certain energy functionals $\mc{E}^s: X^s \rightarrow \R$ (see \cite[Theorem 1.3]{KochLiao}) are defined, which are conserved under every flow \eqref{eqn:GP-n} of the hierarchy. 
    They are conserved because they are constructed from the transmission coefficient, which by Lemma \ref{lem:trans-comm} Poisson commutes with the Hamiltonians $\mc{H}^\GP_n$.
    We therefore obtain, as in \cite[Theorem 1.3]{KochLiao}, the a priori bound
    \begin{align*}
        \sup_{t \in [-T, T]} E^s(q(t)) &\lesssim_{s,q_\ast,p_0} 1 \,.
    \end{align*}
    We conclude the proof with the observation that $p(t) = q(t) - q_\ast$ implies
    \begin{align*}
        \|p(t)\|_{L^\infty_x} + \|p_x(t)\|_{H^{s-1}_x} &\lesssim E^s(q_\ast) + E^s(q(t)) \,.
    \end{align*}
\end{proof}

\appendix

\section{Proof of Lemma \ref{lem:678}}
\label{appendix:678}

\begin{definition}
    Let $X$ be a finite-dimensional complex Banach space. We say that $f: D \subset \R \rightarrow X$ is
    \begin{enumerate}[(i)]
        \item \textbf{smooth on $D$}, if it is infinitely differentiable on $D$.
        \item \textbf{bounded smooth on $D$}, if it is bounded, smooth on $D$, and all derivatives are bounded.
        \item \textbf{polynomially bounded on $D$}, if $p^{-1} f$ is bounded on $D$ for some polynomial $p$.
        \item \textbf{polynomially bounded smooth on $D$}, if it is polynomially bounded, smooth on $D$, and all derivatives are polynomially bounded.
        \item \textbf{Rapidly decreasing at $\pm \infty$}, if $p f$ is bounded on $D \cap \R_\pm$ for every polynomial $p$.
        \item \textbf{Schwartz at $\pm \infty$}, if $p f$ is bounded smooth on $D \cap \R_\pm$ for every polynomial $p$.
    \end{enumerate}
\end{definition} 

\begin{proof}[{Proof of Lemma \ref{lem:678}}]
    Recall for $\Gamma^\pm = \Gamma^\pm_\dg + \Gamma^\pm_\odg$ the system
    \begin{align}
        \tag{\ref{eqn:603}} (\del_x + \del_y) \Gamma^\pm_\dg(x, y) &= \Big(Q(x) + \overarcarrow{Q_\ast(y)}\Big) \Gamma^\pm_\odg(x, y)
        \\ \tag{\ref{eqn:604}} (\del_x - \del_y) \Gamma^\pm_\odg(x, y) &= \Big(Q(x) - \overarcarrow{Q_\ast(y)}\Big) \Gamma^\pm_\dg(x, y)
        \\ \tag{\ref{eqn:605}} \lim_{y \rightarrow \pm \infty} \Gamma^\pm(x, y) &= 0 \qquad \qquad \Gamma^\pm_\odg(x, x) = \frac12 (Q - Q_\ast)_\odg(x) \,.
    \end{align}
    We give the proof for the case $\pm = -$; the plus case is analogous. 
    Recall also that we assume $q - q_\ast$ to be Schwartz, and note that the following ``canceling product'' of multiplication operators
    \begin{align} \label{eqn:606}
        \Big(Q(x) - \overarcarrow{Q_\ast(y)}\Big) \Big(Q(x) + \overarcarrow{Q_\ast(y)}\Big) &= Q(x)^2 - \overarcarrow{Q_\ast(y)^2} = (|q|^2 - 1)(x) - (|q_\ast|^2 - 1)(y)
    \end{align}
    is Schwartz if $x = y$. This is an essential source of integrability in our argument.
    The first step is to write \eqref{eqn:603}--\eqref{eqn:605} as a system of integral equations:
    \begin{align*}
        \Gamma^-_\dg(x, y) &= \int_{-\infty}^x \Big(Q(s) + \overarcarrow{Q_\ast(s + y - x)}\Big) \Gamma^-_\odg(s, s + y - x) \dd s
        \\ \Gamma^-_\odg(x, y) &= \frac12 (Q - Q_\ast)\left( \frac{x+y}{2} \right) + \int_{\frac{x+y}{2}}^x \Big(Q(s) - \overarcarrow{Q_\ast(x+y-s)}\Big) \Gamma^-_\dg(s, x + y - s) \dd s \,. 
    \end{align*}
    We now fix $K \in \N$ and consider the variables
    \begin{align*}
        \Gamma^-_K(x, y) &= (\del_y^k \Gamma^-(x, y))_{1 \leq k \leq K} \in \big(\C^{2 \times 2}\big)^K
        \\ (\del_x + \del_y) \Gamma^-_K(x, y) &= ((\del_x + \del_y) \del_y^k \Gamma^-(x, y))_{1 \leq k \leq K} \in \big(\C^{2 \times 2}\big)^K
    \end{align*}
    for which the extended system of integral equations reads
    \begin{align}
        \label{eqn:410} \del_y^k \Gamma^-_\dg(x, y) &= \int_{-\infty}^x \sum_{l=0}^k \binom{k}{l} \del_y^l \Big(Q(s) + \overarcarrow{Q_\ast(s + y - x)}\Big) (\del_y^{k-l} \Gamma^-_\odg)(s, s + y - x) \dd s
        \\ \label{eqn:411} \del_y^k \Gamma^-_\odg(x, y) &= \frac12 \del_y^k (Q - Q_\ast)\left( \frac{x+y}{2} \right)
        \\ \nonumber &\hspace{-1.6cm}+ \int_{\frac{x+y}{2}}^x \sum_{l=0}^k \binom{k}{l} \del_y^l \Big(Q(s) - \overarcarrow{Q_\ast(x + y - s)}\Big) (\del_y^{k-l} \Gamma^-_\dg)(s, x + y - s) \dd s
        \\ \nonumber &\hspace{-1.6cm}- \sum_{j=0}^{k-1} \sum_{l=0}^{k-1-j} \binom{k-1-j}{l} \frac12 \del_y^j \left( \del_y^l \Big(Q(s) - \overarcarrow{Q_\ast(x+y-s)}\Big) \left(\del_y^{k-1-j-l} \Gamma^-_\dg\right)(s, x + y - s) \Bigg\vert_{s=\frac{x+y}{2}} \right)
    \end{align}
    and
    \begin{align}
        \label{eqn:412} (\del_x + \del_y) \del_y^k \Gamma^-_\dg(x, y) &= \sum_{j=0}^k \binom{k}{l} \del_y^l \Big(Q(x) + \overarcarrow{Q_\ast(y)}\Big) (\del_y^{k-l} \Gamma^-_\odg)(x, y)
        \\ \label{eqn:413} (\del_x + \del_y) \del_y^k \Gamma^-_\odg(x, y) &= \del_y^{k+1} (Q - Q_\ast)\left( \frac{x+y}{2} \right)
        \\ \nonumber &\hspace{-1.6cm}+ 2 \int_{\frac{x+y}{2}}^x \sum_{l=0}^{k+1} \binom{k+1}{l} \del_y^l \Big(Q(s) - \overarcarrow{Q_\ast(x + y - s)}\Big) (\del_y^{k+1-l} \Gamma^-_\dg)(s, x + y - s) \dd s
        \\ \nonumber &\hspace{-1.6cm}- \sum_{j=1}^k \sum_{l=0}^{k-j} \binom{k-j}{l} \del_y^j \left( \del_y^l \Big(Q(s) - \overarcarrow{Q_\ast(x + y - s)}\Big) (\del_y^{k-j-l} \Gamma^-_\dg)(s, x + y - s) \Bigg\vert_{s=\frac{x+y}{2}} \right)
        \\ \nonumber &\hspace{-1.6cm}+ \sum_{l=0}^k \binom{k}{l} \left[ \del_y^l \Big(Q(s) - \overarcarrow{Q_\ast(x+y-s)}\Big) (\del_y^{k-l} \Gamma^-_\dg)(s, s)\right]_{s=\frac{x+y}{2}}^{s=x} \,.
    \end{align}
    Note that for $k \geq 1$ \eqref{eqn:603}--\eqref{eqn:604} imply
    \begin{align*}
        (\del_x + \del_y)^k \Gamma^-_\dg(x, y) &= \sum_{j=0}^{k-1} \binom{k-1}{j} (\del_x + \del_y)^{k-1-j} \Big(Q(x) + \overarcarrow{Q_\ast(y)}\Big) (\del_x + \del_y)^j \Gamma^-_\odg(x, y)
        \\(\del_x - \del_y)^k \Gamma^-_\odg(x, y) &= \sum_{j=0}^{k-1} \binom{k-1}{j} (\del_x - \del_y)^{k-1-j} \Big(Q(x) - \overarcarrow{Q_\ast(y)}\Big) (\del_x - \del_y)^j \Gamma^-_\dg(x, y) \,.
    \end{align*}
    As a result, $\del_x^k \Gamma^-$ can be written as a linear (in $\Gamma^-$) combination of $(\del_x^j \Gamma^-)_{0 \leq j \leq k-1}$ and $(\del_y^l \Gamma^-)_{0 \leq l \leq k}$.
    It therefore suffices to construct $\del_y^k \Gamma^-$ and prove the estimates \eqref{eqn:620}--\eqref{eqn:621} for $m = 0$. Smoothness of $\Gamma^-$ then follows by an inductive construction.
    Furthermore, we only need to prove \eqref{eqn:621}, as the observation that our construction yields $\lim_{s \rightarrow - \infty} \Gamma^-(x + s, y + s) = 0$ together with \eqref{eqn:621} directly implies \eqref{eqn:620}.
    
    We begin by analyzing the behavior of \eqref{eqn:410}--\eqref{eqn:413} on the diagonal $y = x$.
    \begin{claim} \label{claim:1}
        The map $[x \mapsto (\del_y^k \Gamma^-)(x, x)]$ is bounded smooth and Schwartz at $- \infty$. 
        Furthermore, we have $[x \mapsto ((\del_x + \del_y) \del_y^k \Gamma^-)(x, x)] \in \mc{S}(\R; \C^{2 \times 2})$ and $\Big[x \mapsto \Big(Q(x) - \overarcarrow{Q_\ast(x)}\Big)(\del_y^k \Gamma^-_\odg)(x, x)\Big] \in \mc{S}(\R; \C^{2 \times 2})$.
    \end{claim}
    \begin{claimproof}
        Specializing \eqref{eqn:410}--\eqref{eqn:413} to $y = x$ yields
        \begin{align}
            (\del_y^k \Gamma^-_\dg)(x, x) &= \int_{-\infty}^x \sum_{l=0}^k \binom{k}{l} \del_y^l \Big(Q(s) + \overarcarrow{Q_\ast(s + y - x)}\Big) \Bigg\vert_{y=x} (\del_y^{k-l} \Gamma^-_\odg)(s, s) \dd s
            \label{eqn:400} \\ (\del_y^k \Gamma^-_\odg)(x, x) &= \frac{1}{2^{k+1}} \del_x^k (Q - Q_\ast)(x) - \sum_{j=0}^{k-1} \sum_{l=0}^{k-1-j} \binom{k-1-j}{l} \frac12 
            \label{eqn:401} \\ &\del_y^j \left( \del_y^l \Big(Q(s) - \overarcarrow{Q_\ast(x + y - s)}\Big) \left(\del_y^{k-1-j-l} \Gamma^-_\dg\right)(s, x + y - s) \Bigg\vert_{s=\frac{x+y}{2}} \right) \Bigg\vert_{y=x}
            \nonumber
        \end{align}
        and
        \begin{align}
            ((\del_x + \del_y) \del_y^k \Gamma^-_\dg)(x, x) &= \sum_{j=0}^k \binom{k}{l} \del_y^l \Big(Q(x) + \overarcarrow{Q_\ast(y)}\Big) \Bigg\vert_{y=x} (\del_y^{k-l} \Gamma^-_\odg)(x, x)
            \label{eqn:402} \\ ((\del_x + \del_y) \del_y^k \Gamma^-_\odg)(x, x) &= \frac{1}{2^{k+1}} \del_x^{k+1} (Q - Q_\ast)(x) - \sum_{j=1}^k \sum_{l=0}^{k-j} \binom{k-j}{l} 
            \label{eqn:403} \\ &\del_y^j \left( \del_y^l \Big(Q(s) - \overarcarrow{Q_\ast(x + y - s)}\Big) (\del_y^{k-j-l} \Gamma^-_\dg)(s, x + y - s) \Bigg\vert_{s=\frac{x+y}{2}} \right) \Bigg\vert_{y=x} \,.
            \nonumber
        \end{align}
        Recall here that $Q$ and $Q_\ast$ are bounded smooth and $Q - Q_\ast$ is Schwartz. In particular, $\Gamma^-_\odg(x, x)$ is Schwartz.
        Then $\Gamma^-_\dg(x, x)$ is Schwartz at $- \infty$ and bounded smooth on $\R$. 
        By induction over $k \in \N$, using \eqref{eqn:400}--\eqref{eqn:401}, it follows that $(\del_y^k \Gamma^-)(x, x)$ is Schwartz at $- \infty$ and polynomially bounded smooth on $\R$.
        Even though we only know $(\del_y^k \Gamma^-_\odg)(x, x)$ to be Schwartz at $- \infty$, we can show that $((\del_x + \del_y) \del_y^k \Gamma^-_\dg)(x, x)$ is Schwartz.
        To see this, consider that if $l > 0$ in \eqref{eqn:402}, a Schwartz factor $\del_x^l Q_\ast(x)$ is present, so it suffices to consider $l = 0$ and show that $\Big(Q(x) + \overarcarrow{Q_\ast(x)}\Big) (\del_y^k \Gamma^-_\odg)(x, x)$ is Schwartz.
        Substituting with \eqref{eqn:401}, we again obtain a Schwartz factor if any derivative from $\del_y^j$ or $\del_y^l$ falls onto $\Big(Q(s) - \overarcarrow{Q_\ast(x + s - y)}\Big)$. 
        In the final case, where this does not happen, we have as a factor the canceling product \eqref{eqn:606}, which is Schwartz.
        Therefore $((\del_x + \del_y) \del_y^k \Gamma^-_\dg)(x, x)$ is Schwartz.
        Lastly, we show that $((\del_x + \del_y) \del_y^k \Gamma^-_\odg)(x, x)$ is also Schwartz.
        In \eqref{eqn:403} we see again that the only difficult term is the one where no derivative from $\del_y^j$ or $\del_y^l$ falls onto $\Big(Q(s) - \overarcarrow{Q_\ast(x + s - y)}\Big)$.
        Since $j \geq 1$, we only need to consider terms where a derivative from $\del_y^j$ has fallen onto $(\del_y^{k-j-l} \Gamma^-_\dg)(s, x+y-s)\vert_{s=\frac{x+y}{2}}$, i.e. there is a factor $((\del_x + \del_y) \del_y^{k-j-l} \Gamma^-_\dg)\big(\frac{x+y}{2}, \frac{x+y}{2}\big)$. 
        Substituting \eqref{eqn:402}, we again obtain a Schwartz factor either from a derivative or the canceling product \eqref{eqn:606}. 
    \end{claimproof}

    The rest of the proof is more elegantly written using the variables
    \begin{align*}
        \Lambda^-_{K,\dg}(x, y) &= \Gamma^-_{K,\dg}(x+y, x-y) & \Lambda^-_{K,\odg}(x, y) &= \Gamma^-_{K,\odg}(x+y, x-y)
    \end{align*}
    and $\Lambda^-_K = \Lambda^-_{K,\dg} + \Lambda^-_{K,\odg}$.
    The integral equations \eqref{eqn:410}--\eqref{eqn:411} can be succinctly written as
    \begin{align}
        \label{eqn:430} \Lambda^-_{K,\dg}(x_0, y_0) &= \int_{-\infty}^{x_0} \Theta_K^+(x_1, y_0) \Lambda^-_{K,\odg}(x_1, y_0) \dd x_1
        \\ \label{eqn:431} \Lambda^-_{K,\odg}(x_0, y_0) &= \int_0^{y_0} \Theta_K^-(x_0, y_1) \Lambda^-_{K,\dg}(x_0, y_1) \dd y_1 + \Gamma^-_{K,\odg}(x_0, x_0)
    \end{align}
    where
    \begin{align*}
        \Theta_K^\pm(x, y) &= \left( \binom{k}{l} 2^{l-k} (\del_x - \del_y)^{k-l} \Big(Q(x + y) \pm \overarcarrow{Q_\ast(x - y)}\Big) \right)_{0 \leq k, l \leq K} \,.
    \end{align*}
    The system of integral equations for $\del_x \Lambda^-_K = \del_x \Lambda^-_{K,\dg} + \del_x \Lambda^-_{K,\odg}$ is
    \begin{align}
        \label{eqn:440} (\del_x \Lambda^-_{K,\dg})(x_0, y_0) &= \int_0^{y_0} \int_{-\infty}^{x_0} (\del_y \Theta_K^+)(x_0, y_1) (\del_x \Lambda^-_{K,\odg})(x_1, y_1)
        \\ &\qquad\qquad + (\Theta_K^+ \Theta_K^-)(x_0, y_1) (\del_x \Lambda^-_{K,\dg})(x_1, y_1) \dd x_1 \dd y_1 + \Theta_K^+(x_0, 0) \Gamma^-_{K,\odg}(x_0, x_0)
        \\ \label{eqn:441} (\del_x \Lambda^-_{K,\odg})(x_0, y_0) &= \int_0^{y_0} \int_{-\infty}^{x_0} (\del_x \Theta_K^-)(x_0, y_1) (\del_x \Lambda^-_{K,\dg})(x_1, y_1)
        \\ &\qquad\qquad + (\Theta_K^- \Theta_K^+)(x_0, y_1) (\del_x \Lambda^-_{K,\odg})(x_1, y_1) \dd x_1 \dd y_1 + ((\del_x + \del_y) \Gamma^-_{K,\odg})(x_0, x_0) \,.
    \end{align}
    To obtain \eqref{eqn:440}, we compute $\del_{x_0} \eqref{eqn:430}$, apply the Fundamental theorem of calculus (FTC) over the interval $(0, y_0)$, substitute a term with $\del_{y_0} \eqref{eqn:431}$ and apply once more the FTC over the interval $(-\infty, x_0)$.
    Similarly, for \eqref{eqn:441} we compute $\del_{x_0} \eqref{eqn:431}$, and apply the FTC over $(- \infty, x_0)$ for one term while substituting with $\del_{x_0} \eqref{eqn:430}$ for the other.
    
    We shall construct $\del_x \Lambda^-_K = \del_x \Lambda^-_{K,\dg} + \del_x \Lambda^-_{K,\odg}$ via the Neumann series Ansatz
    \begin{align*}
        \del_x \Lambda^-_{K,\dg} &= \sum_{n=0}^\infty \del_x \Lambda^-_{K,\dg,n}
        & \del_x \Lambda^-_{K,\odg} &= \sum_{n=0}^\infty \del_x \Lambda^-_{K,\odg,n}
    \end{align*}
    Then formally defining $\Lambda^-_K(x_0, y_0) = \int_{-\infty}^{x_0} (\del_x \Lambda^-_K)(x_1, y_0) \dd x_1$ yields a solution to \eqref{eqn:430}--\eqref{eqn:431}. 
    This becomes rigorous with the estimates for $\del_x \Lambda^-_K$ that we shall obtain from the construction.
    Recall from Claim \ref{claim:1} that $((\del_x + \del_y) \Gamma^-_{K,\odg})(x, x)$ and $\Theta_K^+(x, 0) \Gamma^-_{K,\odg}(x, x)$ are Schwartz, 
    and observe that there exists a positive rapidly decreasing function $g_K \in C_b(\R; \R_+)$ such that
    \begin{align}
        \label{eqn:442} F_K(x, y) = \left( |\del_y \Theta_K^+| + |\del_x \Theta_K^-| + |\Theta_K^+ \Theta_K^-| + |\Theta_K^- \Theta_K^+| \right)(x, y) &\leq g_K(x + y) + g_K(x - y)  \,.
    \end{align}
    This uses the canceling product \eqref{eqn:606} when no derivatives are present. To be precise, there exists a constant $C_K > 0$ such that
    \begin{align*}
        g_K &= C_K \left(1 + \|(q, q_\ast)\|_{C_b^K(\R)} \right) \left( |q|^2 - 1 + |q_\ast|^2 - 1 + \sum_{k=1}^K |\del_x^k q| + |\del_x^k q_\ast| \right)
    \end{align*}
    is a valid choice for $g_K$.
    We now study the iterated application of the integral operators corresponding to the Neumann series Ansatz for \eqref{eqn:440}--\eqref{eqn:441}.
    For $n \geq 1$ we find that there exist sets
    \begin{align*}
        \mc{F}_{K,n} &\subset \left\{ [(x, y) \mapsto \Theta_K^+(x, 0) \Gamma^-_{K,\odg}(x, x)], [(x, y) \mapsto ((\del_x + \del_y) \Gamma^-_{K,\odg})(x, x)] \right\} 
        \\ &\quad \times \left\{ \del_y \Theta_K^+, \del_x \Theta_K^-, \Theta_K^+ \Theta_K^-, \Theta_K^- \Theta_K^+ \right\}^n \,,
    \end{align*}
    consisting of vectors $(f_0, f_1, \dots, f_n)$ of functions, such that
    \begin{align*}
        (\del_x \Lambda^-_{K,n})(x_0, y_0)
        &= \sum_{(f_0, f_1, \dots, f_n) \in \mc{F}_{K,n}} \int \dots \int \mathds{1}_{\{x_0 > \dots > x_n\}} \mathds{1}_{\{y_0 > \dots > y_n > 0\}}
        \\ &\qquad\qquad \prod_{j=1}^n f_j(x_{j-1}, y_j) f_0(x_n, y_n) \dd x_n \dots \dd x_1 \dd y_n \dots \dd y_1 \,.
    \end{align*}
    We use the integrability of $f_0(x_n, y_n)$ to estimate the integral over $\dd x_n$ in $L^1$. 
    Indeed, $f_0(x_n, y_n)$ is independent of $y_n$ and Schwartz in $x_n$, so there exists a bounded function $h_K$ which rapidly decreases at $- \infty$ such that
    \begin{align*}
        \|f_0(x, y)\|_{L^1_x L^\infty_y((-\infty, x_0] \times \R_+)} &\lesssim h_K(x_0) \,.
    \end{align*}
    If $K = 0$ then a valid choice for $h_0$ is
    \begin{align*}
        h_0 &= \left( \|q\|_{L^\infty} + \|q_\ast\|_{L^\infty} \right) (q - q_\ast) \,.
    \end{align*}
    For $K \geq 1$ an explicit example may be constructed by tracing the proof of Claim \ref{claim:1}.
    We now use the $n!$-fold symmetry of the integral (resp. $(n-1)!$-fold, when the $L^\infty_{x_1}$-norm is used instead of the $L^1_{x_1}$-norm) to obtain
    \begin{align*}
        \|(\del_x \Lambda^-_{K,n})(x_1, y_0)\|_{(L^1 \cap L^\infty)_{x_1}((-\infty, x_0])}
        &\lesssim h_K(x_0) |\mc{F}_{K,n}| \|F_K\|_{L^1_x L^1_y((-\infty, x_0] \times [0, y_0])}^{n-1} 
        \\ &\quad \left( \frac{\|F_K\|_{L^1_x L^1_y((-\infty, x_0] \times [0, y_0])}}{n!} + \frac{\|F_K\|_{L^\infty_x L^1_y((-\infty, x_0] \times [0, y_0])}}{(n-1)!} \right) \,.
    \end{align*}
    Define the integral operator $\mc{I} g_K(x) = \int_{-\infty}^x g_K(s) \dd s$.
    The estimate \eqref{eqn:442} implies
    \begin{align*}
        \|F_K\|_{L^\infty_x L^1_y((-\infty, x_0] \times [0, y_0])} &\lesssim \mc{I} g_K(x_0 + y_0)
        \\ \|F_K\|_{L^1_x L^1_y((-\infty, x_0] \times [0, y_0])} &\lesssim \mc{I}^2 g_K(x_0 + y_0) - \mc{I}^2 g_K(x_0 - y_0) \,.
    \end{align*}
    Including the trivial case $n = 0$, we have shown for all $n \geq 0$ that
    \begin{align*}
        & \|(\del_x \Lambda^-_{K,n})(x, y_0)\|_{(L^1 \cap L^\infty)_x((-\infty, x_0])} 
        \\ &\lesssim C^n h_K(x_0) \left( \frac{(\mc{I}^2 g_K(x_0 + y_0))^n}{n!} + \mathds{1}_{\{n \geq 1\}} \mc{I} g_K(x_0 + y_0) \frac{(\mc{I}^2 g_K(x_0 + y_0) - \mc{I}^2 g_K(x_0 - y_0))^{n-1}}{(n-1)!} \right) \,.
    \end{align*}
    Summing over $n \geq 0$, we arrive at
    \begin{align*}
        \|(\del_x \Lambda^-_K)(x, y_0)\|_{(L^1 \cap L^\infty)_x((-\infty, x_0])} &\lesssim h_K(x_0) \left( 1 + \mc{I} g_K(x_0 + y_0) \right) e^{C (\mc{I}^2 g_K(x_0 + y_0) - \mc{I}^2 g_K(x_0 - y_0))} \,,
    \end{align*}
    which implies \eqref{eqn:621}.
\end{proof}

\section{Linear Estimates}

\label{appendix:linear-estimates}
For $f = f(t, x)$ we write $\ha{f}^{(t)}(\tau, x)$ and $\ha{f}^{(x)}(t, \xi)$ for the Fourier transforms in only one variable.
We define
\begin{align*}
    \mc{D}_\otimes(\R^2) = \mc{D}(\R) \otimes \mc{D}(\R) = \left\{\sum_{n=1}^N f_n \otimes g_n: N \in \N, f, g \in \mc{D}(\R)^N\right\} \,.
\end{align*}
\subsection{The local smoothing estimate}
We consider two symbols $\phi \in C^1(\R; \R)$ and $a \in C^{0,1}(\R; \R)$ which have the following properties:
\begin{enumerate}[(H1)]
    \item $\phi$ has finitely many critical points $\xi_1 < \dots < \xi_N$. We define in addition $\xi_0 = - \infty$ and $\xi_{N+1} = \infty$. 
    Set $\eta_j = \phi(\xi_j)$. We consider the bijective restrictions of $\phi$ to the intervals $(\xi_j, \xi_{j+1})$, on which $\phi$ is strictly monotonic, and denote them by
    \begin{align*}
        \phi_j: (\xi_j, \xi_{j+1}) &\longrightarrow (\ha{\eta}_j, \ch{\eta}_j) = \phi((\xi_j, \xi_{j+1}))
        \\ \psi_j = \phi_j^{-1}: (\ha{\eta}_j, \ch{\eta}_j) &\longrightarrow (\xi_j, \xi_{j+1}) 
    \end{align*}
    for $j \in \{0, \dots, N\}$. Here $\ha{\eta}_j = \eta_j \land \eta_{j+1}$ and $\ch{\eta}_j = \eta_j \lor \eta_{j+1}$. We write $\sigma_j = \sign(\phi_j')$ for the sign.
    \item For every $j \in \{1, \dots, N\}$ there exist $\alpha_j \in (0, 1)$ and $q, Q > 0$ such that for all $\eta \in (0, r)$
    \begin{align}
        \label{eqn:ass-1} q \eta^{\alpha_j} &\leq |\psi_{j-1}(\eta_j - \sigma_{j-1} \eta) - \xi_j|, |\psi_j(\eta_j + \sigma_j \eta) - \xi_j| \leq Q \eta^{\alpha_j} 
        \\ \label{eqn:ass-2} q \eta^{\alpha_j - 1} &\leq |\psi_{j-1}'(\eta_j - \sigma_{j-1} \eta)|, |\psi_j'(\eta_j + \sigma_j \eta)| \leq Q \eta^{\alpha_j - 1} \,.
    \end{align} 
    For simplicity, we write $\alpha = \alpha_j$ for all $j \in \{1, \dots, N\}$.
    \item There exist $\beta \in (0, 1)$ and $p, P > 0$ such that for all $\eta > R$
    \begin{align}
        \label{eqn:ass-3} p \eta^\beta &\leq |\psi_0(- \sigma_0 \eta)|, |\psi_N(\sigma_N \eta)| \leq P \eta^\beta
        \\ \label{eqn:ass-4} p \eta^{\beta - 1} &\leq |\psi_0'(- \sigma_0 \eta)|, |\psi_N'(\sigma_N \eta)| \leq P \eta^{\beta - 1} \,.
    \end{align}
    \item There exist $A, \delta > 0$ such that for all $\eta \in \R$
    \begin{align*}
        |a(\eta)| + \langle \xi \rangle^{1 + \delta} |a'(\eta)| &\leq A \,.
    \end{align*}
    \item There exists $M > 0$ such that for any $\tau \in \R$ we can decompose $\R$ into $M$ intervals on whose interiors $\frac{a(\xi)}{\phi(\xi) - \tau}$ is monotonic.
    In particular, for any $j \in \{0, \dots, N\}$ we can decompose $(\tau \lor \ha{\eta}_j, \ch{\eta}_j)$ and $(\ha{\eta}_j, \ch{\eta}_j \land \tau)$ into $M$ intervals, on whose interiors the function $\xi \longmapsto \frac{a(\psi_j(\eta))}{\eta - \tau}$ is monotonic.
\end{enumerate}
Because of (H3) we have $\ha{\eta}_0 = - \sigma_0 \infty$ and $\ch{\eta}_N = \sigma_N \infty$.
Given any function $f: \R \rightarrow \C$ for which the integrals below exist, we can perform the change of variables $\eta = \phi(\xi)$ by splitting the integral:
\begin{align}
    \int_{\R} f(\xi) \phi'(\xi) \dd \xi 
    &= \sum_{j=0}^N \sigma_j \underset{\ha{\eta}_j < \eta < \ch{\eta}_j}{\int} f(\psi_j(\eta)) \dd \eta \,.
\end{align}
We sometimes assume without loss of generality that $\sigma_0 = -1$.

The goal of this section is to prove the following theorem.
\begin{theorem}[{Generalization of \cite[Theorem 3.4]{KenigPonceVega1993}}] \label{thm:thm3.4}
    Let $\phi \in C^1(\R; \R)$ and $a \in C^{0,1}(\R; \R)$ fulfill (H1)-(H5). Let $\ti{a} \in L^\infty(\R; \R)$. Then
    \begin{align} 
        \label{eqn:smooth-1} \left\| \big(\ti{a} |\phi'|^{\frac12}\big)(D_x) e^{i t \phi(D_x)} u_0 \right\|_{L^\infty_x L^2_t} &\lesssim \|u_0\|_{L^2_x}
        \\ \label{eqn:smooth-2} \left\| (a \phi')(D_x) \int_0^\infty e^{i (t - t') \phi(D_x)} f(t',x) \dd t' \right\|_{L^\infty_x L^2_t} &\lesssim_{\phi,a} \|f\|_{L^1_x L^2_t}
        \\ \label{eqn:smooth-3} \left\| (a \phi')(D_x) \int_0^t e^{i (t - t') \phi(D_x)} f(t',x) \dd t' \right\|_{L^\infty_x L^2_t} &\lesssim_{\phi,a} \|f\|_{L^1_x L^2_t} \,.
    \end{align}
\end{theorem}
\begin{proof}
    We directly prove \eqref{eqn:smooth-1}:
    \begin{align*}
        \int_{\R} \left| \big(\ti{a} |\phi'|^{\frac12}\big)(D_x) e^{i t \phi(D_x)} u_0(x) \right|^2 \dd t 
        &= \frac{1}{2 \pi} \int_{\R} \left| \int_{\R} e^{i (x \xi + t \phi(\xi))} \big(\ti{a} |\phi'|^{\frac12}\big)(\xi) \widehat{u_0}(\xi) \dd \xi \right|^2 \dd t
        \\ &\lesssim_{\|\ti{a}\|_{L^\infty}} \int_{\R} \left| \sum_{j=0}^N \mathds{1}_{(\ha{\eta}_j, \ch{\eta}_j)}(t) e^{i x \psi_j(t)} (\widehat{u_0} \circ \psi_j)(t) |\phi_j' \circ \psi_j|^{-\frac12}(t) \right|^2 \dd t
        \\ &\lesssim_N \sum_{j=0}^N \int_{\xi_j}^{\xi_{j+1}} |\widehat{u_0}(\xi)|^2 |\phi_j'(\xi)|^{-1} |\phi_j'(\xi)| \dd \xi
        \\ &\lesssim_N \left\| \widehat{u_0} \right\|_{L^2_x}^2 \,.
    \end{align*}
    From \eqref{eqn:smooth-1} we know that $\big(\ti{a} |\phi'|^{\frac12}\big)(D_x) e^{i t \phi(D_x)}: L^2_x(\R) \rightarrow L^\infty_x L^2_t(\R \times \R_+)$ is bounded.
    The boundedness of the corresponding adjoint yields the estimate
    \begin{align*}
        \left\| \big(\ti{a} |\phi'|^{\frac12}\big)(D_x) \int_0^\infty e^{i (t - t') \phi(D_x)} f(t',x) \dd t' \right\|_{L^\infty_t L^2_x(\R_+ \times \R)} &\lesssim_{N,A} \|f\|_{L^1_x L^2_t(\R \times \R_+)} \,.
    \end{align*}
    Consider now a test function $g \in \mc{D}_\otimes(\R^2)$. We apply the above with $\ti{a} = |a|^{\frac12}$ and $\ti{a} = \sign(a \phi') |a|^{\frac12}$ to obtain
    \begin{align*}
        & \int_0^\infty \int_{\R} (a \phi')(D_x) \int_0^\infty e^{i (t - t') \phi(D_x)} f(t',x) \dd t' g(t, x) \dd x \dd t
        \\ &= \int_{\R} \big(|a \phi'|^{\frac12}\big)(D_x) \int_0^\infty e^{- i t' \phi(D_x)} f(t',x) \dd t'\, \big(\sign(a \phi') |a \phi'|^{\frac12}\big)(D_x) \int_0^\infty e^{- i t \phi(D_x)} g(t,x) \dd t \dd x
        \\ &\lesssim_{N,A} \|f\|_{L^1_x L^2_t(\R \times \R_+)} \|g\|_{L^1_x L^2_t(\R \times \R_+)} \,,
    \end{align*}
    By duality, this implies \eqref{eqn:smooth-2}.
    For \eqref{eqn:smooth-3} we follow \cite[§3]{KenigPonceVega1993}. 
    This method is also outlined in \cite[Theorem 2.3]{KenigPonceVega-SmallSolutions} and \cite{KenigPonceVega-BenjaminOno}.
    Consider the solution
    \begin{align*}
        u(t, x) &= \int_0^t e^{i (t - t') \phi(D_x)} f(t', x) \dd t' \quad \text{to} \quad D_t u = \phi(D_x) u - i f \quad \text{with initial data} \quad u(0) = 0 \,.
    \end{align*}
    Then
    \begin{align*}
        \ha{u}(\tau, \xi) &= \frac{- i \ha{f}(\tau, \xi)}{\tau - \phi(\xi)} \,.
    \end{align*}
    Define
    \begin{align*}
        v(t, x) &= \PV \frac{1}{2 \pi} \int_{\R} \int_{\R} e^{i t \tau} e^{i x \xi} \frac{- i \ha{f}(\tau, \xi)}{\tau - \phi(\xi)} \dd \xi \dd \tau \,,
    \end{align*}
    which is formally another solution, although it may or may not fulfill $v(0) = 0$. We can characterize the difference between $u$ and $v$ by
    \begin{align*}
        u(t, x) &= v(t, x) - e^{i t \phi(D_x)} v(0, x) \,.
    \end{align*}
    Using Proposition \ref{prop:prop3.1} and $\ha{\sign}(\tau) = - i \sqrt{\frac{2}{\pi}} \PV \left( \frac{1}{\tau} \right)$, we write
    \begin{align*}
        (a \phi')(D_x) v(t, x) &= \frac{1}{\sqrt{2 \pi}} \int_{\R} e^{i x \xi} \PV \frac{1}{\sqrt{2 \pi}} \int_{\R} e^{i t \tau} \frac{- i (a \phi')(\xi) \ha{f}(\tau, \xi)}{\tau - \phi(\xi)} \dd \tau \dd \xi
        \\ &= \frac{1}{\sqrt{2 \pi}} \int_{\R} e^{i x \xi} \frac12 \int_{\R} e^{i \tau (t - t')} \sign(t - t') (a \phi')(\xi) \ha{f}^{(x)}(t', \xi) \dd t' \dd \xi
        \\ &= \frac12 \int_{\R} (a \phi')(D_x) e^{i (t - t') \phi(D_x)} \sign(t - t') f(t', x) \dd t' \,.
    \end{align*}
    We obtain the decomposition
    \begin{align*}
        & (a \phi')(D_x) u(t, x) 
        \\ &= \frac12 \int_{\R} (a \phi')(D_x) e^{i (t - t') \phi(D_x)} (\sign(t - t') - \sign(- t')) f(t', x) \dd t'
        \\ &= (a \phi')(D_x) v(t, x) + \frac12 \int_0^\infty (a \phi')(D_x) e^{i (t - t') \phi(D_x)} f(t', x) \dd t' - \frac12 \int_{-\infty}^0 (a \phi')(D_x) e^{i (t - t') \phi(D_x)} f(t', x) \dd t' \,.
    \end{align*}
    These calculations are identical to \cite[Proposition 3.1, Lemma 3.4]{KenigPonceVega1993} for $\phi(\xi) = \xi^3$.
    We can estimate the integrals at the end with (ii) and a corresponding version of (ii) on $\R_-$, 
    so it remains to show that $\|(a \phi')(D_x) v\|_{L^\infty_x L^2_t} \lesssim_{a, \phi} \|f\|_{L^1_t L^\infty_x}$.
    By Proposition \ref{prop:prop3.2}, we have
    \begin{align*}
        (a \phi')(D_x) v(t, x) &= \frac{1}{2 \pi} \int_{\R} e^{i t \tau} \PV \int_{\R} e^{i x \xi} \frac{- i (a \phi')(\xi) \ha{f}(\tau, \xi)}{\tau - \phi(\xi)} \dd \xi \dd \tau \,,
    \end{align*}
    and by Proposition \ref{prop:prop3.3.1} there exists some kernel $K \in L^\infty(\R^2)$ such that
    \begin{align*}
        (a \phi')(D_x) v(t, x) &= \frac{1}{2 \pi} \int_{\R} e^{i t \tau} \frac{i}{\sqrt{2 \pi}} \int_{\R} K(x - y, \tau) \widehat{f}^{(t)}(y, \tau) \dd y \dd \tau \,.
    \end{align*}
    These propositions correspond to \cite[Proposition 3.2, 3.3]{KenigPonceVega1993}.
    We use Plancherel's theorem and Minkowski's inequality to obtain
    \begin{align*}
        \|(a \phi')(D_x) v\|_{L^\infty_x L^2_t} &= \frac{1}{(2 \pi)^{\frac32}}  \left\| \left( \int_{\R} \left| \int_{\R} K(x - y, \tau) \widehat{f}^{(t)}(y, \tau) \dd y \right|^2 \dd \tau \right)^{\frac12} \right\|_{L^\infty_x}
        \\ &\leq \frac{1}{(2 \pi)^{\frac32}} \|K\|_{L^\infty} \int_{\R} \left( \int_{\R} |\widehat{f}^{(t)}(y, \tau)|^2 \dd \tau \right)^{\frac12} \dd y
        \\ &\lesssim_{\phi,a} \|f\|_{L^1_x L^2_t} \,.
    \end{align*}
\end{proof}
\begin{proposition}[{Generalization of \cite[Proposition 3.1]{KenigPonceVega1993}}] \label{prop:prop3.1}
    Let $f \in \mathcal{D}_\otimes(\R^2)$. Then for all $(t, x) \in \R^2$ we have
    \begin{align*}
        \lim_{\epsilon \searrow 0} \underset{\epsilon < |\phi(\xi) - \tau| < \frac{1}{\epsilon}}{\int \int} e^{i (x \xi + t \tau)} \frac{(a \phi')(\xi)}{\phi(\xi) - \tau} \widehat{f}(\xi, \tau) \dd \xi \dd \tau 
        = \int_{\R} e^{i x \xi} \lim_{\epsilon \searrow 0} \underset{\epsilon < |\phi(\xi) - \tau| < \frac{1}{\epsilon}}{\int} e^{i t \tau} \frac{(a \phi')(\xi)}{\phi(\xi) - \tau} \widehat{f}(\xi, \tau) \dd \tau \dd \xi \,.
    \end{align*}
    Moreover, both terms are finite and both limits exist.
\end{proposition}
\begin{proof}
    The proof is essentially identical to that of \cite[Proposition 3.1]{KenigPonceVega1993}.
    It suffices to assume that $f(t, x) = v(x) w(t)$. Then
    \begin{align*}
        \lim_{\epsilon \searrow 0} \underset{\epsilon < |\phi(\xi) - \tau| < \frac{1}{\epsilon}}{\int} e^{i t \tau} \frac{(a \phi')(\xi)}{\phi(\xi) - \tau} \widehat{f}(\xi, \tau) \dd \tau
        &= (a \phi' \ha{v})(\xi) H[e^{i t \tau} \ha{w}(\tau)](\phi(\xi)) \,,
    \end{align*}
    where $H$ denotes the Hilbert transform. Since $H: \mc{D}(\R) \rightarrow L^\infty(\R)$, the integral
    \begin{align*}
        \int_{\R} e^{i x \xi} (a \phi' \ha{v})(\xi) H[e^{i t \tau} \ha{w}(\tau)](\phi(\xi)) \dd \xi
    \end{align*}
    is absolutely convergent.
    Define
    \begin{align*}
        H^\ast[g](x) &= \sup_{\epsilon > 0} H_\epsilon[g](x) = \sup_{\epsilon > 0} \left| \underset{\epsilon < |x - y| < \frac{1}{\epsilon}}{\int} \frac{g(y)}{x - y} \dd y \right|
    \end{align*}
    and recall that $\lim_{\epsilon \searrow 0} H_\epsilon[g](x) = H[g](x)$ for almost all $x \in \R$, and furthermore that $H^\ast: L^p(\R) \rightarrow L^p(\R)$ is bounded for $1 < p < \infty$.
    Now consider
    \begin{align*}
        \underset{\epsilon < |\phi(\xi) - \tau| < \frac{1}{\epsilon}}{\int \int} e^{i (x \xi + t \tau)} \frac{(a \phi')(\xi)}{\phi(\xi) - \tau} \widehat{f}(\xi, \tau) \dd \xi \dd \tau 
        &= \int_{\R} e^{i x \xi} (a \phi' \ha{v})(\xi) H[e^{i t \tau} \ha{w}(\tau)](\phi(\xi)) \dd \xi
        \\ &- \int_{\R} e^{i x \xi} (a \phi' \ha{v})(\xi) (H - H_\epsilon)[e^{i t \tau} \ha{w}(\tau)](\phi(\xi)) \dd \xi
        \\ &= (I) + (II_\epsilon) \,.
    \end{align*}
    It remains to prove that $\lim_{\epsilon \searrow 0} |(II_\epsilon)| = 0$.
    We have
    \begin{align*}
        (II_\epsilon) &= \sum_{j=0}^N \sigma_j \int_{\ha{\eta}_j}^{\ch{\eta}_j} e^{i x \xi} (a \ha{v})(\psi_j(\eta)) (H - H_\epsilon)[e^{i t \tau} \ha{w}(\tau)](\eta) \dd \eta \,.
    \end{align*}
    Since $a, v \in L^\infty(\R)$ and $|(H_\epsilon - H)[g]| \leq 2 H^\ast[g]$ is bounded on $L^p(\R)$, the claim follows by dominated convergence.
\end{proof}
Define
\begin{align*}
    \mc{K}_\epsilon[g](\tau) &= \underset{\epsilon < |\phi(\xi) - \tau| < \frac{1}{\epsilon}}{\int} \frac{(a \phi')(\xi)}{\phi(\xi) - \tau} g(\xi) \dd \xi
    & \mc{K}[g](\tau) &= \lim_{\epsilon \searrow 0} \mc{K}_\epsilon[g](\tau) 
    & \mc{K}^\ast[g](\tau) &= \sup_{\epsilon > 0} |K[g](\tau)| \,.
\end{align*}
\begin{proposition}[{Generalization of \cite[Proposition 3.2]{KenigPonceVega1993}}] \label{prop:prop3.2}
    Let $f \in \mathcal{D}_\otimes(\R^2)$. Then for all $(t, x) \in \R^2$ we have
    \begin{align*}
        \lim_{\epsilon \searrow 0} \underset{\epsilon < |\phi(\xi) - \tau| < \frac{1}{\epsilon}}{\int \int} e^{i (x \xi + t \tau)} \frac{(a \phi')(\xi)}{\phi(\xi) - \tau} \widehat{f}(\xi, \tau) \dd \xi \dd \tau 
        = \int_{\R} e^{i t \tau} \lim_{\epsilon \searrow 0} \underset{\epsilon < |\phi(\xi) - \tau| < \frac{1}{\epsilon}}{\int} e^{i x \xi} \frac{(a \phi')(\xi)}{\phi(\xi) - \tau} \widehat{f}(\xi, \tau) \dd \xi \dd \tau \,.
    \end{align*}
    Moreover, both integrals are absolutely convergent and both limits exist.
\end{proposition}
\begin{proof}
    The proof is again almost identical to that of \cite[Proposition 3.2]{KenigPonceVega1993}.
    In Propsition \ref{prop:prop3.1} the first limit is shown to exist, using that for $g \in \mc{D}(\R)$ the limit $H(x) = \lim_{\epsilon \searrow 0} H_\epsilon[g](x)$ exists for almost all $x \in \R$ and that $H^\ast[g] \in L^p(\R)$.
    Both the statement and the proof of this proposition are analogous to that of the previous one, except $H, H_\epsilon$ and $H^\ast$ are replaced by $\mc{K}, \mc{K}_\epsilon$ and $\mc{K}^\ast$.
    The key steps in the proof are therefore the existence of the limit $\mc{K}[g]$, and $\mc{K}^\ast[g] \in L^p(\R)$ for $1 < p < \infty$.
    We write
    \begin{align*}
        \mc{K}_\epsilon[g](\tau) &= \sum_{j=0}^N \sigma_j \underset{\epsilon < |\eta - \tau| < \frac{1}{\epsilon}}{\int} 
        \frac{\mathds{1}_{(\ha{\eta}_j, \ch{\eta}_j)}(\eta) (a g)(\psi_j(\eta))}{\eta - \tau} \dd \eta \,.
    \end{align*}
    Since $a \in L^\infty(\R)$, it suffices to prove that $g \circ \psi_j \in L^p((\ha{\eta}_j, \ch{\eta}_j))$ for $1 \leq p \leq \infty$ and use known properties of the Hilbert transform.
    We have
    \begin{align*}
        \sigma_j \int_{\ha{\eta}_j}^{\ch{\eta}_j} |g(\psi_j(\eta))|^p \dd \eta &= \int_{\xi_j}^{\xi_{j+1}} |g(\xi)|^p \phi_j'(\xi) \dd \xi \,,
    \end{align*}
    so our assumptions $\phi \in C^1(\R, \R)$ and $g \in \mc{D}(\R)$ are sufficient.
\end{proof}
Define
\begin{align*}
    K_\epsilon(z, \tau) &= \underset{\epsilon < |\phi(\xi) - \tau| < \frac{1}{\epsilon}}{\int} e^{i z \xi} \frac{(a \phi')(\xi)}{\phi(\xi) - \tau} \dd \xi
    & K(z, \tau) &= \lim_{\epsilon \searrow 0} K_\epsilon(z, \tau) \,.
\end{align*}
\begin{proposition}[{Generalization of \cite[Proposition 3.3]{KenigPonceVega1993}}] \label{prop:prop3.3.2}
    There exists some $\epsilon_0 > 0$ for which the family $\{K_\epsilon\}_{0 < \epsilon < \epsilon_0} \subset L^\infty(\R^2)$ is bounded.
\end{proposition}
Here the proof requires some notable modifications from that of \cite[Proposition 3.3]{KenigPonceVega1993}. 
This stems from the fact that the change of variables $\eta = \phi(\xi)$ is much more involved. 
In fact, after stating and proving the subsequent proposition, we shall devote the rest of this section to the proof of Proposition \ref{prop:prop3.3.2}.
\begin{proposition}[{Generalization of \cite[Proposition 3.3]{KenigPonceVega1993}}] \label{prop:prop3.3.1}
    The limit $K(z, \tau)$ exists for almost all $(z, \tau) \in \R^2$ and $|K(z, \tau)| \leq C(\phi, a)$.
    For any $f \in \mathcal{D}_\otimes(\R^2)$ and $(t, x) \in \R^2$ we have
    \begin{align} \label{eqn:prop3.3}
        \lim_{\epsilon \searrow 0} \frac{1}{\sqrt{2 \pi}} \underset{\epsilon < |\phi(\xi) - \tau| < \frac{1}{\epsilon}}{\int} e^{i x \xi} \frac{(a \phi')(\xi)}{\phi(\xi) - \tau} \widehat{f}(\xi, \tau) \dd \xi 
        = \frac{1}{\sqrt{2 \pi}} \int_{\R} K(x - y, \tau) \widehat{f}^{(t)}(y, \tau) \dd y \,.
    \end{align}
\end{proposition}
\begin{proof}
    We obtain \eqref{eqn:prop3.3} with Plancherel's theorem and dominated convergence. 
    Here the existence of a majorant is given by Proposition \ref{prop:prop3.3.2}, so it remains to show that the limit exists for almost all $(z, \tau) \in \R^2$.
    We decompose
    \begin{align}
        K_\epsilon(z, \tau)  &= \sum_{j=0}^N \sigma_j \underset{\ha{\eta}_j < \eta < \ch{\eta}_j}{\underset{\epsilon < |\eta - \tau| < \frac{1}{\epsilon}}{\int}} 
        e^{i z \psi_j(\eta)} \frac{a(\psi_j(\eta))}{\eta - \tau} \dd \eta \,.
    \end{align}
    For $j \in \{1, \dots, N-1\}$ we are on a finite interval, so by (H4) the limit exists for almost all $\tau \in \R$ as the Hilbert transform of a function in $L^p(\R)$, $1 < p < \infty$.
    For $j \in \{0, N\}$ we are on an infinite interval, where we may assume without loss of generality that $j = 0$ and the interval is $(R, \infty)$ for large $R$.
    If $\tau > R$, i.e. a singularity is present, then we cut out another finite interval and consider the limit again as the Hilbert transform of a function in $L^p(\R)$.
    We may therefore assume $R > \tau + 1$. Transforming back to $\xi = \psi_0(\eta)$, we have to show the existence of
    \begin{align*}
        \lim_{\epsilon \rightarrow 0} \int_{\R} e^{i z \xi} \mathds{1}_{\{|\phi(\xi) - \tau| < \frac{1}{\epsilon}\}} \mathds{1}_{\{\xi > \psi_0(R)\}} \frac{(a \phi')(\xi)}{\phi(\xi) - \tau} \dd \xi \,.
    \end{align*}
    Since $a \in L^\infty(\R)$ it remains to show that $\frac{\phi'}{\phi} \in L^2((\psi_0(R), \infty))$, as then the integrand is a continuous family in $\epsilon$ with values in $L^2(\R)$, and hence by the continuity of the Fourier transform on $L^2(\R)$ the limit exists for almost all $z \in \R$.
    Indeed, we have
    \begin{align*}
        \int_{\psi_0(R)}^\infty \frac{|\phi'(\xi)|^2}{|\phi(\xi)|^2} \dd \xi 
        &= \int_R^\infty \frac{1}{\eta^2 |\psi_0(\eta)|} \dd \eta 
        \lesssim_p \int_1^\infty \frac{1}{\eta^2 \eta^{\beta - 1}} \dd \eta < \infty \,.
    \end{align*}
\end{proof}
To prove Proposition \ref{prop:prop3.3.2}, we need some Van der Corput type lemmas.
\begin{lemma} \label{lem:998}
    Let $f \in C^1([a, b]; \R)$ such that $f'$ is nonzero and monotonic. Then for all $z \in \R$
    \begin{align*}
        \left| \int_a^b e^{i z f(x)} \dd x \right| &\leq \frac{4}{|z|} \sup_{x \in (a, b)} \frac{1}{|f'(x)|} \,.
    \end{align*}
\end{lemma}
\begin{proof}
    We integrate by parts and note that $\frac{1}{f'(x)}$ is monotonic.
\end{proof}
\begin{lemma} \label{lem:999}
    Let $f \in C^1([a, b]; \R)$ such that $f'$ is nonzero and monotonic, 
    and let $g \in C^{0,1}([a, b]; \R)$ have the property that $[a, b]$ can be decomposed into $M$ intervals on whose interiors $g$ is monotonic.
    Then
    \begin{align*}
        \left| \int_a^b e^{i z f(x)} g(x) \dd x \right| &\leq \frac{12 M}{|z|} \sup_{x \in [a, b]} \frac{1}{|f'(x)|} \sup_{x \in [a, b]} |g(x)| \,.
    \end{align*}
\end{lemma}
\begin{proof}
    We integrate by parts and use the previous lemma, as well as the same monotonicity trick to write $\|g'\|_{L^1([a, b])} = |g(b) - g(a)|$.
\end{proof}
\begin{proof}[Proof of Proposition \ref{prop:prop3.3.2}]
    Recall that
    \begin{align*}
        K_\epsilon(z, \tau) &= \underset{\epsilon < |\phi(\xi) - \tau| < \frac{1}{\epsilon}}{\int} e^{i z \xi} \frac{(a \phi')(\xi)}{\phi(\xi) - \tau} \dd \xi
        = \sum_{j=0}^N \sigma_j \underset{\ha{\eta}_j < \eta < \ch{\eta}_j}{\underset{\epsilon < |\eta - \tau| < \frac{1}{\epsilon}}{\int}} e^{i z \psi_j(\eta)} \frac{a(\psi_j(\eta))}{\eta - \tau} \dd \eta \,.
    \end{align*}
    We fix aa small $r > 0$ and consider the strip $(\tau - r, \tau + r)$. 
    For every $0 \leq j \leq N$ exactly one of three cases below holds true. 
    In each case the integral
    \begin{align*}
        \underset{\ha{\eta}_j < \eta < \ch{\eta}_j}{\underset{\epsilon < |\eta - \tau| < \frac{1}{\epsilon}}{\int}} e^{i z \psi_j(\eta)} \frac{a(\psi_j(\eta))}{\eta - \tau} \dd \eta
    \end{align*}
    has to be estimated uniformly in $\tau$, $z$, and $\epsilon$. We fix a large number $R > 0$.
    \begin{enumerate}[(i)]
        \item \textbf{The continuous case} $(\tau - r, \tau + r) \cap (\ha{\eta}_j, \ch{\eta}_j) = \varnothing$.
        \begin{enumerate}[(a)]
            \item \textbf{The continuous bounded case} $j \not\in \{0, N\}$. Here we integrate a continuous function over a bounded domain. This case is trivial, since for any interval $(a, b)$ with $|\eta - \tau|\big\vert_{(a, b)} > r$ we have
            \begin{align*}
                \left| \int_a^b e^{i z \psi_j(\eta)} \frac{\mathds{1}_{\{|\eta - \tau| > \epsilon\}} a(\psi_j(\eta))}{\eta - \tau} \dd \eta \right| &\leq \frac{A}{r} |b - a| \,.
            \end{align*}
            \item \textbf{The continuous oscillatory tail case} $j \in \{0, N\}$. Here we have to combine the oscillation with the $\frac{1}{\eta}$ decay to bound the integral. By using (a), we can furthermore assume without loss of generality that the integral ranges over the intervals $(R, \infty)$ or $(- \infty, - R)$, where $R = R(\phi) > 0$ is arbitrarily large.
        \end{enumerate}
        \item \textbf{The non-critical point case} $(\tau - r, \tau + r) \subset (\ha{\eta}_j, \ch{\eta}_j)$.
        \begin{enumerate}[(a)]
            \setcounter{enumii}{2}
            \item \textbf{The non-critical point singularity bounded case} $j \not\in \{0, N\}$. Here we have to use oscillation as well as cancellation to bound the integral.
            \item \textbf{The non-critical point singularity oscillatory tail case} $j \in \{0, N\}$. Here the singularity is possibly far away in the oscillatory tail. This can be reduced to a combination of (a), (b) and (c), but it requires an improvement of (b) as the oscillation in the integral may become weak far away from the origin.
        \end{enumerate}
        \item \textbf{The critical point case} $(\tau - r, \tau + r) \cap (\ha{\eta}_j, \ch{\eta}_j) \neq \varnothing$ and either $\eta_j \in (\tau - r, \tau + r)$ or $\eta_{j+1} \in (\tau - r, \tau + r)$. 
        This case always comes in pairs, meaning that if $\eta_j \in (\tau - r, \tau + r)$, then $j-1$ is also of the critical point case, 
        and if $\eta_{j+1} \in (\tau - r, \tau + r)$ then $j+1$ is also of the critical point case. 
        We assume without loss of generality that the former is the case and furthermore that $\sigma_{j-1} = 1$. Then $\sigma_j = -1$ if and only if $\xi_j$ is a strict local maximum.
        If $\sigma_j = -1$ we decompose
        \begin{align*}
            & \underset{\ha{\eta}_{j-1} < \eta < \ch{\eta}_{j-1}}{\underset{\epsilon < |\eta - \tau| < \frac{1}{\epsilon}}{\int}} e^{i z \psi_{j-1}(\eta)} \frac{a(\psi_{j-1}(\eta))}{\eta - \tau} \dd \eta 
            - \underset{\ha{\eta}_j < \eta < \ch{\eta}_j}{\underset{\epsilon < |\eta - \tau| < \frac{1}{\epsilon}}{\int}} e^{i z \psi_j(\eta)} \frac{a(\psi_j(\eta))}{\eta - \tau} \dd \eta
            \\ &= \int_{\eta_{j-1}}^{\tau - r} e^{i z \psi_{j-1}(\eta)} \frac{\mathds{1}_{\{|\eta - \tau| > \epsilon\}} a(\psi_{j-1}(\eta))}{\eta - \tau} \dd \eta 
            - \int_{\eta_{j+1}}^{\tau - r} e^{i z \psi_j(\eta)} \frac{\mathds{1}_{\{|\eta - \tau| > \epsilon\}} a(\psi_j(\eta))}{\eta - \tau} \dd \eta
            \\ &+ \int_{\tau - r}^{\eta_j} e^{i z \psi_{j-1}(\eta)} \frac{\mathds{1}_{\{|\eta - \tau| > \epsilon\}} a(\psi_{j-1}(\eta))}{\eta - \tau} \dd \eta 
            - \int_{\tau - r}^{\eta_j} e^{i z \psi_j(\eta)} \frac{\mathds{1}_{\{|\eta - \tau| > \epsilon\}} a(\psi_j(\eta))}{\eta - \tau} \dd \eta \,,
        \end{align*}
        and if $\sigma_j = 1$ we decompose
        \begin{align*}
            & \underset{\ha{\eta}_{j-1} < \eta < \ch{\eta}_{j-1}}{\underset{\epsilon < |\eta - \tau| < \frac{1}{\epsilon}}{\int}} e^{i z \psi_{j-1}(\eta)} \frac{a(\psi_{j-1}(\eta))}{\eta - \tau} \dd \eta 
            + \underset{\ha{\eta}_j < \eta < \ch{\eta}_j}{\underset{\epsilon < |\eta - \tau| < \frac{1}{\epsilon}}{\int}} e^{i z \psi_j(\eta)} \frac{a(\psi_j(\eta))}{\eta - \tau} \dd \eta
            \\ &= \int_{\eta_{j-1}}^{\tau - r} e^{i z \psi_{j-1}(\eta)} \frac{\mathds{1}_{\{|\eta - \tau| > \epsilon\}} a(\psi_{j-1}(\eta))}{\eta - \tau} \dd \eta 
            + \int_{\tau + r}^{\eta_{j+1}} e^{i z \psi_j(\eta)} \frac{\mathds{1}_{\{|\eta - \tau| > \epsilon\}} a(\psi_j(\eta))}{\eta - \tau} \dd \eta
            \\ &+ \int_{\tau - r}^{\eta_j} e^{i z \psi_{j-1}(\eta)} \frac{\mathds{1}_{\{|\eta - \tau| > \epsilon\}} a(\psi_{j-1}(\eta))}{\eta - \tau} \dd \eta 
            + \int_{\eta_j}^{\tau + r} e^{i z \psi_j(\eta)} \frac{\mathds{1}_{\{|\eta - \tau| > \epsilon\}} a(\psi_j(\eta))}{\eta - \tau} \dd \eta \,.
        \end{align*}
        In both cases the integrals in the first line have $|\eta - \tau| \geq r$ and can hence be treated with either (a) or (b), depending on if the integral is over a bounded or unbounded domain. 
        The remaining integrals need further case distinction.
        \begin{enumerate}
            \setcounter{enumii}{4}
            \item \textbf{The singularity near critical point case} $\sigma_j = -1$ and $\tau < \eta_j$, or $\sigma_j = 1$. Here the integral has one or two singularities in close proximity to the critical point.
            \item \textbf{The almost singularity near critical point case} $\sigma_j = -1$ and $\tau > \eta_j$. Here the integral has no singularity, but can be arbitrary close to being singular. This case is easier than the previous one.
            \item \textbf{The singularity at critical point case} $\tau = \eta_j$. This case is strictly more difficult than certain sections of the integrals to estimate in the previous cases, so it serves as the prototypical case.
        \end{enumerate}
    \end{enumerate}

    \begin{figure}[H]
        \caption{Depicted is an example for $\phi$ and some of the aforementioned cases.}
        \centering
        \definecolor{oiBlue}   {HTML}{0062E2}
        \definecolor{oiOrange} {HTML}{E55E00}
        \begin{tikzpicture}
            \pgfplotsset{compat=1.15}
            \pgfmathdeclarefunction{f}{1}{%
            \pgfmathparse{0.3*((#1-0.35)^5 + 1.5*(#1-0.35)^4 - 3.8*(#1-0.35)^3 - 3.5*(#1-0.35)^2 + 3.4*(#1-0.35) + 5)}%
            }
            \begin{axis}[
                axis lines = middle,
                xlabel = $\xi$,
                ylabel = {$\eta$},
                domain=-2.15:2.15,
                samples=200,
                width=16cm,
                height=9cm,
                xtick=\empty,
                ytick=\empty,
            ]
            \addplot [
                  thick,
                  oiBlue,
                  domain=-2.075:2.075,
            ] {f(x)};
            \node at (axis cs:0,1.9) [anchor=south west] {$\tau$};
            \addplot [
                oiOrange,
                semithick,
                domain=-2.15:2.15,
            ] {1.9};
            \node at (axis cs:0,2.2) [anchor=south west] {$\tau + r$};
            \addplot [
                oiOrange,
                semithick,
                dashed,
                domain=-2.15:2.15,
            ] {2.2};
            \node at (axis cs:0,1.6) [anchor=south west] {$\tau - r$};
            \addplot [
                oiOrange,
                semithick,
                dashed,
                domain=-2.15:2.15,
            ] {1.6};
            \node at (axis cs:-1.56,0) [anchor=north west] {$\xi_1$};
            \addplot [
                black,
                semithick,
                dashed
            ] coordinates {(-1.56,-0.15) (-1.56,2.55)}; 
            \node at (axis cs:-0.5,0) [anchor=north west] {$\xi_2$};
            \addplot [
                black,
                semithick,
                dashed
            ] coordinates {(-0.5,-0.15) (-0.5,2.55)}; 
            \node at (axis cs:0.69,0) [anchor=north west] {$\xi_3$};
            \addplot [
                black,
                semithick,
                dashed
            ] coordinates {(0.69,-0.15) (0.69,2.55)}; 
            \node at (axis cs:1.58,0) [anchor=north west] {$\xi_4$};
            \addplot [
                black,
                semithick,
                dashed
            ] coordinates {(1.58,-0.15) (1.58,2.55)}; 
            \node at (axis cs:-1.95,0) [anchor=north west] {$- R$};
            \addplot [
                black,
                semithick,
                dashed
            ] coordinates {(-1.95,-0.15) (-1.95,2.55)}; 
            \node at (axis cs:2.01,0) [anchor=north west] {$R$};
            \addplot [
                black,
                semithick,
                dashed
            ] coordinates {(2.01,-0.15) (2.01,2.55)}; 
            \pgfmathparse{f(-1.95)}
            \edef\tempA{\pgfmathresult}
            \node at (axis cs:-1.95,\tempA) [anchor=north west, text=oiBlue] {$\phi_0$};
            \pgfmathparse{f(-1.0)}
            \edef\tempB{\pgfmathresult}
            \node at (axis cs:-1.0,\tempB) [anchor=north east, text=oiBlue] {$\phi_1$};
            \pgfmathparse{f(0.1)}
            \edef\tempC{\pgfmathresult}
            \node at (axis cs:0.1,\tempC) [anchor=north west, text=oiBlue] {$\phi_2$};
            \pgfmathparse{f(1.2)}
            \edef\tempD{\pgfmathresult}
            \node at (axis cs:1.2,\tempD) [anchor=north east, text=oiBlue] {$\phi_3$};
            \pgfmathparse{f(1.75)}
            \edef\tempE{\pgfmathresult}
            \node at (axis cs:1.75-0.02,\tempE+0.02) [anchor=north west, text=oiBlue] {$\phi_4$};
            \pgfmathparse{f(-1.95-0.02)}
            \edef\tempF{\pgfmathresult}
            \pgfmathparse{f(-2.085+0.02)}
            \edef\tempG{\pgfmathresult}
            \draw[black] (axis cs:-2.085+0.02,\tempG) rectangle (axis cs:-1.95-0.02,\tempF);
            \pgfmathparse{f(-1.95+0.02)}
            \edef\tempH{\pgfmathresult}
            \pgfmathparse{f(-1.7925-0.02)}
            \edef\tempI{\pgfmathresult}
            \draw[black] (axis cs:-1.95+0.02,\tempH) rectangle (axis cs:-1.7925-0.02,\tempI);
            \pgfmathparse{f(-1.817+0.02)}
            \edef\tempJ{\pgfmathresult}
            \pgfmathparse{f(-1.56)+0.02}
            \edef\tempK{\pgfmathresult}
            \draw[black] (axis cs:-1.817+0.02,\tempJ) rectangle (axis cs:-1.215-0.02,\tempK);
            \pgfmathparse{f(-1.23+0.02)}
            \edef\tempL{\pgfmathresult}
            \pgfmathparse{f(-0.5-0.02)}
            \edef\tempM{\pgfmathresult}
            \draw[black] (axis cs:-1.23+0.02,\tempL) rectangle (axis cs:-0.5-0.02,\tempM);
            \pgfmathparse{f(-0.5+0.02)}
            \edef\tempN{\pgfmathresult}
            \pgfmathparse{f(0.45-0.02)}
            \edef\tempO{\pgfmathresult}
            \draw[black] (axis cs:-0.5+0.02,\tempN) rectangle (axis cs:0.45-0.02,\tempO);
            \pgfmathparse{f(0.47+0.02)}
            \edef\tempP{\pgfmathresult}
            \pgfmathparse{f(0.69)+0.02}
            \edef\tempQ{\pgfmathresult}
            \draw[black] (axis cs:0.47+0.02,\tempP) rectangle (axis cs:0.9-0.02,\tempQ);
            \pgfmathparse{f(0.923+0.02)}
            \edef\tempR{\pgfmathresult}
            \pgfmathparse{f(1.58-0.02)}
            \edef\tempS{\pgfmathresult}
            \draw[black] (axis cs:0.923+0.02,\tempR) rectangle (axis cs:1.58-0.02,\tempS);
            \pgfmathparse{f(1.58+0.02)}
            \edef\tempT{\pgfmathresult}
            \pgfmathparse{f(1.916-0.02)}
            \edef\tempU{\pgfmathresult}
            \draw[black] (axis cs:1.58+0.02,\tempT) rectangle (axis cs:1.916-0.02,\tempU);
            \pgfmathparse{f(1.886+0.02)}
            \edef\tempV{\pgfmathresult}
            \pgfmathparse{f(2.0158-0.02)}
            \edef\tempW{\pgfmathresult}
            \draw[black] (axis cs:1.886+0.02,\tempV) rectangle (axis cs:2.0158-0.02,\tempW);
            \pgfmathparse{f(2.003+0.02)}
            \edef\tempX{\pgfmathresult}
            \pgfmathparse{f(2.088-0.02)}
            \edef\tempY{\pgfmathresult}
            \draw[black] (axis cs:2.003+0.02,\tempX) rectangle (axis cs:2.088-0.02,\tempY);
            \node at (axis cs:-1.9,-0.4) [anchor=center] {\small $(b)$};
            \node at (axis cs:-1.75,0.82) [anchor=center] {\small $(a)$};
            \node at (axis cs:-1.17,2.05) [anchor=center] {\small $(e)$};
            \node at (axis cs:-1.13,0.58) [anchor=center] {\small $(a)$};
            \node at (axis cs:0.35,0.58) [anchor=center] {\small $(a)$};
            \node at (axis cs:0.95,1.75) [anchor=center] {\small $(f)$};
            \node at (axis cs:1.03,0.82) [anchor=center] {\small $(a)$};
            \node at (axis cs:1.83,0.82) [anchor=center] {\small $(a)$};
            \node at (axis cs:1.83,2.05) [anchor=center] {\small $(c)$};
            \node at (axis cs:1.95,2.68) [anchor=center] {\small $(b)$};
          \end{axis}
        \end{tikzpicture}
    \end{figure}
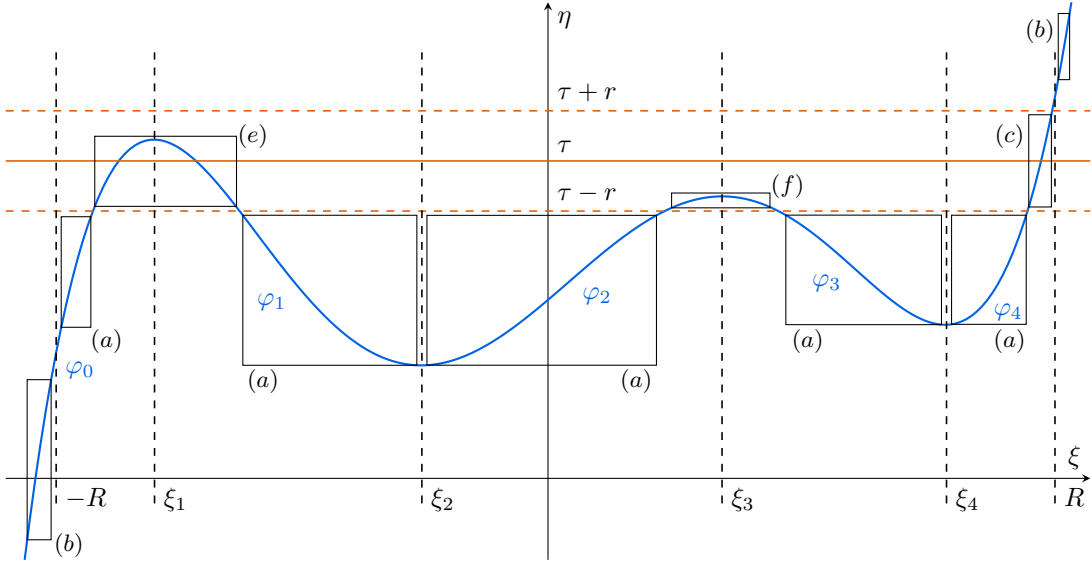
    We now enumerate all nontrivial cases and perform the necessary estimates, assuming here without loss of generality that $z \geq 0$.

    \textbf{(b) The continuous oscillatory tail case.} 
    For high frequencies it suffices to estimate each tail individually, but for low frequencies there is little oscillation, and instead cancellation between the two tails has to be exploited. We assume without loss of generality that $\sigma_0 = -1$ and distinguish the cases $\sigma_N = 1$ and $\sigma_N = -1$. 
    If $\sigma_N = 1$, then we have to bound the quantity
    \begin{align*}
        \int_R^\infty e^{i z \psi_0(\eta)} \frac{a(\psi_0(\eta))}{\eta - \tau} \dd \eta 
        - \int_R^\infty e^{i z \psi_N(\eta)} \frac{a(\psi_N(\eta))}{\eta - \tau} \dd \eta \,.
    \end{align*}
    On the other hand, if $\sigma_N = - 1$ we have to bound the quantity
    \begin{align*}
        \int_R^\infty e^{i z \psi_0(\eta)} \frac{a(\psi_0(\eta))}{\eta - \tau} \dd \eta 
        + \int_{-\infty}^{-R} e^{i z \psi_N(\eta)} \frac{a(\psi_N(\eta))}{\eta - \tau} \dd \eta \,.
    \end{align*}
    We can simultaneously treat both cases by writing
    \begin{align*}
        \int_R^\infty e^{i z \psi_0(\eta)} \frac{a(\psi_0(\eta))}{\eta - \tau} \dd \eta 
        - \int_R^\infty e^{i z \psi_N(\sigma_N \eta)} \frac{a(\psi_N(\sigma_N \eta))}{\eta - \sigma_N \tau} \dd \eta \,.
    \end{align*}
    Since we are in the case where there is no singularity in the oscillatory tail, we have $\tau < R - r$ and $\sigma_N \tau < R - r$. Then for all $\eta > R $ we have
    \begin{align} \label{eqn:850}
        \left| \frac{1}{\eta - \tau} \right| &\leq \left( 1 + \mathds{1}_{\{\tau \geq 0\}} \left| \frac{1}{\frac{\eta}{\tau} - 1} \right| \right) \frac{1}{\eta} 
        \leq \left( 1 + \frac{1}{\frac{R}{R-r} - 1} \right) \frac{1}{\eta}
        \lesssim_{R,r} \frac{1}{\eta}
        \\ \label{eqn:851}
        \left| \frac{1}{\eta - \sigma_N \tau} \right| &\leq \left( 1 + \mathds{1}_{\{\sigma_N \tau \geq 0\}} \left| \frac{1}{\frac{\eta}{\sigma_N \tau} - 1} \right| \right) \frac{1}{\eta} 
        \leq \left( 1 + \frac{1}{\frac{R}{R-r} - 1} \right) \frac{1}{\eta}
        \lesssim_{R,r} \frac{1}{\eta} \,.
    \end{align}
    We decompose
    \begin{align*}
        & \int_R^\infty e^{i z \psi_0(\eta)} \frac{a(\psi_0(\eta))}{\eta - \tau} \dd \eta 
        - \int_R^\infty e^{i z \psi_N(\sigma_N \eta)} \frac{a(\psi_N(\sigma_N \eta))}{\eta - \sigma_N \tau} \dd \eta
        \\ &= \int_R^\infty e^{i z \psi_N(\sigma_N \eta)} \left( \frac{a(\psi_0(\eta))}{\eta - \tau} - \frac{a(\psi_N(\sigma_N \eta))}{\eta - \sigma_N \tau} \right) \dd \eta
        + \int_R^\infty \left( e^{i z \psi_0(\eta)} - e^{i z \psi_N(\sigma_N \eta)} \right) \frac{a(\psi_0(\eta))}{\eta - \tau} \dd \eta
        \\ &= (I) + (II)
    \end{align*}
    and directly estimate
    \begin{align*}
        |(I)| &\lesssim_A \int_R^\infty \left| \frac{\sigma_N \tau - \tau}{(\eta - \tau) (\eta - \sigma_N \tau)} \right| + \left| \frac{a(\psi_0(\eta)) - a(\psi_N(\sigma_N \eta))}{(\eta - \sigma_N \tau)} \right| \dd \eta
        \\ &\lesssim_{A,R,P,r} \int_R^\infty \frac{1}{\eta^2} + \frac{(\eta^\beta)^{- 1 - \delta} \eta^{\beta - \delta}}{\eta} \dd \eta
        \lesssim_{R,\beta,\delta} 1 \,.
    \end{align*}
    If $z = 0$ then $(II) = 0$, so we assume $z \neq 0$. We split the interval $(R, \infty)$ into $(R, R \lor z^{- \frac{1}{\beta}})$ and $(R \lor z^{- \frac{1}{\beta}}, \infty)$. 
    For the finite interval, we estimate
    \begin{align*}
        \left| \int_R^{R \lor z^{- \frac{1}{\beta}}} \left( e^{i z \psi_0(\eta)} - e^{i z \psi_N(\sigma_N \eta)} \right) \frac{a(\psi_0(\eta))}{\eta - \tau} \dd \eta \right| 
        &\lesssim_A z \int_R^{R \lor z^{- \frac{1}{\beta}}} \frac{|\psi_0(\eta)| + |\psi_N(\sigma_N \eta)|}{\eta - \tau} \dd \eta
        \\ &\lesssim_{R,P,r} \, z \int_R^{R \lor z^{- \frac{1}{\beta}}} \frac{\eta^\beta}{\eta} \dd \eta
        \\ &= \, z \frac{1}{\beta} \left( (R \lor z^{- \frac{1}{\beta}})^\beta - R^\beta \right)
        \lesssim_\beta 1 \,.
    \end{align*}
    We decompose the other interval into dyadic blocks
    \begin{align*}
        [R \lor z^{-\frac{1}{\beta}}, \infty) &= \bigcup_{n=L}^\infty [2^l, 2^{l+1}) \cup [R \lor z^{-\frac{1}{\beta}}, 2^L) \,,
    \end{align*}
    where $L \in \N$ is the unique integer such that $2^{L-1} \leq R \lor z^{-\frac{1}{\beta}} < 2^L$.
    Now we use the oscillation with Lemma \ref{lem:999} to estimate
    \begin{align*}
        &\sum_{l=L}^\infty \left| \int_{2^l}^{2^{l+1}} \left( e^{i z \psi_0(\eta)} - e^{i z \psi_N(\sigma_N \eta)} \right) \frac{a(\psi_0(\eta))}{\eta - \tau} \dd \eta \right|
        \\ &+ \left| \int_{R \lor z^{-\frac{1}{\beta}}}^{2^L} \left( e^{i z \psi_0(\eta)} - e^{i z \psi_N(\sigma_N \eta)} \right) \frac{a(\psi_0(\eta))}{\eta - \tau} \dd \eta \right|
        \\ &\lesssim_{M,R,A,r,p} \frac{1}{z} \left( \sum_{l=L}^\infty 2^{(l+1) (1 - \beta)} 2^{- l} + 2^{L (1 - \beta)} (R \lor z^{-\frac{1}{\beta}})^{-1} \right)
        \\ &\lesssim \left( \frac{1}{z} \sum_{l=L}^\infty 2^{- \beta l} + z^{- 1} 2^{(1 - \beta) L} (R \lor z^{-\frac{1}{\beta}})^{-1} \right)
        \\ &\lesssim \left( \left( 2^L (R \lor z^{-\frac{1}{\beta}})^{-1} \right)^{- \beta} + \left( 2^L (R \lor z^{-\frac{1}{\beta}})^{-1} \right)^{1 - \beta} \right)
        \lesssim 1 \,.
    \end{align*}
    Here (H5) implies the monotonicity condition that the Lemma requires.

    \textbf{(c) The non-critical point singularity bounded case.}
    Due to (a) it suffices to estimate the integral on $(\tau - r, \tau + r)$. 
    In fact, the interval $(\tau - r', \tau + r')$ where $r' = \frac{r}{2}$ suffices, since the case (iii) takes care of the integral near a critical point. 
    We have ensured that we are not close to a critical point and hence $\sup_{\eta \in (\tau - r', \tau + r')} |\psi_j'(\eta)|$ is bounded uniformly for small $\tau$, 
    i.e. $|\tau| \leq R$. On the other hand, (H2) ensures boundedness for $|\tau| > R$. 
    Another control quantity is $\sup_{\eta \in (\tau - r', \tau + r')} \frac{1}{|\psi_j'(\eta)|}$, which is bounded uniformly in $\tau$ for $|\tau| \leq R$ only. 
    We have
    \begin{align*}
        \int_{\tau - r'}^{\tau + r'} e^{i z \psi_j(\eta)} \frac{\mathds{1}_{\{|\eta - \tau| > \epsilon\}} a(\psi_j(\eta))}{\eta - \tau} \dd \eta
        &= \int_\epsilon^{r'} \left( e^{i z \psi_j(\tau + \eta)} a(\psi_j(\tau + \eta)) - e^{i z \psi_j(\tau - \eta)} a(\psi_j(\tau - \eta)) \right) \frac{1}{\eta} \dd \eta \,.
    \end{align*}
    The Lipschitz property yields
    \begin{align*}
        & \left| \int_\epsilon^{r' \land \frac{1}{z}} \left( e^{i z \psi_j(\tau + \eta)} a(\psi_j(\tau + \eta)) - e^{i z \psi_j(\tau - \eta)} a(\psi_j(\tau - \eta)) \right) \frac{1}{\eta} \dd \eta \right|
        \lesssim_{A,R,\delta,\beta} \sup_{\eta \in (\tau - r', \tau + r')} |\psi_j'(\eta)| \,.
    \end{align*}
    If $z r' \leq 1$, then this is the whole integral, so we assume $z r' > 1$ and aim to estimate
    \begin{align*}
        & \int_{\frac{1}{z}}^{r'} \left( e^{i z \psi_j(\tau + \eta)} a(\psi_j(\tau + \eta)) - e^{i z \psi_j(\tau - \eta)} a(\psi_j(\tau - \eta)) \right) \frac{1}{\eta} \dd \eta \,.
    \end{align*}
    Here we apply Lemma \ref{lem:999} to each summand, using that $\psi_j$ is strictly monotonic and that $\frac{a(\psi_j(\tau \pm \eta))}{\eta}$ only changes monotonicity less than $M$ times by (H5).
    This yields
    \begin{align*}
        \left| \int_{\frac{1}{z}}^{r'} e^{i z \psi_j(\tau \pm \eta)} a(\psi_j(\tau \pm \eta)) \frac{1}{\eta} \dd \eta \right|
        &\lesssim_A \sup_{\eta \in (\tau - r', \tau + r')} \frac{1}{|\psi_j'(\eta)|} \,.
    \end{align*}
    Note that the right hand side may grow with $\tau$. This is acceptable as we are in the case $j \not\in \{0, N\}$, i.e. there exists some fixed large $R$ such that $|\tau| < R$, but it needs to be improved for the case (d).
    
    \textbf{(d) The non-critical point singularity oscillatory tail case.} 
    In this case the singularity is far from the origin in the oscillatory tail, i.e. $|\tau| > R$. 
    Note that now $j \in \{0, N\}$. 
    We assume without loss of generality that $\tau > \frac{R}{1 - r} > R$, and correspondingly that $\sigma_j = - 1$ if $j = 0$, and $\sigma_j = 1$ if $j = N$, as otherwise there is no singularity in the oscillatory tail.
    We would like to separate the singularity from the oscillatory tail. This requires us to revisit the proofs of (b) and (c). 
    Note that in (c) we have treated the singularity on an interval $(\tau - r', \tau + r')$. 
    We now instead use the interval $(\tau - \tau r, \tau + \tau r) \subset (R, \infty)$. 
    We decompose
    \begin{align*}
        & \int_\epsilon^{\tau r} \left( e^{i z \psi_j(\tau + \eta)} a(\psi_j(\tau + \eta)) - e^{i z \psi_j(\tau - \eta)} a(\psi_j(\tau - \eta)) \right) \frac{1}{\eta} \dd \eta
        \\ &= \int_\epsilon^{\tau r} e^{i z \psi_j(\tau + \eta)} (a(\psi_j(\tau + \eta)) - a(\psi_j(\tau - \eta))) \frac{1}{\eta} \dd \eta
        + \int_\epsilon^{\tau r} \left( e^{i z \psi_j(\tau + \eta)} - e^{i z \psi_j(\tau - \eta)} \right) a(\psi_j(\tau - \eta)) \frac{1}{\eta} \dd \eta \,.
    \end{align*}
    The Lipschitz property yields
    \begin{align*}
        \left| \int_\epsilon^{\tau r} e^{i z \psi_j(\tau + \eta)} (a(\psi_j(\tau + \eta)) - a(\psi_j(\tau - \eta))) \frac{1}{\eta} \dd \eta \right|
        \lesssim_{A,P} \tau r \sup_{\eta \in (\tau - \tau r, \tau + \tau r)} (\eta^\beta)^{-1-\delta} \eta^{\beta - 1} \lesssim_R 1
    \end{align*}
    and
    \begin{align*}
        \left| \int_\epsilon^{\tau r \land \frac{\tau^{1-\beta}}{z}} \left( e^{i z \psi_j(\tau + \eta)} - e^{i z \psi_j(\tau - \eta)} \right) a(\psi_j(\tau - \eta)) \frac{1}{\eta} \dd \eta \right|
        &\lesssim_{A,P} z \frac{\tau^{1-\beta}}{z} \sup_{\eta \in (\tau - \tau r, \tau + \tau r)} \eta^{\beta - 1} 
        \lesssim 1 \,.
    \end{align*}
    If $\frac{\tau^{1-\beta}}{z} < \tau r$ then there is a remainder, which we estimate using the oscillation, i.e. Lemma \ref{lem:999}:
    \begin{align*}
        \left| \int_{\frac{\tau^{1-\beta}}{z}}^{\tau r} \left( e^{i z \psi_j(\tau + \eta)} - e^{i z \psi_j(\tau - \eta)} \right) a(\psi_j(\tau - \eta)) \frac{1}{\eta} \dd \eta \right|
        &\lesssim_{M,A,P} \frac{1}{z} \sup_{\eta \in (\tau - \tau r, \tau + \tau r)} \eta^{1 - \beta} \frac{z}{\tau^{1-\beta}}
        \\ &\lesssim 1 \,.
    \end{align*}
    As before, (H5) supplies the required monotonicity assumption.
    We now modify the proof of (b) to deal with the remaining oscillatory tail, where the interval $(\tau - \tau r, \tau + \tau r)$ around the singularity has been removed. 
    We again assume that $\sigma_0 = - 1$ and distinguish two cases.
    If $\sigma_N = 1$ we have to bound the quantity
    \begin{align*}
        & \int_R^{\tau - \tau r} e^{i z \psi_0(\eta)} \frac{a(\psi_0(\eta))}{\eta - \tau} \dd \eta 
        - \int_R^{\tau - \tau r} e^{i z \psi_N(\eta)} \frac{a(\psi_N(\eta))}{\eta - \tau} \dd \eta
        \\ &+ \int_{\tau + \tau r}^\infty e^{i z \psi_0(\eta)} \frac{a(\psi_0(\eta))}{\eta - \tau} \dd \eta 
        - \int_{\tau + \tau r}^\infty e^{i z \psi_N(\eta)} \frac{a(\psi_N(\eta))}{\eta - \tau} \dd \eta \,,
    \end{align*}
    while if $\sigma_N = - 1$ we have to bound the quantity
    \begin{align*}
        & \int_R^{\tau - \tau r} e^{i z \psi_0(\eta)} \frac{a(\psi_0(\eta))}{\eta - \tau} \dd \eta 
        + \int_{- \tau + \tau r}^{- R} e^{i z \psi_N(\eta)} \frac{a(\psi_N(\eta))}{\eta - \tau} \dd \eta
        \\ &+ \int_{\tau + \tau r}^\infty e^{i z \psi_0(\eta)} \frac{a(\psi_0(\eta))}{\eta - \tau} \dd \eta 
        + \int_{- \infty}^{- \tau - \tau r} e^{i z \psi_N(\eta)} \frac{a(\psi_N(\eta))}{\eta - \tau} \dd \eta 
        \\ &+ \int_{- \tau - \tau r}^{- \tau + \tau r} e^{i z \psi_N(\eta)} \frac{a(\psi_N(\eta))}{\eta - \tau} \dd \eta \,.
    \end{align*}
    We can directly estimate all the integrals over finite intervals, as the length of the intervals is always bounded by $\tau + \tau r$, while the integrand is bounded by $\frac{A}{\tau r}$.
    Unifying the remaining integrals for $\sigma_N \in \{-1, 1\}$ as before, it remains to bound the quantity
    \begin{align*}
        & \int_{\tau + \tau r}^\infty e^{i z \psi_0(\eta)} \frac{a(\psi_0(\eta))}{\eta - \tau} \dd \eta - \int_{\tau + \tau r}^\infty e^{i z \psi_N(\sigma_N \eta)} \frac{a(\psi_N(\sigma_N \eta))}{\eta - \sigma_N \tau} \dd \eta
        \\ &= \int_{\tau + \tau r}^\infty e^{i z \psi_N(\sigma_N \eta)} \left( \frac{a(\psi_0(\eta))}{\eta - \tau} - \frac{a(\psi_N(\sigma_N \eta))}{\eta - \sigma_N \tau} \right) \dd \eta
        + \int_{\tau + \tau r}^\infty \left( e^{i z \psi_0(\eta)} - e^{i z \psi_N(\sigma_N \eta)} \right) \frac{a(\psi_0(\eta))}{\eta - \tau} \dd \eta
        \\ &= (I) + (II) \,.
    \end{align*}
    This requires adapting \eqref{eqn:850}--\eqref{eqn:851}. For all $\eta \in (\tau + \tau r, \infty)$ we have
    \begin{align*}
        \left| \frac{1}{\eta - \tau} \right| &\leq \left( 1 + \left| \frac{1}{\frac{\eta}{\tau} - 1} \right| \right) \frac{1}{\eta} 
        \leq \left( 1 + \left| \frac{1}{\frac{\tau + \tau r}{\tau} - 1} \right| \right) \frac{1}{\eta}
        \leq \left( 1 + \frac{1}{r} \right) \frac{1}{\eta}
        & \left| \frac{1}{\eta + \tau} \right| &\leq \frac{1}{\eta} \,.
    \end{align*}
    Both terms $(I)$ and $(II)$ can now be estimated with the same method as in (b).
    
    \textbf{(iii) The critical point case}
    Outside the interval $(\tau - r, \tau + r)$ we have $|\eta - \tau| \geq r$, so the integrals over these sections can be estimaed with (a) if the interval is finite and (b) if not.
    It therefore suffices to consider the integral over $(\tau - r, \tau + r)$. Here we shall describe our approach to the cases (e), (f) and (g) with the help of Figure \ref{fig:2}.
    We assume without loss of generality that $\sigma_{j-1} = 1$.

    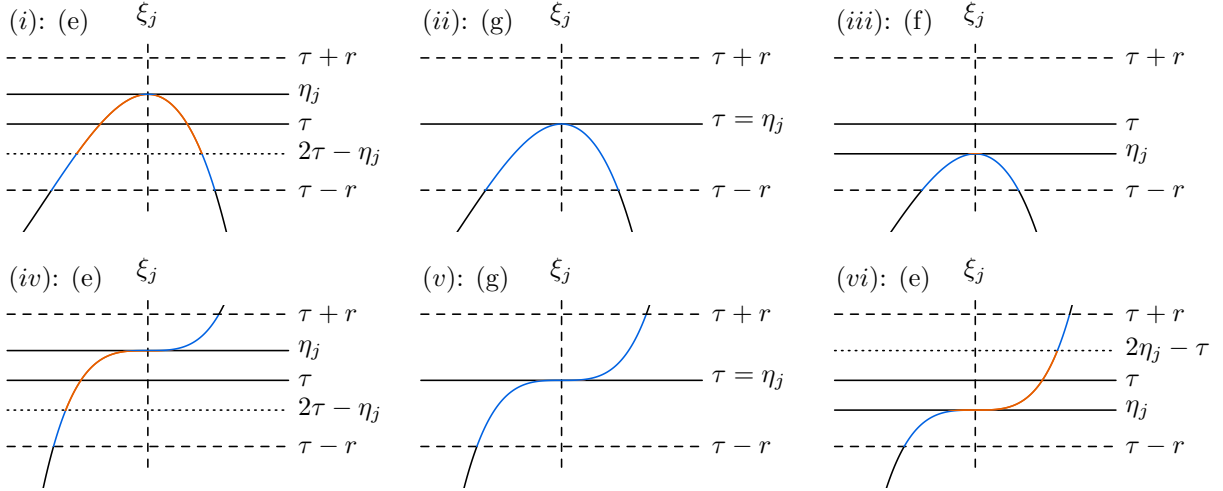
\begin{figure}[H]
        \caption{The 6 possibilities for cases (e), (f) and (g) when $\sigma_{j-1} = 1$.}
        \centering

    \definecolor{oiBlue}   {HTML}{0062E2}
    \definecolor{oiOrange} {HTML}{E55E00}
    
    \begin{tikzpicture}[
        scale=1.1,
        line cap=round,
        line join=round,
        semithick,
        samples=160,
        declare function={
        tauVal=1.0;
        rVal=0.8;
        apex=tauVal+rVal;
        f(\x)=apex -30*((0.2*\x+0.75)^4 - (0.2*\x+0.75)^3 - 0.75^4 + 0.75^3);
        g(\x)= + 0.75^3 * \x^3 + 0.75^4 * \x^4 + tauVal;
        }
    ]
        \def\xmin{-1.7}
        \def\xmax{1.7}
    
        \newcommand*\DrawSegF[4]{%
        \begin{scope}
            \clip (\xmin,{#1}) rectangle (\xmax,{#2});
            \draw[#3] plot[domain=\xmin:\xmax] (\x,{f(\x)-rVal+\bVal});
        \end{scope}
        }
    
        \newcommand*\DrawSegG[4]{%
        \begin{scope}
            \clip (\xmin,{#1}) rectangle (\xmax,{#2});
            \draw[#3] plot[domain=\xmin:\xmax] (\x,{g(\x)+\bVal});
        \end{scope}
        }
    
        \newcommand\PanelF[4]{%
        \pgfmathsetmacro\bVal{#2}%
        \begin{scope}[xshift=#1 cm]
            \draw[dashed] (0,{tauVal-rVal-0.25}) -- (0,{tauVal+rVal+0.25}) node[above]{$\xi_j$};
            \draw[dashed] (\xmin,{tauVal+rVal}) -- (\xmax,{tauVal+rVal}) node[right]{$\tau+r$};
            \draw[dashed] (\xmin,{tauVal-rVal}) -- (\xmax,{tauVal-rVal}) node[right]{$\tau-r$};
            \ifdim\bVal pt=0pt
            \draw (\xmin,{tauVal}) -- (\xmax,{tauVal}) node[right]{$\tau=\eta_j$};
            \else
            \draw (\xmin,{tauVal}) -- (\xmax,{tauVal}) node[right]{$\tau$};
            \draw (\xmin,{tauVal+\bVal}) -- (\xmax,{tauVal+\bVal}) node[right]{$\eta_j$};
            \fi
            \ifdim\bVal pt>0pt
            \draw[dotted] (\xmin,{tauVal-\bVal}) -- (\xmax,{tauVal-\bVal}) node[right]{$2\tau-\eta_j$};
            \fi
            \DrawSegF{-0.3}{tauVal-rVal}{black}{left}
            \DrawSegF{tauVal-rVal}{tauVal-abs(\bVal)}{oiBlue}{left}
            \DrawSegF{tauVal-\bVal}{tauVal+\bVal}{oiOrange}{above left}
            \DrawSegF{tauVal-\bVal}{tauVal+\bVal}{oiOrange}{above right}
            \DrawSegF{tauVal+abs(\bVal)}{tauVal+rVal}{oiBlue}{right}
            \DrawSegF{tauVal+rVal}{apex+0.1}{black}{right}
            \node[anchor=north west] at (\xmin-0.1,{apex+0.7}) {#4: #3};
        \end{scope}
        }
    
        \newcommand\PanelG[4]{%
        \pgfmathsetmacro\bVal{#2}%
        \begin{scope}[xshift=#1 cm]
            \draw[dashed] (0,{tauVal-rVal-0.25}) -- (0,{tauVal+rVal+0.25}) node[above]{$\xi_j$};
            \draw[dashed] (\xmin,{tauVal+rVal}) -- (\xmax,{tauVal+rVal}) node[right]{$\tau+r$};
            \draw[dashed] (\xmin,{tauVal-rVal}) -- (\xmax,{tauVal-rVal}) node[right]{$\tau-r$};
            \ifdim\bVal pt=0pt
            \draw (\xmin,{tauVal}) -- (\xmax,{tauVal}) node[right]{$\tau=\eta_j$};
            \else
            \draw (\xmin,{tauVal}) -- (\xmax,{tauVal}) node[right]{$\tau$};
            \draw (\xmin,{tauVal+\bVal}) -- (\xmax,{tauVal+\bVal}) node[right]{$\eta_j$};
            \fi
            \ifdim\bVal pt>0pt
            \draw[dotted] (\xmin,{tauVal-\bVal}) -- (\xmax,{tauVal-\bVal}) node[right]{$2\tau-\eta_j$};
            \fi
            \ifdim\bVal pt<0pt
            \draw[dotted] (\xmin,{tauVal-\bVal}) -- (\xmax,{tauVal-\bVal}) node[right]{$2\eta_j-\tau$};
            \fi
            \DrawSegG{-0.3}{tauVal-rVal}{black}{left}
            \DrawSegG{tauVal-rVal}{tauVal-abs(\bVal)}{oiBlue}{left}
            \DrawSegG{tauVal-\bVal}{tauVal+\bVal}{oiOrange}{above left}
            \DrawSegG{tauVal-\bVal}{tauVal+\bVal}{oiOrange}{above right}
            \DrawSegG{tauVal+abs(\bVal)}{tauVal+rVal}{oiBlue}{right}
            \DrawSegG{tauVal+rVal}{apex+0.1}{black}{right}
            \node[anchor=north west] at (\xmin-0.1,{apex+0.7}) {#4: #3};
        \end{scope}
        }
    
        \PanelF{0}{0.36}{(e)}{$(i)$}
        \PanelF{5}{0.00}{(g)}{$(ii)$}
        \PanelF{10}{-0.36}{(f)}{$(iii)$}
    
        \begin{scope}[yshift=-3.1cm]
        \PanelG{0}{0.36}{(e)}{$(iv)$}
        \PanelG{5}{0.00}{(g)}{$(v)$}
        \PanelG{10}{-0.36}{(e)}{$(vi)$}
        \end{scope}
    
    \end{tikzpicture}
    
    \label{fig:2}
    \end{figure}
    We refer to different parts of the integral to estimate as ``sections''.
    The orange sections contain the singularity, so cancellation of the sections before and after the singularity is crucial.
    Note that in $(i)$ cancellation between the orange section before $\xi_j$ and the one after $\xi_j$ in is not necessary.
    We perform the proof only for the section to the right of $\xi_j$ in $(i)$, as it is representative of all other orange sections.
    All blue sections need to be considered in pairs consisting of a section to the left of $\xi_j$ and one to the right, in order to exploit a cancellation of the form $|e^{i x} - e^{i y}| \leq |x| + |y|$.
    Here the cases $(ii)$ and $(v)$ are strictly more difficult than the other ones, because the sections we integrate over go all the way up to the singularity.
    We treat them as representative of all the other blue sections.

    Define $b = \eta_j - \tau$.
    In the case $(i)$ we have $\tau - r < \tau - b < \tau < \tau + b = \eta_j < \tau + r$, and the integral to estimate is
    \begin{align*}
        & \int_{\tau - b}^{\tau + b} e^{i z \psi_{j-1}(\eta)} \frac{\mathds{1}_{\{|\eta - \tau| > \epsilon\}} a(\psi_{j-1}(\eta))}{\eta - \tau} \dd \eta 
        + \int_{\tau - b}^{\tau + b} e^{i z \psi_j(\eta)} \frac{\mathds{1}_{\{|\eta - \tau| > \epsilon\}} a(\psi_j(\eta))}{\eta - \tau} \dd \eta
    \end{align*}
    As mentioned before, we estimate these two terms individually and with the same technique, so it suffices to consider the case of $\psi_j$.
    Here
    \begin{align*}
        & \int_{\tau - b}^{\tau + b} e^{i z \psi_j(\eta)} \frac{\mathds{1}_{\{|\eta - \tau| > \epsilon\}} a(\psi_j(\eta))}{\eta - \tau} \dd \eta
        \\ &= \int_\epsilon^b \left( e^{i z \psi_j(\tau + \eta)} a(\psi_j(\tau + \eta)) - e^{i z \psi_j(\tau - \eta)} a(\psi_j(\tau - \eta)) \right) \frac{1}{\eta} \dd \eta
        \\ &= \int_\epsilon^b \left( e^{i z \psi_j(\eta_j - (b - \eta))} a(\psi_j(\eta_j - (b - \eta))) - e^{i z \psi_j(\eta_j - (b + \eta))} a(\psi_j(\eta_j - (b + \eta))) \right) \frac{1}{\eta} \dd \eta \,.
    \end{align*}
    To simplify notation, we assume without loss of generality that $\eta_j = 0$ and replace $\psi_j(\eta)$ by $\psi_j(- \eta)$. 
    This means that $\psi_j(\eta_j - (b \pm \eta))$ is replaced by $\psi_j(b \pm \eta)$.
    We decompose
    \begin{align*}
        & \int_\epsilon^b \frac{e^{i z \psi_j(b + \eta)} a(\psi_j(b + \eta)) - e^{i z \psi_j(b - \eta)} a(\psi_j(b - \eta))}{\eta} \dd \eta
        \\ &= \int_\epsilon^b e^{i z \psi_j(b - \eta)} \frac{a(\psi_j(b + \eta)) - a(\psi_j(b - \eta))}{\eta} \dd \eta + \int_\epsilon^b \frac{e^{i z \psi_j(b + \eta)} - e^{i z \psi_j(b - \eta)}}{\eta} a(\psi_j(b + \eta)) \dd \eta \,.
    \end{align*}
    The first term can be directly estimated:
    \begin{align*}
        & \int_\epsilon^b \frac{|a(\psi_j(b + \eta)) - a(\psi_j(b - \eta))|}{\eta} \dd \eta \leq \int_0^b \frac{1}{\eta} \int_{-\eta}^{\eta} |a'(\psi_j(b + y)) \psi_j'(b + y)| \dd y \dd \eta 
        \\ &= \int_{-b}^b |a'(\psi_j(b + y)) \psi_j'(b + y)| \int_{|y|}^b \frac{1}{\eta} \dd \eta \dd y
        \lesssim_{A,Q} \int_{-b}^b (b + y)^{\alpha-1} \ln\left(\frac{b}{|y|}\right) \dd y 
        \\ &\leq b^\alpha \int_0^2 \frac{1}{y^{1-\alpha}} \ln\left(\frac{1}{|y - 1|}\right) \dd y 
        \lesssim_\alpha 1 \,.
    \end{align*}
    To estimate the second term ,we use the Van der Corput type Lemma \ref{lem:999} for the regime $\left(\epsilon, \frac{b^{1-\alpha}}{z}\right)$. 
    Specifically, suppose that $\frac{b^{1-\alpha}}{z} < b$.
    Then (H5) permits application of Lemma \ref{lem:999}, which yields
    \begin{align*}
        \left| \int_{\frac{b^{1-\alpha}}{z}}^b e^{i z \psi_j(b \pm \eta)} \frac{a(\psi_j(b \pm \eta))}{\eta} \dd \eta \right| 
        &\lesssim_{M,A,q} \frac{1}{z} \sup_{\eta \in \left(\frac{b^{1-\alpha}}{z}, b\right)} (b \pm \eta)^{1 - \alpha} \sup_{\eta \in \left(\frac{b^{1-\alpha}}{z}, b\right)} \frac{1}{\eta}
        \\ &\leq \frac{1}{z} (2 b)^{1 - \alpha} \frac{z}{b^{1-\alpha}} \lesssim 1 \,.
    \end{align*}
    The remaining interval is $\left( \epsilon, \epsilon \lor \frac{b^{1-\alpha}}{z} \right)$, so we assume $\epsilon < \frac{b^{1-\alpha}}{z}$ and estimate
    \begin{align*}
        \left| \int_\epsilon^{\frac{b^{1-\alpha}}{z}} \frac{e^{i z \psi_j(b + \eta)} - e^{i z \psi_j(b - \eta)}}{\eta} a(\psi_j(b + \eta)) \dd \eta \right|
        &= \left| \int_\epsilon^{\frac{b^{1-\alpha}}{z}} \frac{a(\psi_j(b + \eta))}{\eta} \int_{-\eta}^\eta e^{i z \psi_j(b + y)} i z \psi_j'(b + y) \dd y \, \dd \eta \right|
        \\ &= \left| \int_{-\frac{b^{1-\alpha}}{z}}^{\frac{b^{1-\alpha}}{z}} e^{i z \psi_j(b + y)} i z \psi_j'(b + y) \int_{\epsilon \lor |y|}^{\frac{b^{1-\alpha}}{z}} \frac{a(\psi_j(b + \eta))}{\eta} \dd \eta \, \dd y \right|
        \\ &\lesssim_{A, Q} \frac{z}{b^{1-\alpha}} \int_{-\frac{b^{1-\alpha}}{z}}^{\frac{b^{1-\alpha}}{z}} \frac{b^{1-\alpha}}{(b + y)^{1-\alpha}} \left|\log\left(\frac{b^{1-\alpha}}{z} \frac{1}{|y|}\right)\right| \dd y
        \\ &= \int_{-1}^1 \frac{1}{\left(1 + \frac{t}{b^\alpha z}\right)^{1-\alpha}} \left|\log\left(\frac{1}{|t|}\right)\right| \dd t
        \lesssim_\alpha 1 \,.
    \end{align*}
    Note that at the end we have used $b^\alpha z > 1$. It remains to deal with the case $b^\alpha z \leq 1$:
    \begin{align*}
        \left| \int_\epsilon^b \frac{e^{i z \psi_j(b + \eta)} - e^{i z \psi_j(b - \eta)}}{\eta} a(\psi_j(b + \eta)) \dd \eta \right|
        &= \left| \int_{-b}^b e^{i z \psi_j(b + y)} i z \psi_j'(b + y) \int_{\epsilon \lor |y|}^b \frac{a(\psi_j(b + \eta))}{\eta} \dd \eta \, \dd y \right|
        \\ &\lesssim_{A, Q} \frac{1}{b} \int_{-b}^b \frac{b^{1-\alpha}}{(b + y)^{1-\alpha}} \left|\log\left(\frac{b}{|y|}\right)\right| \dd y
        \lesssim_\alpha 1 \,.
    \end{align*}

    In the cases $(ii)$ and $(v)$ we have $\tau - r < \tau = \eta_j < \tau + r$ and $\sigma_{j-1} = 1$ while $\sigma_j \in \{-1, 1\}$. 
    The integral to estimate can be written for both cases as
    \begin{align*}
        & \int_{\tau - r}^{\tau - \epsilon} e^{i z \psi_{j-1}(\eta)} \frac{a(\psi_{j-1}(\eta))}{\eta - \tau} \dd \eta 
        + \int_{\tau + \sigma_j \epsilon}^{\tau + \sigma_j r} e^{i z \psi_j(\eta)} \frac{a(\psi_j(\eta))}{\eta - \tau} \dd \eta
        \\ &= \int_\epsilon^r \left( e^{i z \psi_j(\eta_j + \sigma_j \eta)} - e^{i z \psi_{j-1}(\eta_j - \eta)} \right) \frac{a(\psi_j(\eta_j + \sigma_j \eta))}{\eta} \dd \eta
        \\ &+ \int_\epsilon^r e^{i z \psi_{j-1}(\eta_j - \eta)} \frac{a(\psi_j(\eta_j + \sigma_j \eta)) - a(\psi_{j-1}(\eta_j - \eta))}{\eta} \dd \eta \\
        &= (I) + (II) \,.
    \end{align*}
    We have
    \begin{align*}
        |(II)| &\lesssim_A \int_\epsilon^r \frac{|\psi_j(\eta_j + \sigma_j \eta) - \xi_j + \xi_j - \psi_{j-1}(\eta_j - \eta)|}{\eta} \dd \eta
        \lesssim_Q \int_0^r \frac{\eta^\alpha}{\eta} \dd \eta 
        \lesssim_\alpha 1 \,.
    \end{align*}
    It remains to estimate $(I)$.
    Assuming without loss of generality that $\epsilon < z^{-\frac{1}{\alpha}}$, we can immediately estimate
    \begin{align*}
        \left| \int_\epsilon^{z^{-\frac{1}{\alpha}}} \left( e^{i z \psi_j(\eta_j + \sigma_j \eta)} - e^{i z \psi_{j-1}(\eta_j - \eta)} \right) \frac{a(\psi_j(\eta_j + \sigma_j \eta))}{\eta} \dd \eta \right|
        &\lesssim_{A, Q} z \int_0^{z^{-\frac{1}{\alpha}}} \eta^\alpha \frac{1}{\eta} \dd \eta
        \\ &\lesssim_\alpha z (z^{-\frac{1}{\alpha}})^\alpha 
        = 1 \,.
    \end{align*}
    If $z^{-\frac{1}{\alpha}} \geq r$ then this is the whole integral, so we assume $z^{-\frac{1}{\alpha}} < r$ and decompose
    \begin{align*}
        [z^{-\frac{1}{\alpha}}, r) &= \bigcup_{l=0}^{L-1} [r 2^{-l-1}, r 2^{-l}) \cup [z^{-\frac{1}{\alpha}}, r 2^{-L}) \,,
    \end{align*}
    where $L \in \N$ is the unique integer such that $2^{-L-1} \leq \frac{z^{-\frac{1}{\alpha}}}{r} < 2^{-L}$.
    With (H5) we can use the oscillation through Lemma \ref{lem:999} to estimate
    \begin{align*}
        &\sum_{l=0}^{L-1} \left| \int_{r 2^{-l-1}}^{r 2^{-l}} \left( e^{i z \psi_j(\eta_j + \sigma_j \eta)} - e^{i z \psi_{j-1}(\eta_j - \eta)} \right)  \frac{a(\psi_j(\eta_j + \sigma_j \eta))}{\eta} \dd \eta \right|
        \\ &+ \left| \int_{z^{-\frac{1}{\alpha}}}^{r 2^{-L}} \left( e^{i z \psi_j(\eta_j + \sigma_j \eta)} - e^{i z \psi_{j-1}(\eta_j - \eta)} \right) \frac{a(\psi_j(\eta_j + \sigma_j \eta))}{\eta} \dd \eta \right|
        \\ &\lesssim_{M, A, q} \frac{1}{z} \left( \sum_{l=0}^{L-1} (r 2^{-l})^{1-\alpha} (r 2^{-l-1})^{-1} + (r 2^{-L})^{1-\alpha} (z^{-\frac{1}{\alpha}})^{-1} \right)
        \\ &\lesssim \left( \frac{1}{z} \sum_{l=0}^{L-1} 2^{\alpha l} + z^{\frac{1}{\alpha} - 1} 2^{(\alpha - 1) L} \right)
        \lesssim_\alpha \left( \left( 2^{- L} z^{\frac{1}{\alpha}} \right)^{- \alpha} + \left( 2^{- L} z^{\frac{1}{\alpha}} \right)^{1 - \alpha} \right)
        \lesssim_{r, \alpha} 1 \,.
    \end{align*}
\end{proof}

\subsection{The maximal function estimate}
Let $\phi \in C^3(\R; \R)$. We assume that there exists a large $R > 0$ for which the following hold.
\begin{enumerate}[(J1)]
    \item There exists some $A > 0$ and $n > 1$ such that 
    \begin{align*}
        |\phi'(\xi)| \geq R \Longrightarrow |\phi'(\xi)| \leq A |\xi|^{n-1} \,.
    \end{align*}
    \item There exists some $B > 0$ such that 
    \begin{align*}
        |\phi'(\xi)| \geq R \Longrightarrow |\phi'(\xi)| \leq \frac{1}{B} |\xi \phi''(\xi)| \,.
    \end{align*}
    \item We have
    \begin{align*}
        \left\| \frac{\phi'''}{\phi'} \right\|_{L^1(\{|\phi'(\xi)| \geq R\})} + \left\| \frac{\phi''}{\phi'} \right\|_{(L^2 \cap L^\infty)(\{|\phi'(\xi)| \geq R\})} < \infty \,.
    \end{align*}
    \item There exists $N \in \N$ such that $|(\phi')^{-1}(\{\eta\}) \cap \{|\phi'(\xi)| \geq R\}| \leq N$ for all $\eta \in \R$. 
    \item There exist $E, r > 0$ such that for all $R' > R$ we have
    \begin{align*}
        \frac12 R' \leq |\phi'(\xi)| \leq \frac32 R' \Longrightarrow \sup_{|\zeta| < r} |\phi''(\xi + \zeta)| \leq E R' \,.
    \end{align*}
    \item We have $\lim_{\xi \rightarrow \infty} |\phi'(\xi)| = \infty$.
\end{enumerate}
\begin{corollary}[Generalization of {\cite[Corollary 2.9]{KenigPonceVega1991-2}}] \label{cor:34}
    Let $\phi \in C^3(\R; \R)$ fulfill (J1)-(J6) for some $n > 1$ and let $a \in L^\infty(\R; \R)$.
    Let $b = b(t, \xi) \in L^\infty(\R^2; \C)$ have bounded support in $\xi$, uniformly in $t$.
    For any $T > 0$, $u_0 \in \mc{D}(\R)$ and $s \geq \frac12 \lor \frac{n-1}{4}$ we have
    \begin{align*}
        \|(a(D_x) e^{i t \phi(D_x)} + b(t, D_x)) u_0\|_{L^2_{x \in \R} L^\infty_{t \in [-T,T]}} &\lesssim_{T,s,\phi,a,b} \|u_0\|_{H^s} \,.
    \end{align*}
\end{corollary} 
\begin{proof}
    This was shown for $\phi(\xi) = \xi |\xi|^{n-1}$, $n \geq 2$ in \cite[Corollary 2.9]{KenigPonceVega1991-2} with explicit time growth bound $C(\phi, s, T) = (1 + T)^\rho \ti{C}(s, n-1)$, $\rho > \frac34$. 
    The corollary is a consequence of \cite[Corollary 2.8]{KenigPonceVega1991-2}, which itself is a consequence of \cite[Theorem 2.7]{KenigPonceVega1991-2}. 
    This theorem crucially relies on \cite[Proposition 2.6]{KenigPonceVega1991-2}. 
    We now state and prove generalized versions of these results, noting that the proofs are largely identical. 
    This corollary is then a direct consequence of Theorem \ref{thm:cor2.8}.
\end{proof}
We recall another version of the Van der Corput lemma:
\begin{lemma}[Van der Corput lemma] \label{lem:31}
    Let $\psi \in C_0^\infty(\R; \R)$ and assume that $\Phi \in C^2(\R; \R)$ with $\Phi''\big\vert_{\supp \psi} > \lambda > 0$. Then
    \begin{align} \label{eqn:31}
         \left| \int_{\R} e^{i \Phi(\xi)} \psi(\xi) \dd \xi \right| &\leq 10 \lambda^{-\frac12} (\|\psi\|_{L^\infty} + \|\psi'\|_{L^1}) \,.
    \end{align}
\end{lemma}
\begin{proof}
    See \cite[pp. 309-311]{Stein1986}.
\end{proof}
\begin{lemma}[{Generalization of \cite[Proposition 2.6]{KenigPonceVega1991-2}}] \label{lem:32}
    Let $\phi \in C^3(\R; \R)$ fulfill (J1)-(J6) and let $\psi \in C^\infty(\R; \R)$ with $\supp \psi \subseteq [2^{k-1}, 2^{k+1}]$ for some $k \in \N$. 
    Let $T > 0$. There exists a constant $c = c(\phi) > 0$ for which the function 
    \begin{align}
        \label{eqn:1234} H_k^{n-1}(x) &= \begin{cases}
            2^k &, |x| \leq 2 T R
            \\ 2^{\frac{k}{2}} |x|^{- \frac12} &, 2 T R < |x| \leq c T 2^{(n-1) k}
            \\ (1 + |x|^2)^{-1} &, |x| > c T  2^{(n-1) k}
        \end{cases}
    \end{align}
    fulfills
    \begin{align*}
        \sup_{|t| \leq T} \int_{\R} e^{i (t \phi(\xi) + x \xi)} \psi(\xi) \dd \xi &\lesssim_{T,\phi} H_k^{n-1}(x) (\|\psi''\|_{L^1} + \|\psi'\|_{L^1 \cap L^\infty} + \|\psi\|_{L^\infty}) \,.
    \end{align*}
\end{lemma}
\begin{proof}
    We may assume $t \in [0, T]$, $x \in \R$. If $|x| \leq 2 T R$ we estimate
    \begin{align*}
        \left| \int_{\R} e^{i (t \phi(\xi) + x \xi)} \psi(\xi) \dd \xi \right|
        &\leq 2 (2^{k+1} - 2^{k-1}) \|\psi\|_{L^\infty} = 2^k 3 \|\psi\|_{L^\infty} \,.
    \end{align*}
    We therefore assume $|x| > 2 T R > 1$ and define
    \begin{align*}
        \Omega &= \left\{ \xi \in [2^{k-1}, 2^{k+1}]: |t \phi'(\xi) + x| \leq \frac{|x|}{2} \right\}
        \\ \ti{\Omega} &= \left\{ \xi \in [2^{k-1}, 2^{k+1}]: |t \phi'(\xi) + x| \leq \frac{|x|}{3} \right\} \,.
    \end{align*}
    Note that $\xi \in \Omega$ implies
    \begin{align*}
        \frac12 \frac{|x|}{t} &\leq |\phi'(\xi)| \leq \frac32 \frac{|x|}{t} \,,
    \end{align*}
    and furthermore that $|x| > 2 T R$ and $t \leq T$ ensure $\frac12 \frac{|x|}{t} \geq R$.
    Since $R > 1$ can be chosen arbitrarily large, we are in a regime where $|\phi'(\xi)|$ is large and hence (J1)-(J6) may be applied. 
    Let $\eta \in C^\infty(\R; \R)$ with $\supp \eta \subseteq \Omega$ and $\eta\big\vert_{\ti{\Omega}} = 1$. 
    We know that $\|\eta\|_{L^\infty} \leq 1$ and now prove furthermore that $\|\eta'\|_{L^\infty} \lesssim_\phi 1$. 
    We do this by showing $\ti{\Omega} + B_r(0) \subseteq \Omega$ for some small radius $r > 0$ independent of $t$. Let $0 < |\zeta| < r$ and $\xi \in \ti{\Omega}$. Then
    \begin{align*}
        |t \phi'(\xi + \zeta) + x| &\leq \frac{|x|}{3} + t r \sup_{|a| < r} |\phi''(\xi + a)| \,.
    \end{align*}
    With (J5) we obtain
    \begin{align*}
        t r \sup_{|a| < r} |\phi''(\xi + a)| &\leq t r E \frac{|x|}{t} = E r |x| \,,
    \end{align*}
    so a choice of $r$ with $E r < \frac16$ suffices. Since (J4) implies that $\Omega$ and $\ti{\Omega}$ can each written as unions of $2 N + 1$ or fewer closed intervals, 
    there exists a choice of $\eta$ for which also $\|\eta'\|_{L^1}, \|\eta''\|_{L^1} \lesssim_\phi 1$.
    Now if $\xi \in \supp \eta \subseteq \Omega$, then (J2) implies
    \begin{align*}
        |(t \phi(\xi) + x \xi)''| &= t |\phi''(\xi)| \geq t B \left| \frac{\phi'(\xi)}{\xi} \right| \geq \frac{B}{2^{k+1}} \frac{|x|}{2} = \frac{B}{4} 2^{-k} |x| \,.
    \end{align*}
    We apply \eqref{eqn:31} to obtain
    \begin{align*}
        \left| \int_{\R} e^{i (t \phi(\xi) + x \xi)} (\eta \psi)(\xi) \dd \xi \right| &\lesssim_\phi 2^{\frac{k}{2}} |x|^{-\frac12} (\|(\eta \psi)'\|_{L^1} + \|\eta \psi\|_{L^\infty}) \lesssim_\phi 2^{\frac{k}{2}} |x|^{-\frac12} (\|\psi'\|_{L^1} + \|\psi\|_{L^\infty}) \,.
    \end{align*}
    On the other hand if $\xi \in \supp (1 - \eta)$, then
    \begin{align*}
        |(t \phi(\xi) + x \xi)'| = |t \phi'(\xi) + x| \geq \frac{|x|}{3} \,.
    \end{align*}
    Due to (J6) we can perform the integration by parts below without any boundary terms appearing. We estimate
    \begin{align*}
        & \left| \int_{\R} e^{i (t \phi(\xi) + x \xi)} ((1 - \eta) \psi)(\xi) \dd \xi \right|
        \\ &= \left| \int_{\R} e^{i (t \phi(\xi) + x \xi)} \del_\xi \left( \frac{1}{t \phi'(\xi) + x} \del_\xi \left( \frac{1}{t \phi'(\xi) + x} ((1 - \eta) \psi)(\xi) \right) \right) \dd \xi \right|
        \\ &\leq \int_{\supp (1 - \eta)} \left| \frac{((1 - \eta) \psi)'' (t \phi' + x) - t \phi''' (1 - \eta) \psi}{(t \phi' + x)^3} \right| \dd \xi
        \\ &+ \int_{\supp (1 - \eta)} \left| \frac{3 t \phi'' \big(((1 - \eta) \psi)' (t \phi' + x) - t \phi'' (1 - \eta) \psi\big)}{(t \phi' + x)^4} \right| \dd \xi
        \\ &\lesssim \frac{1}{|x|^2} \left( 1 + \left\| \frac{t \phi'''}{t \phi' + x} \right\|_{L^1(\R \setminus \ti{\Omega})} + \left\| \frac{t \phi''}{t \phi' + x} \right\|_{L^\infty(\R \setminus \ti{\Omega})} + \left\| \frac{(t \phi'')^2}{(t \phi' + x)^2} \right\|_{L^1(\R \setminus \ti{\Omega})} \right)
        \\ & (\|((1 - \eta) \psi)''\|_{L^1} + \|((1 - \eta) \psi)'\|_{L^1} + \|(1 - \eta) \psi\|_{L^\infty})
        \\ &\lesssim_\phi \frac{1}{|x|^2} \Bigg( \left( 1 + \left\| \frac{\phi'''}{\phi'} \right\|_{L^1(\{|\phi'(\xi)| \geq R\})} + \left\| \frac{\phi''}{\phi'} \right\|_{(L^2 \cap L^\infty)(\{|\phi'(\xi)| \geq R\})}^2 \right) \left\| \left\langle 1 - \frac{x}{t \phi' + x} \right\rangle^2 \right\|_{L^\infty(\R \setminus \ti{\Omega})}
        \\ &+ \left\| \phi''' \right\|_{L^1(\{|\phi'(\xi)| \leq R\})} + \left\| \phi'' \right\|_{(L^2 \cap L^\infty)(\{|\phi'(\xi)| \leq R\})}^2 \Bigg) (\|\psi''\|_{L^1} + \|\psi'\|_{L^1 \cap L^\infty} + \|\psi\|_{L^\infty}) \,.
    \end{align*}
    We apply (J3) to estimate the norms with $|\phi'(\xi)| \geq R$, and (J6) together with continuity of $\phi'''$ to estimate the norms with $|\phi'(\xi)| \leq R$.
    We have shown that
    \begin{align*}
        \left| \int_{\R} e^{i (t \phi(\xi) + x \xi)} ((1 - \eta) \psi)(\xi) \dd \xi \right| &\lesssim_\phi |x|^{-2} (\|\psi''\|_{L^1} + \|\psi'\|_{L^1 \cap L^\infty} + \|\psi\|_{L^\infty}) \,.
    \end{align*}
    In summary,
    \begin{align*}
        \left| \int_{\R} e^{i (t \phi(\xi) + x \xi)} \psi(\xi) \dd \xi \right| 
        &\lesssim_\phi 2^{\frac{k}{2}} |x|^{-\frac12} (\|\psi'\|_{L^1} + \|\psi\|_{L^\infty}) + |x|^{-2} (\|\psi''\|_{L^1} + \|\psi'\|_{L^1} + \|\psi\|_{L^\infty})
        \\ &\lesssim_{T,\phi} 2^{\frac{k}{2}} |x|^{-\frac12} (\|\psi''\|_{L^1} + \|\psi'\|_{L^1 \cap L^\infty} + \|\psi\|_{L^\infty}) \,.
    \end{align*}
    Note that
    \begin{align*}
        |t \phi'(\xi) + x| \leq \frac12 |x| &\Longrightarrow \frac{t |\phi'(\xi)|}{|x|} \geq \frac12 \xRightarrow{\text{(J1)}} |x| \leq 2 T A 2^{(n-1) (k + 1)} \,,
    \end{align*}
    so we can choose a constant $c(\phi)$ so that in the case $|x| > c(\phi) T 2^{(n-1) k}$ only the term with decay $|x|^{-2}$ appears. Since in thise case $|x| > c T$, we have $|x|^{-2} \lesssim_{T,\phi} (1 + |x|^2)^{-1}$.
\end{proof}
\begin{theorem}[{Generalization of \cite[Theorem 2.7, Corollary 2.8]{KenigPonceVega1991-2}}] \label{thm:cor2.8}
    Let $\phi \in C^3(\R; \R)$ fulfill (J1)-(J6) for some $n > 1$ and let $a \in C^\infty(\R; \R)$.
    Let $b = b(t, \xi) \in L^\infty(\R^2; \C)$ have bounded support in $\xi$, uniformly in $t$.
    For any $T > 0$, $u_0 \in \mc{D}(\R)$ and $s \geq \frac12 \lor \frac{n-1}{4}$ we have
    \begin{align*}
        \left( \sum_{j \in \Z} \sup_{|t| \leq T} \sup_{j \leq x \leq j+1} |(a(D_x) e^{i t \phi(D_x)} + b(t, D_x)) u_0(x)|^2 \right)^{\frac12} &\lesssim_{T,s,\phi,a,b} \|u_0\|_{H^{\frac12 \lor \frac{n-1}{4}}} \,.
    \end{align*}
\end{theorem}
\begin{proof}
    The proof is identical to that of \cite[Theorem 2.7]{KenigPonceVega1991-2}, except \cite[Proposition 2.6]{KenigPonceVega1991-2} is replaced by our Lemma \ref{lem:32}, and the time interval $|t| \leq 1$ is replaced with $|t| \leq T$.
    The extension of of the time interval is performed separately in \cite[Corollay 2.8]{KenigPonceVega1991-2} using a scaling argument, which yields an explicit algebraic growth bound $(1 + T)^{\rho}, \rho > \frac34$ for the constant in the time $T$. 
    Since $\phi(\xi)$ is not necessarily homogeneous in our situation, this scaling argument can not be trivially generalized. 
    As we do not need to know a growth bound for the constant in $T$, we can simply skip this argument and work directly with $|t| \leq T$.
    The additional operators $a(D_x)$ and $b(t, D_x)$ are not present in the reference, but in the proof the claim is reduced to estimating an $L^2_x(\R)$-norm of a linear function of $u$, and at that point the weight $a(\xi)$ disappears.
    Similarly, $b(t, D_x)$ disappears because an inhomogeneous dyadic decomposition using frequency projectors $\psi_k(D_x)$, $k \in \N$ is used, which for the term with $b(t, D_x)$ leaves only the case $k = 0$. 
    In this case the oscillation from the semigroup is not necessary, i.e. Lemma \ref{lem:32} is not used.

\end{proof}

\printbibliography

\begin{minipage}[t]{0.48\textwidth}
\raggedright\small
\textbf{Xian Liao}\\
School of Mathematical Sciences\\
Dalian University of Technology\\
Linggong Road 2 \\
116024 Dalian, P.R.China\\
\href{mailto:liao@dlut.edu.cn}{liao@dlut.edu.cn}
\end{minipage}\hfill
\begin{minipage}[t]{0.48\textwidth}
\raggedright\small
\textbf{Robert Wegner}\\
Institute for Analysis\\
Karlsruhe Institute of Technology\\
Englerstraße 2\\
76131 Karlsruhe, Germany\\
\href{mailto:robert.wegner@kit.edu}{robert.wegner@kit.edu}
\end{minipage}

\end{document}